%% file: main.tex
\documentclass[reqno,oneside]{amsart} %oneside for margins

\usepackage{amsthm}

\usepackage{%srcltx,
amssymb,amsmath,enumerate,stackrel,
mathrsfs,verbatim,stmaryrd,
leftidx}
\allowdisplaybreaks[2]
\SetSymbolFont{stmry}{bold}{U}{stmry}{m}{n}
\usepackage{xcolor} % für Hintergrundfarben
\usepackage{xr-hyper}
\usepackage{nicefrac,stackrel}
\usepackage[all,cmtip]{xy}
\usepackage[utf8]{inputenc} % for BibTeX from arXiv
\chardef\cprime"7E % ditto
\usepackage{cases}

\usepackage[T1]{fontenc} % to get \k{e} in Kestutis (in bib file) right...

\usepackage[hyperindex,linktocpage=true]{hyperref}
\hypersetup{
    colorlinks,
    linkcolor={red!50!black},
    citecolor={red!50!black}, %{blue!50!black},
    urlcolor={red!50!black}, % {blue!80!black},
    filecolor={red!50!black} % black
}

\usepackage{tikz-cd}%Für Bildchen
\usetikzlibrary{matrix,arrows}

\usepackage{xcolor}

\usepackage{theoremref}

\definecolor{labelkey}{rgb}{1,0,0}

\makeatletter
\@addtoreset{equation}{section}
\makeatother

\numberwithin{equation}{section}

\theoremstyle{definition}
\newtheorem{Def}{Definition}[section] \newcommand{\defi}{\begin{Defi}} \newcommand{\xdefi}{\end{Defi}} %\crefname{Defi}{Def.}{Defs.}
\newtheorem{DefiLemm}[Def]{Definition and Lemma} \newcommand{\defilemm}{\begin{DefiLemm}} \newcommand{\xdefilemm}{\end{DefiLemm}} %\crefname{DefiLemm}{Definition~and~Lem.}{Defs.~and~Lems.}
\newtheorem{Exp}[Def]{Example} \newcommand{\exam}{\begin{Exp}} \newcommand{\xexam}{\end{Exp}} %\crefname{Bsp}{Exam.}{Exam.}
\newtheorem{Syno}[Def]{Synopsis} \newcommand{\syno}{\begin{Syno}} \newcommand{\xsyno}{\end{Syno}} %\crefname{Syno}{Syno.}{Syno.}
\newtheorem{Rem}[Def]{Remark} \newcommand{\rema}{\begin{Bem}} \newcommand{\xrema}{\end{Bem}} %\crefname{Bem}{Rem.}{Rems.}
\newtheorem{Notation}[Def]{Notation} \newcommand{\nota}{\begin{Notation}} \newcommand{\xnota}{\end{Notation}} %\crefname{Bem}{Rem.}{Rems.}
\newtheorem{Warning}[Def]{Warning} \newcommand{\warn}{\begin{Warning}} \newcommand{\xwarn}{\end{Warning}} %\crefname{Bem}{Rem.}{Rems.}
\newtheorem{Situation}[Def]{Situation} \newcommand{\situ}{\begin{Situation}} \newcommand{\xsitu}{\end{Situation}} %\crefname{Bem}{Rem.}{Rems.}

\theoremstyle{plain}
\newtheorem{Theo}[Def]{Theorem} \newcommand{\theo}{\begin{Theo}} \newcommand{\xtheo}{\end{Theo}} %\crefname{Theo}{Thm.}{Thms.}
\newtheorem{Prop}[Def]{Proposition} \newcommand{\prop}{\begin{Prop}} \newcommand{\xprop}{\end{Prop}} %\crefname{Satz}{Prop.}{Props.}
\newtheorem{Lem}[Def]{Lemma} \newcommand{\lem}{\begin{Lem}} \newcommand{\xlemm}{\end{Lem}} %\crefname{Lemm}{Lem.}{Lems.}
\newtheorem{Cor}[Def]{Corollary} \newcommand{\coro}{\begin{Cor}} \newcommand{\xcoro}{\end{Cor}}
\newtheorem{Ques}[Def]{Question} \newcommand{\ques}{\begin{Ques}} \newcommand{\xques}{\end{Ques}}
\newtheorem{Conj}[Def]{Conjecture} \newcommand{\conj}{\begin{Conj}} \newcommand{\xconj}{\end{Conj}}

\newcommand{\eqn}{\begin{equation}} \newcommand{\xeqn}{\end{equation}}
\newcommand{\eqnarr}{\begin{eqnarray*}} \newcommand{\xeqnarr}{\end{eqnarray*}}
\newcommand{\eqnarra}{\begin{eqnarray}} \newcommand{\xeqnarra}{\end{eqnarray}}

\newcommand{\pf}{\begin{proof}} \newcommand{\xpf}{\end{proof}}

\usepackage[top=1.5cm, bottom=1.5cm, left=2cm, right=2cm, heightrounded, marginparwidth=3.5cm, marginparsep=.2cm]{geometry}

\usepackage{color}

\newcommand{\nc}{\newcommand}%Faulheit
\nc{\StP}[1]{\cite[Tag~\href{http://stacks.math.columbia.edu/tag/#1}{#1}]{StacksProject}}
\nc{\StPd}[2]{\cite[Tags~\href{http://stacks.math.columbia.edu/tag/#1}{#1}, \href{http://stacks.math.columbia.edu/tag/#2}{#2}]{StacksProject}} % Stacks-Project

\renewcommand{\l}{\left(}
\renewcommand{\r}{\right)}
\nc{\lb}{\left[}
\nc{\rb}{\right]}
\nc{\lsb}{\begin{bmatrix}}
\nc{\rsb}{\end{bmatrix}}
\def\ul{\ensuremath{\underline}}

\def\mf{\ensuremath\mathfrak}
\def\mc{\ensuremath\mathcal}
\def\mb{\ensuremath\mathbb}
\def\f{\ensuremath\textbf}

\nc{\category}[1]{\mathrm{#1}}
\nc{\op}[1]{\operatorname{#1}}

\nc{\Set}{\category{Set}}
\nc{\Map}{\category{Map}}
\nc{\Dat}{\category{Dat}}
\nc{\GlobDat}{\category{GlobDat}}
\nc{\Seq}{\category{Seq}}
\nc{\rk}{\category{rk}}
\nc{\size}{\category{size}}
\nc{\type}{\category{type}}
\nc{\ParDat}{\category{ParDat}}
\nc{\TwParDat}{\category{TwParDat}}
\nc{\CompDat}{\category{CompDat}}
\nc{\TwCompDat}{\category{TwCompDat}}
\nc{\im}{\category{im}}
\nc{\Loc}{\category{Loc}}
\nc{\MC}{\category{MC}}
\nc{\Rep}{\category{Rep}}

\nc{\UL}{\category{UL}}
\nc{\ULD}{\category{ULD}}
\nc{\ParMod}{\category{ParMod}}
\nc{\ParCoh}{\category{ParCoh}}
\nc{\ParQCoh}{\category{ParQCoh}}
\nc{\ParBun}{\category{ParBun}}
\nc{\QCoh}{\category{QCoh}}
\nc{\SUM}{\category{SUM}}
\nc{\Inh}{\category{Inh}}
\nc{\secmul}{\category{multseq}}

\nc{\sub}{\category{sub}}
\nc{\Res}{\category{Res}}

\nc{\fin}{\mathrm{fin}}
\nc{\oro}{\category{or}}
\nc{\p}{^{\prime}}
\renewcommand{\b}{^{\bullet}}
\nc{\wreg}{\category{wreg}}
\nc{\reg}{\category{reg}}
\nc{\Reg}{\category{Reg}}
\nc{\degen}{\category{degen}}
\nc{\ext}{\category{ext}}
\nc{\irec}{\category{irec}}
\nc{\glob}{\category{glob}}
\nc{\para}{\category{par}}
\nc{\tw}{\category{tw}}
\nc{\inde}{\category{idcp}}
\nc{\geominde}{\category{geom-idcp}}
\nc{\mini}{\category{min}}
\nc{\irrel}{\category{irrel}}
\nc{\noncusp}{\category{noncusp}}
\nc{\inter}{\category{inter}}
\nc{\internull}{{\underline{0},\category{inter}}}
\nc{\interone}{{1,\category{inter}}}
\nc{\minnull}{{0,\category{min}}}
\nc{\gennull}{{0,\category{gen}}}
\nc{\Rnull}{{0,\category{R}}}
\nc{\dec}{\category{dec}}
\nc{\spc}{\category{spc}}
\nc{\gen}{\category{gen}}
\nc{\qs}{\category{qs}}
\nc{\et}{\category{\acute{e}t}}
\nc{\liset}{\category{lis-\acute{e}t}}
\nc{\red}{\category{red}}
\nc{\interc}{{\category{inter}\subset}}
\nc{\interonec}{{1,\category{inter}\subset}}
\nc{\internullec}{{0,\category{inter}\subset}}
\nc{\onegentc}{{1,\l\category{gen}+\frac{t-1}{t}\r\subset}}
\nc{\new}{\category{new}}
\nc{\comp}{\category{comp}}
\nc{\sem}{\category{ss}}
\nc{\uni}{\category{u}}
\nc{\conjug}{\category{conj}}
\nc{\irr}{\category{irr}}
\nc{\ana}{\category{an}}
\nc{\loc}{\category{loc}}

\def\FL{\ensuremath\text{Fl}}
\nc{\supp}{\category{supp}}
\nc{\Bun}{\op{\boldsymbol{\mathcal{B}un}}}
\nc{\PIC}{\op{\boldsymbol{\mathcal{P}ic}}}
\nc{\FlagBun}{\op{\boldsymbol{\mathcal{F}lag\mathcal{B}un}}}
\nc{\AlgStk}{\category{AlgStk}}
\nc{\Spec}{\category{Spec}}
\nc{\Pic}{\category{Pic}}
\nc{\mindim}{\category{mindim}}
\nc{\Stab}{\category{Stab}}
\nc{\dom}{\category{dom}}
\nc{\pr}{\category{pr}}
\nc{\Flat}{\category{Flat}}
\nc{\Push}{\category{Push}}
\nc{\rigid}{ \hspace{-.1cm} \fatslash \hspace{.10cm}}

\nc{\GL}{\mathrm{GL}}
\nc{\PGL}{\mathrm{PGL}}
\nc{\coker}{\mathrm{coker}}
\nc{\Mod}{\mathrm{Mod}}
\nc{\Ext}{\mathrm{Ext}}
\nc{\Tor}{\mathrm{Tor}}
\nc{\Aut}{\mathrm{Aut}}
\nc{\End}{\mathrm{End}}
\nc{\Hom}{\mathrm{Hom}}
\nc{\id}{\mathrm{id}}
\nc{\act}{\mathrm{act}}

\begin{document}

\author{Paul Siemon}

\title{Stacks of Parabolic Vector Bundles \\ in the Geometric Langlands Program}

\address{Technical University of Darmstadt, Department of Mathematics, 64289 Darmstadt, Germany}
\email{%
  \href{mailto:siemon@mathematik.tu-darmstadt.de}{%
    \mbox{siemon@mathematik.tu-darmstadt.de}%
  }%
}

\begin{abstract}
I define a locus inside stacks of parabolic vector bundles on smooth projective curves which is important for the construction of twisted Hecke eigensheaves in the ramified geometric Langlands program for $\GL_n$. I give an explicit description of this locus in the case of rank $3$ parabolic vector bundles on $\mb P^1$ with parabolic structure at $4$ points. This example involves partial flags in the level structure as well as a non-trivial twisted equivariance condition.
\end{abstract}

\maketitle

\tableofcontents

%%%%%%%%%%Introduction%%%%%%%%%%%%%%%%%%%

\input{1_Introduction}
\input{2_Local_and_Global_Data}

\input{3_Stacks_of_Parabolic_Bundles}
\input{4_Interesting_Locus_Classification}

\input{5_Geometry}

%\newpage
%\input{6_To_Do}

% Inhalt der Datei references.tex einfügen

\input{references.tex}\end{document}

%% file: 1_Introduction.tex
\section{Introduction}
In this article we study stacks of parabolic vector bundles on a connected smooth projective curve $X$ over an algebraically closed field $k$. Such stacks arise naturally in the Global Geometric Langlands Program for $\GL_n$. The contents of this article can be subdivided in two parts:

First, we define several loci inside stacks of parabolic vector bundles, some of which are naturally defined for any algebraic stack, whereas others play a prominent role in the construction of Hecke eigensheaves in the geometric Langlands program for $\GL_n$. We wish to understand their precise relation. 

Second, we study a specific example where $X = \mb P^1$ and a certain parabolic structure at four points is present. We will give explicit descriptions of the aforementioned loci in this example. We plan to use this description in forthcoming work, to construct Hecke eigensheaves in this specific example.

\subsection{Background} \label{Motivational part of intro}
The unramified geometric Langlands program for split reductive groups over $\mb C$ is concerned with the following problem: Let $X$ be a connected smooth projective algebraic curve over $\mb C$ and $G$ a split connected reductive $\mb C$-group. Let $^LG$ be the Langlands dual group of $G$ over $\mb C$. For an irreducible $^LG$-local system $\mc L$ on $X$ find a \textit{Hecke eigensheaf} $\mc H_{\mc L}$ with eigenvalue $\mc L$ on the stack $\Bun_G$ of $G$-bundles on $X$. \\
More precisely, there is a collection $\{\mb H_V\}_{V \in \Rep(^LG)}$ of functors, called \textit{Hecke operators}, from the category of sheaves on $\Bun_G$ to the category of sheaves on $X \times \Bun_G$. This collection is indexed by the category $\Rep(^LG)$ of finite-dimensional algebraic representations of $^LG$. The datum of an eigensheaf for $\mc L$ consists of a sheaf $\mc H_{\mc L}$ on $\Bun_G$ together with compatible isomorphisms
\begin{align*}
    \mb H_V(\mc H_{\mc L}) \simeq V_{\mc L} \boxtimes \mc H_{\mc L},
\end{align*}
which are compatible with the tensor structure on $\Rep(^LG)$. The expression $V_{\mc L}$ denotes the twist of $V$ by $\mc L$.

Note that the term $^LG$\textit{-local system} can mean at least two different objects: A principal $^LG$-bundle with connection or a representation of the fundamental group of $X^\ana$ with values in $^LG(\mb C)$. We refer to the former case as the \textit{de Rham setting} and to the latter as the \textit{Betti setting}. The term \textit{sheaf} has to be understood as \textit{object in a derived category of sheaves} and which category to take depends on the context. For an introduction to the topic see e.g. \cite{GaitsgoryInformal} or \cite{FrenkelCFT}.

In both settings, de Rham and Betti, the problem of finding Hecke eigensheaves has recently been solved by D. Arinkin, D. Beraldo, J. Campbell, L. Chen, J. Faergeman, D. Gaitsgory, K. Lin, S. Raskin and N. Rozenblyum in a series of five papers \cite{GLC1},\cite{GLC2},\cite{GLC3},\cite{GLC4},\cite{GLC5}. Even better, they prove in the de Rham resp. Betti setting an equivalence of categories between the category of $D$-modules resp. ind-coherent $\mb C$-sheaves with nilpotent singular support on $\Bun_G$ and the category of ind-coherent sheaves with nilpotent singular support on the derived stack of de Rham resp. Betti $^LG$-local systems. The eigensheaf for a (de Rham resp. Betti) $^LG$-local system $\mc L$ can be obtained by taking the preimage under this equivalence of the skyscraper sheaf at $\mc L$.

In contrast to the unramified picture described above, in this article we concentrate on the \textit{ramified situation}. On the spectral side allowing ramification means to allow local systems which are defined on non-empty Zariski-open subsets of $X$, i.e. on $X \setminus S$, where $S$ is a finite set of closed points of $X$. With such a local system $\mc L$ on $X \setminus S$ one can associate for each point $s$ in $S$ an isomorphism class of local systems on the punctured local disc around $s$ which is called the \textit{local monodromy} of $\mc L$ at $s$.

In studying the global geometric Langlands program with ramification, it is convenient to fix the set $S$ as well as, for the spectral category, specified local monodromies $\{\mc L_s\}_{s \in S}$. This means that on the spectral side, we are interested in the collection of $^LG$-local systems on $X \setminus S$ with monodromy at $s$ isomorphic to $\mc L_s$.

Let us for the sake of concreteness restrict our setting: From now on we consider the case $G = \GL_n$ and the Betti setting. Then $^LG = \GL_n$ and the datum of $\mc L_s$ is an equivalence class of rank $n$ local systems of $\mb C$-vector spaces on a punctured complex analytic disc around $s$ inside $X^\ana$. Given a collection of equivalence classes $\{\mc L_s\}_{s \in S}$ we denote by $\Loc_\irr(X \setminus S,\{\mc L_s\})$ the category of irreducible local systems on $X \setminus S$ which have monodromy $\mc L_s$ at $s$. The datum of $\mc L_s$ is equivalent to a conjugacy class in $\GL_n(\mb C)$, which is uniquely described by its Jordan canonical form.

On the automorphic side, we note that $\Bun_{\GL_n}$ is naturally identified with the stack of rank $n$ vector bundles on $X$. When passing to the above ramified situation, the stack of rank $n$ vector bundles has to be replaced by a stack of rank $n$ parabolic vector bundles on $X$. See e.g. \cite{HeinlothThesis} for an instance of this idea, which goes back to Drinfeld \cite{DrinfeldGL2}. A \textit{rank $n$ parabolic vector bundle on}\footnote{Some authors use the term \textit{quasi-parabolic vector bundle} for what is called a \textit{parabolic vector bundle} in this article and reserve the latter for a vector bundle equipped with weighted flags in some fibers.} $X$ is a rank $n$ vector bundle $\mc E$ on $X$ plus the datum of a flag $F_s$ in the fiber $\mc E \otimes k(s)$ for each point $s \in S$, cf. Definition \ref{Definition Parabolic Sheaves}. The flag $F_s$ is supposed to have a certain type, which is specified by $\mc L_s$ as follows: Write a representative of $\mc L_s$ in Jordan canonical form as
\begin{align*}
    \mc L_s \widehat{=} 
\begin{pmatrix}
J_{\lambda_1} &        &        & 0 \\
                    & \ddots &        &   \\
                    &        & \ddots &   \\
0                   &        &        & J_{\lambda_m}
\end{pmatrix},
\end{align*}
where $\lambda_1, \ldots, \lambda_m$ are the pairwise distinct eigenvalues of $\mc L_s$, each $J_{\lambda_i}$ is in Jordan canonical form and has a single eigenvalue $\lambda_i$. Each $J_{\lambda_i}$ determines an unordered partition $p_i$ of the size $\nu_i$ of $J_{\lambda_i}$. Then the type of flag\footnote{To the knowledge of the author it is not stated anywhere in the literature that this is the correct recipe in general to obtain the flag type from the local monodromies. But see \cite{YunEpipelagic} for a situation, where (a more general form of) this principle works.} at $s$ is $(p_1^\conjug, \ldots, p_m^\conjug)$, where $p_i^\conjug$ denotes the conjugate partition of $p_i$.

Denote by $D$ the collection $\{D_s\}_{s \in S}$ of types of flags specified by $\{\mc L_s\}_{s \in S}$ as above and let $\Bun_D$ be the stack of parabolic rank $n$ vector bundles of parabolic type $D_s$ at $s \in S$. If all $m = 1$ and $\lambda_1 = 1$, then objects in $\Loc_\irr(X \setminus S, \{\mc L_s\})$ have unipotent monodromy at all points in $S$ and Hecke eigensheaves of objects in \\$\Loc_\irr(X \setminus S, \{\mc L_s\})$ should be objects in the derived category of sheaves on $\Bun_D$. 

If non-trivial eigenvalues occur, one has to consider the natural map $\pi_D^+ : \Bun_D^+ \to \Bun_D$, where $\Bun_D^+$ parametrizes parabolic bundles of type $D$ plus a basis of the subquotients. Writing $D = \{D_s\}_{s \in S}$ as above, the map $\pi_D^+$ is a torsor under the group $L := \prod_{s \in S} L_s$, where $L_s \subset \GL_n$ is the standard Levi subgroup of type $D_s$. Then Hecke eigensheaves for objects in $\Loc_\irr(X \setminus S, \{\mc L_s\})$ should be objects in the derived category of $\chi_\Lambda$-equivariant sheaves on $\Bun_D^+$. Here, $\chi_\Lambda$ is a character sheaf on $L$ associated with the collection $\Lambda$ of all eigenvalues of all $\mc L_s$, $s \in S$. See for example \cite[Section~4.5]{BenZviNadler} for an appearance of this idea in the Betti setting.

\subsection{Contents and Motivation}
In the first part of this article (sections \ref{section data} and \ref{stacks of parabolic vector bundles}) we study stacks of parabolic vector bundles on a smooth projective curve over an algebraically closed field. In section \ref{section data} the notation for various sets of and operations on local and global data is fixed. This subsection may be skipped on a first reading. In section \ref{stacks of parabolic vector bundles} we define stacks of parabolic vector bundles (Definition \ref{Definition of stack of parabolic vector bundles depending on D}), define the \textit{interesting locus} inside these stacks (Definition \ref{Definition of interesting locus}), which are important for the construction of Hecke eigensheaves, and study the interplay of this locus with certain Hecke operators (Theorem \ref{Proposition interesting locus and hecke operators}).

In the second part of this article (section \ref{section on classification} and \ref{Last section}), we study a specific example of rank $3$ on $X = \mb P^1$, where ramification is present at four distinct points $0,1,t,\infty \in \mb P^1$. On the spectral side, the local monodromy classes $\mc L_s$ are given for three points by
\begin{align*}
    \begin{pmatrix}
        \lambda & & \\
        & 1 & \\
        & & 1
    \end{pmatrix}
\end{align*}
and for the fourth point by
\begin{align*}
    \begin{pmatrix}
        \lambda^{-1} & 1&  \\
        & \lambda^{-1} & 1 \\
        & & \lambda^{-1}
    \end{pmatrix},
\end{align*}
where $\lambda \in \mb C^{\times} \setminus \{\pm 1\}$. In this example, we explicitly calculate the loci introduced in the first part of the article.

The choice of example is motivated as follows: In constructing eigensheaves for irreducible local systems with specified local monodromy, a good measure for the complexity of the geometry of the situation at hand is the dimension of the coarse moduli space of such objects. For given data, this space may be empty. The zero-dimensional case is the rigid case and has been considered by several authors (\cite{YunRigidity}, \cite{JakobKamgarpourYi}, \cite{HeinlothNgoYun}).

We are interested in the next step: The two-dimensional case (the coarse moduli of irreducible local systems carries a symplectic structure, hence its dimension must be even). One example of such a situation, studied by uit de Bos in his thesis \cite{uitdeBosThesis}, is the case of $\mb P^1$ and, on the spectral side, irreducible rank $2$ local systems with regular unipotent monodromy at four points. One may ask how to produce further examples. 

One way to do this is to use Katz' middle convolution algorithm \cite{KatzRigid}: It uses the structure of the additive group $\mb G_a \subset \mb P^1$ to produce for any finite set $S \subset \mb P^1$ with $\infty \in S$ and $\#S \ge 3$ and any $\lambda \in \mb C^{\times} \setminus \{1\}$ a functor $\MC_\lambda$ from the category of tame irreducible local systems\footnote{Here we mean local systems on $(\mb P^1)^\ana$.} of rank $\ge 2$ on $\mb P^1 \setminus S$ to itself. Further, $\MC_\lambda$ induces a functor $\MC_\lambda^\loc$ on local monodromy data. Our example arises by applying $\MC_\lambda^\loc$ to the rank $2$ situation with four points and unipotent monodromy.

One prominent feature of this example is that the geometry of the stack $\Bun_D$ can be approached computationally in a very explicit way since the geometry of $\mb P^1$ is not too complicated. This is an advantage of the ramified situation: One may consider non-trivial examples with $X = \mb P^1$ (since $\mb P^1$ is simply connected, the unramified situation is trivial in this case). In section \ref{section on classification} we explicitly classify the points in $\vert \Bun_D \vert^\inter$ by writing them down in coordinates. In section \ref{Last section} we produce explicit charts of various loci in the stack $\Bun_D$ and compute their Picard groups.

The present article is heavily influenced by the results obtained in the ramified setting by uit de Bos \cite{uitdeBosThesis}, Arinkin \cite{ArinkinOrthogonality}, Yun \cite{YunRigidity} and Heinloth \cite{HeinlothThesis}. We plan to use the results in the present paper to construct Hecke eigensheaves in the explicit rank $3$ example described above. To the knowledge of the author, no non-rigid example of ramified geometric Langlands involving partial flags as well as a non-trivial and non-generic equivariance condition has yet been studied in the literature.

\subsection{Main results I: Geometry of $\Bun_D$}
The main definition in the first part of this article (sections \ref{section data} and \ref{stacks of parabolic vector bundles}) is the \textit{interesting locus} $\vert \Bun_D \vert^\inter$, a subspace of the topological space $\vert \Bun_D \vert$ of the stack $\Bun_D$. It is related to the support of \textit{twisted equivariant cuspidal} eigensheaves on $\Bun_D^+$. Such a sheaf corresponds under the Langlands correpsondence to an \textit{irreducible} local system $\mc L$. The adjective \textit{cuspidal} means that the sheaf on $\Bun_D^+$ does not come via geometric parabolic induction from a stack corresponding to a Levi subgroup, cf. section 6 in \cite{uitdeBosThesis} for a formulation of this condition. It restricts the support of cuspidal sheaves: There are points $P \in \vert \Bun_D \vert$ such that the stalk at any point $P^+ \in \vert \Bun_D^+ \vert$ over $P$ of every cuspidal sheaf on $\Bun_D$ vanish. We define a class of \textit{noncuspidal} points $P$ with this property (Definition \ref{Main Definition noncupsidality recognizing}).

Further, the equivariance condition puts a similar constraint on the support of \textit{any} twisted equivariant sheaf on $\Bun_D^+$, see \cite[Lemma~2.7.9]{YunRigidity}. As done in \cite{YunRigidity}, we define a class of \textit{irrelevant points} $P \in \vert \Bun_D \vert$ such that the stalk at each point $P^+ \in \vert \Bun_D^+ \vert$ over $P$ of every $\chi_\Lambda$-twisted equivariant sheaf on $\Bun_D$ vanishes (Definition \ref{Main Definition of irrelevance recognizing}).

Combining these two conditions we arrive at the following main definition:

\begin{Def}[Definition \ref{Definition of stack of parabolic vector bundles depending on D}]
    The \textit{interesting locus} $\vert \Bun_D \vert^\inter$ is the set of all those points $P \in \vert \Bun_D \vert$ which are neither irrelevant nor noncuspidal.
\end{Def}

Note that this locus depends not only on the geometric situation $\pi_D^+ : \Bun_D^+ \to \Bun_D$, but also on the collection of eigenvalues $\Lambda$.

We then study the properties of $\vert \Bun_D \vert^\inter$. It always contains the open locus of points which have automorphism group $\mb G_m$ (Corollary \ref{minimal locus inside interesting locus}). Further, its intersection $\vert \Bun_D^d \vert^\inter$ with the connected component $\Bun_D^d$ of $\Bun_D$ consisting of parabolic bundles having underlying bundle of degree $d$ is independent of $d$ provided the pair $(D,\Lambda)$ has a particular form. To make this precise, we assume the presence of parabolic structure at points $x_1,\ldots, x_s$ and write the datum $(D,\Lambda)$ of rank $n$ as $(D_1^\tw,\ldots, D_s^\tw)$, where each $D_i^\tw$ is a pair $(D_i,\Lambda_i)$, where $D_i = (r_{i1},\ldots, r_{in_i})$ is a partition of $n$ of length $n_i$ and $\Lambda_i = (\lambda_{i1},\ldots,\lambda_{in_i}) \in \l \mb C^{\times} \r^{n_i}$ is a $n_i$-dimensional vector of non-zero numbers. Then we show:

\begin{Theo}[Theorem \ref{Proposition interesting locus and hecke operators}]
    Let $(D,\Lambda) = (D_1^\tw,\ldots, D_s^\tw) \in \TwParDat_{s,\ul n}$ and $\lambda \in \mb C^{\times}$. Assume that \\
    $D_s^\tw = ((1,\ldots, 1),(\lambda,\ldots, \lambda))$. Then $\mb T^{(l)}_s$ induces a homeomorphism
    \begin{align*}
        \vert \Bun_D^d \vert^\inter &\simeq \vert \Bun_D^{d-1} \vert^\inter.
    \end{align*}
\end{Theo}

On the spectral side, the condition on $D_s^\tw$ corresponds to regular monodromy with a single eigenvalue.

\subsection{Main Results II: An example of $\mb P^1$ minus four points}
In the second part of this article, we study the structure and geometry of $\vert \Bun_D \vert^\inter$ in a specific example: $X = \mb P^1$, $\rk = 3$ and parabolic structure being present at four distinct points $0,1,t,\infty \in \mb P^1(k)$. At $0, 1$ and $t$ partial flags of type $(1,2)$ (i.e. consisting of a line in a three-dimensional vector space) and equivariance datum $(\lambda,1)$ are considered and the structure at $\infty$ consists of $(1,1,1)$ (a full flag) and equivariance datum $(\lambda^{-1},\lambda^{-1},\lambda^{-1})$. These data constitute a \textit{twisted parabolic datum} in the sense of Definition \ref{definition twisted parabolic datum}.

The first main result is an explicit description of the space $\vert \Bun_D^0 \vert^\inter$ of those parabolic bundles in $\vert \Bun_D \vert^\inter$, which are defined over $k$ and whose underlying vector bundle has degree $0$:

\begin{Theo}[Theorem \ref{Classification Theorem for degree zero}] \label{Intro Classification Theorem for degree zero}
    The degree $0$ interesting locus $\vert \Bun_D^0(k) \vert^\inter$ defined over $k$ consists of 
    \begin{itemize}
        \item[a)] a family of points $P_\alpha^0$ indexed by $\alpha \in \l \mb P^1 \setminus \{0,1,\frac{t-1}{t},\infty\}\r(k)$
        \item[b)] three points $P_i^{0,\vee}$, $i \in \{0,1,\infty\}$,
        \item[c)] three points $P_i^{0,\wedge}$, $i \in \{0,1,\infty\}$,
        \item[d)] three points points $P_i^{0,\dec}$, $i \in \{0,1,\infty\}$,
        \item[e)] one point $P_{\frac{t-1}{t}}^{0,\spc}$.
    \end{itemize}
    Each of these points except the points in d) and e) has automorphism group isomorphic to $\mb G_m$. The points in d) have automorphism group isomorphic to $\mb G_m \times \mb G_m$. The point $P_{\frac{t-1}{t}}^{0,\spc}$ has automorphism group isomorphic to $\mb G_m \times \mb G_a$.
\end{Theo}

The points mentioned in Theorem \ref{Intro Classification Theorem for degree zero} are explicitly defined in coordinates in section \ref{Subsection about degree zero locus summary}. The index set in a) is chosen compatibly with the morphism $\mf d^{\interc}$ in part 3) of Theorem \ref{Big Theorem Introduction} below. Using rank one lower Hecke modifications at the point $\infty$, we obtain a similar descriptions for the degree $1$ component (Theorem \ref{Classification Theorem for degree one}) from which one reads off:

\begin{Prop}[Proposition \ref{First Proposition about Geometry}]\label{underlying bundle statement}
    The underlying bundle of each two points in $\vert \Bun_D^1(k) \vert^\inter$ is isomorphic to $\mc O(1) \oplus \mc O \oplus \mc O$. 
\end{Prop}

An analogous statement for the degree $0$ component is false. We use Proposition \ref{underlying bundle statement} to write down a convenient chart around the point $P_{\frac{t-1}{t}}^{1,\spc}$ corresponding to $P_{\frac{t-1}{t}}^{0,\spc}$ via the heomeomorphism $\vert \Bun_D^1 \vert^\inter \simeq \vert \Bun_D^0 \vert^\inter$ (Lemma \ref{Lemma theta one}). Combining the description of the degree $0$ interesting locus (Theorem \ref{Intro Classification Theorem for degree zero}) and the corresponding statement for the degree $1$ locus, we construct an open substack $\Bun_D^\interonec$ of $\Bun_D^1$ such that $\vert \Bun_D^\interonec \vert \supset \vert \Bun_D^1 \vert^\inter$ (Definition \ref{Definition Buninterc}). We obtain the following results on the geometry of $\Bun_D^\interonec$ and $\vert \Bun_D \vert^\inter$:

\begin{Theo}[Propositions \ref{Long Proposition about intereseting locus in degree one compoenent}, \ref{First Proposition about Geometry}, \ref{Proposition coarse moduli}, Lemma \ref{Lemma about gluing of morphisms}] \label{Big Theorem Introduction}
$\Bun_D$ has the following properties:
    \begin{itemize}
        \item[1)] $\vert \Bun_D \vert^\inter$ is not locally closed inside $\vert \Bun_D \vert$.
        \item[2)] The topological space $\vert \Bun_D^0 \vert^\inter$ contains the topological space $\vert \Bun_D^\internull \vert$ of an open substack $\Bun_D^\internull$ of $\Bun_D$ and the complement $\vert \Bun_D^0 \vert^\inter \setminus \vert \Bun_D^\internull \vert$ consists of a single point, namely $P_{\frac{t-1}{t}}^{0,\spc}$ in Theorem \ref{Intro Classification Theorem for degree zero}. The point $P_{\frac{t-1}{t}}^{0,\spc}$ is up to isomorphism the unique point of $\Bun_D^0$ which is indecomposable and has automorphism group of dimension $> 1$.
        \item[3)] $\Bun_D^0$ admits a morphism
            \begin{align*}
                    \mf d^{\interc} : \Bun_D^{\interc} \to \mb P^1,
            \end{align*}
            whose restriction to the open substack $\Bun_D^{0,\mini}$ of $\Bun_D^0$ consisting of points with automorphism group isomorphic to $\mb G_m$ is smooth and fits into a commutative diagram
            \[
            \begin{tikzcd}
\Bun_D^\minnull \arrow[r, "\mf c^{(0)}"] \arrow[rd, "{\mf d^{0,\interc}}"'] & \mc P\p \arrow[d, "p\p"] \\
                                                                                  & \mb P^1                 
\end{tikzcd}
            \]
    where $\mc P\p = \mb P^1 \sqcup_{\mb P^1 \setminus \{0,1,\infty\}} \mb P^1$ is a non-separated version of $\mb P^1$ with three double points and $p\p : \mc P\p \to \mb P^1$ is the natural map.
    \item[4)] The morphism $\mf c^{(0)}$ is a $\mb G_m$-gerbe and exhibits $\mc P\p$ as the coarse moduli space of $\Bun_D^{0,\mini}$.
    \end{itemize}
\end{Theo}

Finally, we compute the Picard group of the rigidification $\Bun_D^{0,\interc} \rigid \mb G_m$ of $\Bun_D^{0,\interc}$ for the central embedding of $\mb G_m$ into the automorphism groups of points of $\Bun_D^{0,\interc}$:

\begin{Theo} [Theorem \ref{Theorem Picard group of interc}]
    The Picard group of $\Bun_D^{0,\interc} \rigid \mb G_m$ is isomorphic to $\mb Z^{\oplus 6}$.
\end{Theo}

\subsection{Acknowledgements}
This article is heavily influenced by the results obtained in the ramified setting by uit de Bos \cite{uitdeBosThesis}, Arinkin \cite{ArinkinOrthogonality}, Yun \cite{YunRigidity} and Heinloth \cite{HeinlothThesis}. I thank Konstantin Jakob for countless helpful discussions and Timo Richarz for his general advise and his interest in the subject of this article. 
I acknowledge support (through Timo Richarz) by the European Research Council (ERC) under Horizon Europe (grant agreement nº 101040935), by the Deutsche Forschungsgemeinschaft (DFG, German Research Foundation) TRR 326 \textit{Geometry and Arithmetic of Uniformized Structures}, project number 444845124 and the LOEWE professorship in Algebra, project number LOEWE/4b//519/05/01.002(0004)/87. This project is further supported (through Konstantin Jakob) by the DFG Research Grant „Katz’s middle convolution for automorphic sheaves“, project number 566801746.

\subsection{AI Declaration}
The research questions and overall strategy of the paper were developed by the author. During the preparation of this manuscript, the author used ChatGPT (OpenAI) to review and improve the exposition. In addition, ChatGPT was used to assist with certain cohomology calculations at the end of Section \ref{Last section}. All other mathematical results were obtained without AI assistance. In particular, the calculations in Section \ref{section on classification} were carried out entirely without AI assistance. The author takes full responsibility for the mathematical correctness and content of the paper.

\subsection{Notation}
\begin{itemize}
    \item The symbol $\simeq$ is used to denote isomorphisms.
    \item Our notation involves lower and upper indices. The $n$-fold self-product resp. self-tensor product of an object is denoted by $M^{\times n}$ resp. $M^{\otimes n}$.
    \item If $\mf X$ is an algebraic stack over a field $k$ and $K$ is a field extension of $k$, we denote by $\vert \mf X \vert$ the topological space of $\mf X$ and by $\vert \mf X(K) \vert$ the image of the natural functor $\mf X(K) \to \vert \mf X \vert$.
    \item Let $k$ be a field. We will work with coordinates $[x:y]$ on $\mb P^1 = \mb P^1_k$. Whenever we embed $\mb A^1 = \mb A^1_k$ into $\mb P^1$ and do not specify the embedding we use $x \mapsto [x:1]$.
    \item For a $k$-vector space $V$ we denote by $\mb P(V)$ the projective space of $V$ (parametrizing one-dimensional subspaces) and by $\FL(V)$ the flag variety (parametrizing full flags in $V$). We set $\mb P^2 := \mb P(k^{\oplus 3})$ and $\FL_3 := \FL(k^{\oplus 3})$.
    \item For a scheme $T$ and a $T$-valued point $x$ of a stack $\mf X$, we denote by $\ul{\Aut}(x)$ the automorphism sheaf of $x$ on the big étale site over $T$. We denote by $\Aut(x)$ the sections of $\ul{\Aut}(x)$ over $T$. If $\mf X$ is an algebraic stack with quasi-separated diagonal and $T$ the spectrum of a field, $\ul{\Aut}(x)$ is representable by a scheme \cite[Corollary~4.5.31]{Alper}.
    \item If $X$ is a curve over a field $k$, $x \in X(k)$ a $k$-valued point and $T$ any $k$-scheme, we denote the natural inclusion of the ideal sheaf of $\{x\} \times T$ into the structure sheaf by $\iota_x : \mc O_{X_T}(-x) \hookrightarrow \mc O_{X_T}$. 
\end{itemize}

%% file: 2_Local_and_Global_Data.tex
\section{Numerical data associated with parabolic bundles} \label{section data}
Let $X$ be a proper smooth connected curve over an algebraically closed field $k$ and $\mc E$ a vector bundle on $X$, i.e. a sheaf of $\mc O_X$-modules, which is locally on $X$ isomorphic to $\mc O_X^{\oplus m}$ for some $m \in \mb Z_{\ge 0}$. With $\mc E$ we can associate its rank $\rk(\mc E)$ and its Harder-Narasimhan type. Now let $\mc E\b$ be a parabolic vector bundle (see \ref{Definition Parabolic Sheaves} for the definition). It has an underlying vector bundle $\UL(\mc E\b)$ (ibid.) and hence we can define the rank and Harder-Narasimhan type of $\mc E\b$ as the rank and Harder-Narasimhan type of $\UL(\mc E\b)$. We refer to these data as \textit{global data} associated with $\mc E\b$. 

On the other hand, one can associate with $\mc E\b$ certain data, which describe its parabolic structure and to which we refer as its \textit{parabolic datum}.

If a global datum $D_\glob$ and a parabolic datum $D_\para$ are associated with the same parabolic vector bundle $\mc E\b$ the rank of $D_\glob$ is compatible with $D_\para$: One may think of the parabolic vector bundle $\mc E\b$ as the datum of the vector bundle $\mc E := \UL(\mc E\b)$ plus a collection of flags in the fibers of $\mc E$ at finitely many closed points in $X$. Then $D_\para$ consists of collections of dimensions of subquotients of these flags. In particular, for a fixed point, the sum of these numbers has to be equal to $\rk(\mc E)$, i.e. the rank of $D_\glob$. If this condition is satisfied for an arbitrary global datum $D_\glob$ and an arbitrary parabolic datum $D_\para$, we call $D_\glob$ and $D_\para$ \textit{compatible}.

Consider a parabolic datum $D$ of rank $r$. With $D$ we will associate the stack $\Bun_D$ of parabolic bundles with datum $D$. Associated with $D$ there are standard Levi subgroups $L_1, \ldots, L_s$ of $\GL_r$ and a natural torsor $\Bun_D^+ \to \Bun_D$ under the group $L := \prod_{i=1}^s L_i$. We will be interested in twisted sheaves for this torsor, i.e. sheaves on $\Bun_D^+$, which are twisted equivariant for some character sheaf $\chi$ on $L$, since the category of twisted sheaves is the place to look for Hecke eigensheaves (cf. \ref{Motivational part of intro}). In order to also specify the equivariance condition, we will extend the parabolic datum $D_\para$ to what we will refer to as a \textit{twisted parabolic datum}. 

The aim of this section is to set up notation for all these data and describe their relations. Our definition of parabolic data and twisted parabolic data makes sense for any curve. The definition of global data is already adapted to the case $X = \mb P^1$, our main application: By the Grothendieck-Birkhoff theorem, every vector bundles on $\mb P^1$ is a direct sum of line bundles, hence the Harder-Narasimhan filtration of a vector bundle of rank $r$ is uniquely determined by a non-decreasing sequence of integers $(d_1, \ldots, d_r)$.

\subsection{Global data} \label{subsection Global Data}
We start with the definition of the set of \textit{global data}. For a map $f : \mb Z \to \mb Z_{\ge 0}$ we define its \textit{support} by 
\begin{align*}
    \supp(f) := \{z \in \mb Z \ \vert \ f(z) \neq 0\}.
\end{align*}

\begin{Def}
    The \textit{set of global data} is the set $\GlobDat := \Map_\fin(\mb Z, \mb Z_{\ge 0})$ of maps $\ul d : \mb Z \to \mb Z_{\ge 0}$ with finite support.
\end{Def}

The map
\begin{align*}
    \rk : \GlobDat &\to \mb Z_{\ge 0} \\
    \ul d &\mapsto \sum_{z \in \mb Z} \ul d(z)
\end{align*}
is called the \textit{rank map}. For $r \in \mb Z_{\ge 0}$ we define the set
\begin{align*}
    \GlobDat^r := \rk^{-1}(r)
\end{align*}
of \textit{global data of rank} $r$. 

The map
\begin{align*}
    \deg : \GlobDat &\to \mb Z \\
    \ul d &\mapsto \sum_{z \in \mb Z} z \ul d(z)
\end{align*}
is called the \textit{degree map}.

We endow $\GlobDat$ with the partial order
\begin{align*}
    \ul d \le \ul d\p \ \Leftrightarrow \ \forall z \in \mb Z \ : \ \ul d(z) \le \ul d\p(z).
\end{align*}

For an element $\ul d \in \GlobDat$ we think of its support $\supp(\ul d)$ as the set of degrees of subquotients occuring in the Harder-Narasimhan filtration of a given vector bundle on the curve $X$. For $z \in \mb Z$ the non-negative number $\ul d(z)$ records the multiplicity of the degree $z$ in the filtration. We formalize this intuition as follows:

Let $\Seq^+ \subset \bigsqcup_{m \in \mb Z_{\ge 1}} \mb Z^m$ be the set of strictly increasing finite sequences (of variable length $m \ge 1$) and $\Seq_{> 0} := \bigsqcup_{m \in \mb Z_{\ge 1}}  \l\mb Z_{> 0} \r^m$ the set of finite sequences (with variable length) in $\mb Z_{>0}$. Consider the map
\begin{align*}
    \size : \bigsqcup_{m \in \mb Z_{\ge 1}} \mb Z^m \to \mb Z_{\ge 1},
\end{align*}
which sends a finite sequence to its length. We denote the restriction of size to $\Seq^+$ and $\Seq_{> 0}$ by $\size$ again.

\begin{Lem} \label{Lemma global data formalism}
    There are mutually inverse bijections
    \begin{align*}
        \GlobDat\ & \substack{\longrightarrow\\[-0.5ex]\longleftarrow}\ \Seq^+ \times_{\size,\mb Z_{\ge 1}, \size} \Seq_{> 0} \\
        \ul d &\mapsto \l(d_1,\ldots,d_m),(t_1,\ldots,t_m) \r \\
        \sum_{i = 1}^m t_i\p \delta_{d_i\p} &\mapsfrom \l(d_1\p,\ldots,d_m\p),(t_1\p,\ldots,t_m\p) \r,
    \end{align*}
    where
    \begin{itemize}
        \item $(d_1,\ldots,d_m) \in \Seq^+$ is the unique element such that $\supp(\ul d) = \{d_1,\ldots, d_m\}$,
        \item $t_i = \ul d(d_i)$,
        \item $\delta_z$ denotes the Kronecker delta at $z \in \mb Z$ (i.e. $\delta_z(z\p) = 1$ if $z = z\p$ and $0$ otherwise).
    \end{itemize}
\end{Lem}

If $\ul d \in \GlobDat$ maps to  $\l(d_1,\ldots,d_m),(t_1,\ldots, t_m)\r$ under the map in Lemma \ref{Lemma global data formalism}, we denote $\ul d$ also by
\begin{align*}
    \ul d = \l d_1^{[t_1]}, \ldots, d_m^{[t_m]} \r
\end{align*}
and call $\secmul(d_1,\ldots,d_m) \in \Seq^+$ the \textit{sequence of degrees} and $\secmul(t_1,\ldots,t_m) \in \Seq_{>0}$ the \textit{sequence of multiplicities} of $\ul d$.

\subsection{Parabolic data} \label{subsection parabolic data}
For this subsection we fix $s \in \mb Z_{\ge 1}$ and $\ul n = (n_1,\ldots, n_s) \in (\mb Z_{\ge 2})^{\times s}$. One should think of $s$ as the number of points at which a non-trivial parabolic structure is placed and of $n_i$ as being the number of steps in the parabolic structure at the $i$-th point. For each $n \in \mb Z_{\ge 1}$ we define
\begin{align*}
    \rk_\para^{(n)} : \mb Z_{\ge 0}^n \setminus \{0\} &\to \mb Z_{\ge 1} \\
    (z_1,\ldots,z_n) &\mapsto \sum_{i=1}^n z_i.
\end{align*}

\begin{Def}
The set
\begin{align*}
    \ParDat_{s,\ul n} := \left\{(D_1,\ldots, D_s) \in \prod_{i=1}^s \l (\mb Z_{\ge 0})^{\times n_i} \setminus \{0\} \r \vert \ \rk_\para^{(n_1)}(D_1) = \ldots = \rk_\para^{(n_s)}(D_s)\right\}
\end{align*}
    is called the \textit{set of parabolic data of type} $(s,\ul n)$.
\end{Def}

Each vector $D_i$ should be thought of as storing the datum of a \textit{flag type}, i.e. the sequence of dimensions of subquotients of a flag. 

We define the \textit{parabolic rank map} by
\begin{align*}
    \rk_\para : \ParDat_{s,\ul n} &\to \mb Z_{\ge 1} \\
    D = (D_1,\ldots, D_s) &\mapsto \rk_\para^{(n_1)}(D_1)
\end{align*}
and set for $r \in \mb Z_{\ge 1}$
\begin{align*}
    \ParDat_{s,\ul n}^r := \rk_\para^{-1}(r).
\end{align*}

We endow $(\mb Z,+)$ with the structure of a partially ordered semigroup and $\prod_{i=1}^s \mb Z^{\times n_i} = \mb Z^{\times \l\sum_i n_i \r}$ with the product structure. Along the natural inclusion $\ParDat_{s,\ul n} \subset \prod_{i=1}^s \mb Z^{n_i}$, $\ParDat_{s,\ul n}$ inherits the structure of a partially ordered semigroup. Note that for $D, D\p \in \ParDat_{s,\ul n} \subset \prod_{i=1}^s \mb Z^{\times n_i}$ with $D\p < D$ the difference $D - D\p$ in $\prod_{i=1}^s \mb Z^{n_i}$ lies again in $\ParDat_{s,\ul n}$. We obtain an involution
\begin{align*}
    \sub_D : \ParDat_{s,\ul n, < D} &\to \ParDat_{s,\ul n, <D} \\
    D\p &\mapsto D - D\p,
\end{align*}
where $\ParDat_{s,\ul n, < D}$ denotes the set of $D\p \in \ParDat_{s,\ul n}$ with $D\p < D$.

Finally, note that in the definition of $\ParDat_{s,\ul n}$ we allowed entries in the vectors $D_i$ constituting $D \in \ParDat_{s,\ul n}$ to be $0$. This means that not every inclusion in the corresponding flag subquotients of dimension $0$ may occur. If this is not the case, we call the datum $D$ \textit{regular}:

\begin{Def}
    Let $r \in \mb Z_{\ge 1}$.
    \begin{itemize}
        \item[i)] The set of \textit{regular parabolic data of type} $(s,\ul n)$ is
    \begin{align*}
        \ParDat_{s,\ul n}^\reg := \ParDat_{s,\ul n} \cap \prod_{i=1}^s \mb Z_{> 0}^{n_i}.
    \end{align*}
    \item[ii)] The set of \textit{regular parabolic data of type} $(s,\ul n)$ \textit{of rank} $r$ is
    \begin{align*}
        \ParDat_{s,\ul n}^{r,\reg} := \ParDat_{s,\ul n}^r \cap \ParDat_{s,\ul n}^\reg.
    \end{align*}
    \end{itemize}
\end{Def}

\begin{Rem}
    Let $s \in \mb Z_{\ge 1}$ and $r \in \mb Z_{\ge 1}$. For arbitrary $\ul n \in (\mb Z_{\ge 2})^s$, the set $\ParDat_{s,\ul n}^{r}$ is non-empty. This is no longer true for $\ParDat_{s,\ul n}^{\reg,r}$: Its non-emptiness is equivalent to $n_i \le r$ for all $i$.
\end{Rem}

\subsection{Associated regular data} \label{Associated regular data}
In this subsection, we want to assign to $D \in \ParDat_{s,\ul n}$ a regular datum $D^\reg$ in a natural way. Roughly, $D^\reg$ is obtained from $D = (D_1,\ldots,D_s)$ by "deleting all zeros in the vectors $D_i$". In doing so, $s$ and $\ul n$ will change.

More precisely, let $B$ be the set over $\mb Z_{\ge 1}$ with fiber over $s \in \mb Z_{\ge 1}$ equal to $\l \mb Z_{\ge 2} \r^{\times s}$. We write as typical element in $B$ as $(s,\ul n)$. Let $\ParDat$ be the set over $B$ whose fiber over $(s,\ul n)$ is $\ParDat_{s,\ul n}$. We denote a typical element of $\ParDat$ as $(s,\ul n,D)$. Let $B^{\reg,\emptyset}$ be the set obtained by formally adjoining another element $\emptyset$ to $B$:
\begin{align*}
    B^{\reg,\emptyset} := B \sqcup \{\emptyset\}.
\end{align*}
Let $\ParDat^{\reg,\emptyset}$ be the set over $B^{\reg,\emptyset}$ whose fiber over $(s,\ul n) \in B$ is $\ParDat_{s,\ul n}^\reg$ and whose fiber over $\emptyset$ is $\emptyset$.

\begin{Def} \label{Definition associated regular datum}
    We define a map
    \begin{align*}
        \Reg : \ParDat \to \ParDat^{\reg,\emptyset}
    \end{align*}
    as follows: Let $(s,\ul n, D) \in \ParDat$. Write $D = (D_1,\ldots, D_s)$ with $D_i = (r_{i1},\ldots, r_{in_i})$. For each $i \in \{1,\ldots, s\}$ define    \begin{align*}
        Z_i &:= \{j \in \{1,\ldots, n_i\} \ \vert \ r_{ij} \neq 0\} \\
        z_i &:= \#Z_i \in \mb Z_{\ge 1} 
    \end{align*}
    and set
    \begin{align*}
        S^\reg &:= \{i \in \{1,\ldots,s\} \ \vert \ z_i \ge 2\}, \\
        s^\reg &:= \# S^\reg \in \mb Z_{\ge 0}.
    \end{align*}
    If $s^\reg \ge 1$ let
    \begin{itemize}
        \item $\omega_D : \{1,\ldots, s^\reg\} \to \{1,\ldots,s\}$ be the unique order-preserving map with image $S^\reg$,
        \item $n_i^\reg := z_{\omega_D(i)}$ for $i \in \{1,\ldots, s^\reg\}$,
        \item $\omega_{D,i} : \{1,\ldots,n_i^\reg\} \to \{1,\ldots, n_i\}$ for $i \in \{1,\ldots, s^\reg\}$ the unique order-preserving map with image $Z_i$,
        \item $r^\reg_{ij} := r_{\omega_D(i) \omega_{D,i}(j)}$ for $i \in \{1,\ldots,s^\reg\}$, $j \in \{1,\ldots,n_i^\reg\}$,
        \item $D_i^\reg := \l d^\reg_{i1},\ldots, d^\reg_{in_i^\reg} \r \in \mb Z_{\ge 1}^{n_i^\reg}$ for $i \in \{1,\ldots,s^\reg\}$.
    \end{itemize}
    Define
    \begin{align*}
        \Reg(s,\ul n,D) := \begin{cases}
            \emptyset, & s^\reg = 0 \\
            \l s^\reg, \l n_1^\reg,\ldots, n_{s^\reg}^\reg \r, \l D_1^\reg, \ldots, D_{s^\reg}^\reg \r \r, & s^\reg > 0.
        \end{cases}
    \end{align*}.
\end{Def}

Non-regular parabolic data occur naturally when one wants to study parabolic subbundles, as e.g. in Definitions \ref{Main Definition of irrelevance recognizing} and \ref{Main Definition noncupsidality recognizing}: If one considers a parabolic bundle which has regular parabolic type\footnote{The parabolic type of a parabolic vector bundle is given by a parabolic datum.}, a parabolic subbundle will in general have non-regular parabolic type. However, the stack of parabolic vector bundles of a given parabolic type depends only on the associated regular datum, see Remark \ref{Remark about isomorphisms between stacks having different data}.

\subsection{Twisted parabolic data}
As in subsection \ref{subsection parabolic data} we fix  $s \in \mb Z_{\ge 1}$ and $\ul n = (n_1,\ldots, n_s) \in (\mb Z_{\ge 2})^{\times s}$.

\begin{Def} \label{definition twisted parabolic datum}
    The set of \textit{twisted parabolic data of type} $(s,\ul n)$ is
    \begin{align*}
        \TwParDat_{s,\ul n} := \left\{ (D_i,\Lambda_i)_i \in \prod_{i=1}^s \l\l (\mb Z_{\ge 0})^{\times n_i} \setminus \{0\} \r \times (\mb C^\times)^{\times n_i}\r \ \vert \ (D_1,\ldots,D_s) \in \ParDat_{s,\ul n} \right\}.
    \end{align*}
\end{Def}

The set $\TwParDat_{s,\ul n}$ admits per definition a natural projection
\begin{align*}
    \pr_{s,\ul n} : \TwParDat_{s,\ul n} \to \ParDat_{s,\ul n}.
\end{align*}
For $r \in \mb Z_{\ge 1}$ and $D \in \ParDat_{s,\ul n}$ we set
\begin{align*}
    \TwParDat_{s,\ul n}^r &:= \pr_{s,\ul n}^{-1}\l \ParDat_{s,\ul n}^r \r \\
    \TwParDat_{s,\ul n, < D} &:= \pr_{s,\ul n}^{-1}\l \ParDat_{s,\ul n, < D} \r.
\end{align*}

\begin{Def}
    We define a map
    \begin{align*}
        \mf R : \TwParDat_{s,\ul n} \to \mb C^{\times}
    \end{align*}
    as follows: For $D^\tw = (D,\Lambda) \in \TwParDat_{s,\ul n}$ with $D = (r_{ij})_{ij}$ and $\Lambda = (\lambda_{ij})_{ij}$ we set
    \begin{align*}
        \mf R\l D^\tw \r := \prod_{i=1}^s \prod_{j=1}^{n_i} \lambda_{ij}^{r_{ij}}.
    \end{align*}
\end{Def}

We use this map to define the notion of \textit{regularity} for twisted data:

\begin{Def}
    Let $D^\tw = (D,\Lambda) \in \TwParDat_{s,\ul n}$.
    \begin{itemize}
        \item[i)] $D^\tw$ is called \textit{weakly regular} iff $\mf R(D^\tw) = 1$.
        \item[ii)] $D^\tw$ is called \textit{regular} iff it is weakly regular and $D$ is regular.
    \end{itemize}    
\end{Def}

We denote by $\TwParDat_{s,\ul n}^\reg$ resp. $\TwParDat_{s,\ul n}^\wreg$ the set of regular resp. weakly regular twisted parabolic data of type $(s, \ul n)$. Similarly we define $\TwParDat_{s,\ul n}^{\reg, r}$ resp. $\TwParDat_{s,\ul n}^{\wreg, r}$ to be the respective subset of those twisted data having rank $r$.

Given $D^\tw = (D,\Lambda) \in \TwParDat_{s,\ul n}$, every $D\p \in \ParDat_{s,\ul n,<D}$ inherits the structure of a twisted parabolic datum:

\begin{Def}
    Given $D^\tw = (D,\Lambda) \in \TwParDat_{s,\ul n}$, we define a map
    \begin{align*}
    \Res_{D^\tw} : \ParDat_{s,\ul n, < D} &\to \TwParDat_{s,\ul n, < D} \\
    D\p &\mapsto (D\p,\Lambda).
\end{align*}
\end{Def}

Finally, the map $\Res_{D^\tw}$ can be used to define the condition of being \textit{irrelevance recognizing}. It will be used in Definition \ref{Main Definition of irrelevance recognizing} to define the corresponding notion for parabolic vector bundles.

\begin{Def}
    Let $D^\tw \in \TwParDat_{s,\ul n}$. We call a datum $D\p \in \ParDat_{s,\ul n, < D}$ $D^\tw$-\textit{irrelevance recognizing} iff $\Res_{D^\tw}(D\p)$ is not weakly regular, i.e. $\mf R(\Res_{D^\tw}(D\p)) \neq 1$. We denote the set of $D^\tw$-irrelevance recognizing data by $\ParDat_{s,\ul n, < D}^{D^\tw-\irec}$.
\end{Def}

%% file: 3_Stacks_of_Parabolic_Bundles.tex
\section{Stacks of parabolic vector bundles} \label{stacks of parabolic vector bundles}
Let $k$ be an algebraically closed field and $X$ a connected smooth projective curve over $k$, which will be fixed throughout this section and not occur in the notation. Let $s \in \mb Z_{\ge 1}$ and $\ul n = (n_1,\ldots, n_s) \in \mb (\mb Z_{\ge 2})^{\times s}$. We further fix an ordered length $s$ tuple $S^\oro = (x_1,\ldots, x_s)$ of distinct $k$-points of $X$. We will denote by $S = \{x_1,\ldots, x_s\} \subset X(k)$ the corresponding unordered set.

Depending on these data we define for each $D \in \ParDat_{s,\ul n}$ an algebraic stack $\Bun_{S^\oro,D}$ over $k$, which classifies parabolic vector bundles on $X$ having parabolic structure $D$ at the points $x_i$. We start by defining the notion of \textit{parabolic vector bundle} in the case of parabolic structure at a single point:

\begin{Def}
    Let $T$ be a $k$-scheme, $x \in X(k)$ and $n \in \mb Z_{\ge 2}$. 
    \begin{itemize}
        \item[1)] An $(x,n)$-parabolic $\mc O_{X_T}$-module consists of the following data:
    \begin{itemize}
        \item[i)] A diagram
        \begin{align*}
            \mc F^{(-n+1)} \xrightarrow{\phi^{(-n+1)}} \mc F^{(-n+2)} \xrightarrow{\phi^{(-n+2)}} \ldots \xrightarrow{\phi^{(-1)}} \mc F^{(0)}
        \end{align*}
        in the category of $\mc O_{X_T}$-modules,
        \item[ii)] A morphism $\phi^{(-n)} : \mc F^{(0)}(-x) \to \mc F^{(-n+1)}$ of $\mc O_{X_T}$-modules
    \end{itemize}
    such that
    \begin{itemize}
        \item[a)] The composition $\phi^{(-1)} \circ \ldots \circ \phi^{(-n)}$ is isomorphic to\footnote{Recall that $\iota_x : \mc O_{X_T}(-x) \hookrightarrow \mc O_{X_T}$ denotes the inclusion of the ideal sheaf of $\{x\} \times T$.} $\mc F^{(0)} \otimes \iota_x$.
        \item[b)] For each $j \in \{1,\ldots, n-1\}$ the composition $\phi^{(-(j+1))} \circ \ldots \circ \phi^{(-n)} \circ (\phi^{(-1)} \otimes \id_{\mc O(-x)}) \circ \ldots \circ (\phi^{(-j)} \otimes \id_{\mc O(-x)})$ is isomorphic to $\id \otimes \iota_x : \mc F^{(-j)} \otimes \mc O(-x) \to \mc F^{(-j)}$.
    \end{itemize}
    We will use the notation $(\mc F\b,\phi\b)$ or just $\mc F\b$ for parabolic sheaves and always define
    \begin{align*}
        \mc F^{(-n)} := \mc F^{(0)}(-x).
    \end{align*}
    \item[2)] An $(x,n)$-parabolic $\mc O_{X_T}$-module is called \textit{quasi-coherent} resp. \textit{coherent} resp. \textit{locally free} iff for all $j \in \{0,\ldots, n-1\}$, $\mc F^{(-j)}$ has the respective property. 
    \item[3)] Let $\mc F\b$ be an $(x,n)$-parabolic $\mc O_{X_T}$-module. Then we call $\mc F^{(0)}$ its \textit{underlying sheaf} and denote it by $\UL(\mc F\b).$
    \end{itemize}
\end{Def}

Now we consider the general case:

\begin{Def} \label{Definition Parabolic Sheaves}
     Let $T$ be a $k$-scheme. \\
    An $(S^\oro,\ul{n})$-parabolic $\mc O_{X_T}$-module is a family $(\mc F\b_i)_{i \in \{1,\ldots,s\}}$ where each $\mc F\b_i$ is a $(x_i,n_i)$-parabolic $\mc O_{X_T}$-module such that $\UL(\mc F\b_i) = \UL(\mc F\b_j)$ for all $i,j \in \{1,\ldots, s\}$. We denote this sheaf by $\UL((\mc F\b_i)_i)$. The $(S^\oro,\ul{n})$-parabolic $\mc O_{X_T}$-module $(\mc F\b_i)_i$ is called \textit{quasi-coherent} resp. \textit{coherent} resp. \textit{locally free} iff for all $i \in \{1,\ldots, s\}$, $\mc F\b_i$ has the respective property. We denote the corresponding full subcategories by
    \begin{align*}
        \ParMod_{S^\oro,\ul n}(T) \supset \ParQCoh_{S^\oro,\ul n}(T) \supset \ParCoh_{S^\oro,\ul n}(T) \supset \ParBun_{S^\oro,\ul n}(T).
    \end{align*}
    For $\mc F\b = (\mc F\b_i)_i \in \ParMod_{S^\oro,\ul n}(T)$, the $\mc O_{X_T}$-module $\UL(\mc F\b ) := \UL(\mc F\b_i)$, which is independent of the choice of $i$, is called the \textit{underlying sheaf} of $\mc F\b$. If $\mc F\b$ is a parabolic vector bundle, we call this sheaf the \textit{underlying bundle}.
\end{Def}

We will call a locally free $(S^\oro,\ul n)$-parabolic $\mc O_{X_T}$-module also a $T$\textit{-family of} $(S^\oro,\ul n)$\textit{-parabolic vector bundles} on $X$. In the case $T = \Spec(k)$ we just call it an $(S^\oro,\ul n)$\textit{-parabolic vector bundle} on $X$. \\
By a \textit{monomorphism of parabolic vector bundles} we will mean a monomorphism in $\ParBun_{S^\oro,\ul n}(T)$ and we will call a morphism in $\ParBun_{S^\oro,\ul n}(T)$ whose cokernel is again a parabolic vector bundle an \textit{injection of parabolic vector bundles}.

\begin{Def}
    For a $k$-scheme $T$ we define 
    \begin{align*}
        \Bun_{S^\oro,\ul n}(T) := \ParBun_{S^\oro,\ul n}(T)^{\simeq},
    \end{align*}
    where $\ParBun_{S^\oro,\ul n}(T)^{\simeq}$ is the largest subgroupoid of the category $\ParBun_{S^\oro,\ul n}(T)$.
\end{Def}

The stack of vector bundles on $X$ (of arbitrary rank) without parabolic structure will be denoted by $\Bun_{\emptyset}$. 

Let $x \in X(k)$, $T = \Spec(K)$ for a field extension $K$ of $k$ and let $\mc F\b$ be $T$-family of $(x,n)$-parabolic vector bundle on $X$. We set
\begin{align*}
    \rk(\mc F\b) := \rk(\UL(\mc F\b)).
\end{align*}
The natural map of sheaves $\iota_x : \mc O_{X_T}(-x) \to \mc O_{X_T}$ is injective and since $\mc F^{(0)}$ is locally free, also $\mc F^{(0)} \otimes \iota_x$ is injective. Thus, all maps $\phi^{(-i)}$, $i \in \{1,\ldots, n\}$, are injective. Further, $\mc F^{(0)} \otimes \iota_x$ is an isomorphism away from $x_T$, so the locally free sheaves $\mc F^{(-i)}$ do all have the same rank and $\coker(\phi^{(-i)})$ is a torsion sheaf supported at $x_T$ for all $i \in \{1,\ldots, n\}$, in particular
\begin{align*}
    \chi\l\coker(\phi^{(-i)})\r = \dim_{k(x_T)}\l\coker(\phi^{(-i)} \otimes k(x_T))\r \ge 0.
\end{align*}
The exactness of the sequences
\begin{align*}
    0 \to \mc F^{(-i)} \to \mc F^{(-i+1)} \to \coker(\phi^{(-i)}) \to 0,
\end{align*}
the additivity of the Euler characteristic in exact sequences and the Riemann-Roch theorem imply
\begin{align*}
    \sum_{i=1}^n \chi(\coker(\phi^{(-i)})) &= \chi\l \mc F^{(0)} \r - \chi\l\mc F^{(0)}(-x)\r \\
    &= \deg\l\mc F^{(0)}\r + \rk\l\mc F^{(0)}\r(1-g) - \l\deg\l\mc F^{(0)}(-x)\r + \rk\l\mc F^{(0)}(-x)\r(1-g) \r \\
    &= \deg\l\mc F^{(0)}\r -  \deg\l\mc F^{(0)}(-x)\r  \\
    &= \rk(\mc F\b).
\end{align*}
We obtain an element
\begin{align*}
    D_x = \l\chi\l\coker(\phi^{(-n)})\r, \ldots, \chi\l\coker(\phi^{(-1)})\r \r \in \ParDat_{1,(n)}
\end{align*}
with $\rk(D_x) = \rk(\mc F\b)$. \\
If $\mc F\b$ an $(S^\oro,\ul n)$-parabolic vector bundle, performing the above construction for every point $x_i$ in $S^\oro$ we obtain an element
\begin{align*}
    D(\mc F\b) = (D_{x_1},\ldots, D_{x_s}) \in \ParDat_{s,\ul n}.
\end{align*}

\begin{Def}
    The parabolic datum $D(\mc F\b) \in \ParDat_{s,\ul n}$ is called the \textit{parabolic type} of $\mc F\b$.
\end{Def}

To understand how the datum $D(\mc F\b)$ behaves in families, we record the following Lemma:

\begin{Lem} \label{Lemma about flatness in families}
    Let $T$ be a $k$-scheme, $\mc E$, $\mc E\p$ vector bundles on $X_T$ of the same rank and $f : \mc E \to \mc E\p$ a monomorphism of quasi-coherent sheaves. Assume there is a finite subset $S \subset X(k)$ such that $\supp(\coker(f)) \subset S_T$. Then $\coker(f)$ is flat over $T$.
\end{Lem}
\begin{proof}
    The statement is local on $T$, so we may assume that $T$ is affine. Since $X_T$ is of finite presentation over $T$ and $\mc E$, $\mc E\p$ are locally free, we may use Noetherian approximation to assume that $T$ is affine locally Noetherian. Let $y \in X_T$ and denote its image in $T$ by $t$. Set $R := \mc O_{T,t}$ and denote by $\kappa$ the residue field of $R$. To show that the stalk $M := \coker(f)_y = \coker(f_y)$ is flat over $R$, it suffices to show $\Tor_1^R(\kappa, M) = 0$ \cite[Tag~\href{https://stacks.math.columbia.edu/tag/00MK}{00MK}]{StacksProject}. Since $X_T \to T$ is flat and $\mc E$, $\mc E\p$ are flat over $X_T$, we obtain $\Tor_1^R(\kappa,\mc E\p_y) = 0$, hence by the long exact Tor-sequence
    \begin{align*}
        \Tor_1^R(\kappa, M) = \ker\l\mc E_y \otimes_R \kappa \xrightarrow{f_y \otimes_R \kappa} \mc E\p_y \otimes_R \kappa\r.
    \end{align*}
    We have a natural identification
    \begin{align*}
        f_y \otimes_R \kappa = \l(\mc E\vert_{X_t})_y \xrightarrow{(f \vert_{X_t})_y} \mc (E\p\vert_{X_t})_y\r,
    \end{align*}
    so it suffices to show that $f \vert_{X_t}$ is injective. 
    
    Now $f \vert_{X_t}$ is a morphism between vector bundles of the same rank on the connected smooth projective curve $X_t$. The assumptions imply that it is generically an isomorphism. Since $X_t$ is integral, $f \vert_{X_t}$ is injective.
\end{proof}

Lemma \ref{Lemma about flatness in families} together with the fact that the Euler characteristic is locally constant in flat families implies:

\begin{Cor}
    Let $T$ be a $k$-scheme and $\mc F\b \in \ParBun_{S^\oro,\ul n}(T)$. The map
    \begin{align*}
        \vert T \vert &\to \ParDat_{s,\ul n} \\
        t &\mapsto D(\mc F\b \otimes k(t))
    \end{align*}
    is locally constant.
\end{Cor}

\begin{Def} \label{Definition of stack of parabolic vector bundles depending on D}
    Let $D \in \ParDat_{s,\ul n}$ and $T$ a $k$-scheme. We define a subgroupoid $\Bun_{S^\oro, D}(T)$ of $\Bun_{S^\oro,\ul n}(T)$ by
    \begin{align*}
        \Bun_{S^\oro, D}(T) := \{\mc F\b \in \Bun_{S^\oro,\ul n}(T) \ \vert \ \forall t \in \vert T \vert \ : \ D(\mc F\b \otimes k(t)) = D\}.
    \end{align*}
\end{Def}

If $S^\oro$ is clear from the context, we just write $\Bun_D$ instead of $\Bun_{S^\oro, D}$.

\begin{Rem} \label{Remark about isomorphisms between stacks having different data}
    Let $s$, $S^\oro$, $\ul n$ be as before and let $D \in \ParDat_{s,\ul n}$. Recall that in Definition \ref{Definition associated regular datum} we constructed a datum $(s^\reg,\ul n^\reg, D^\reg)$, which in the case $s^\reg \neq 0$ came together with an injective order-preserving map $\omega_D : \{1,\ldots, s^\reg\} \to \{1,\ldots, s\}$. We then define
    \begin{align*}
        S^{\oro,\reg} := \l \{1,\ldots, s^\reg\} \xrightarrow{\sigma_\reg} \{1,\ldots, s\} \xrightarrow{S^\oro} X(k) \r \in X(k)^{\times s^\reg}.
    \end{align*}
    Then we have natural isomorphisms of algebraic stacks
    \begin{align*}
        \Bun_{S^\oro, D} \simeq \Bun_\emptyset
    \end{align*}
    if $s^\reg = 0$ and
    \begin{align*}
        \Bun_{S^\oro, D} \simeq \Bun_{S^{\oro,\reg}, D^\reg}
    \end{align*}
    if $s^\reg \neq 0$. This follows from the fact that the groupoid of flags in a fixed vector space of type $(r_1,\ldots,r_n)$ with $r_i \ge 0$ and $r_{i_0} = 0$ is equivalent to the groupoid of flags of type $(r_1,\ldots, r_{i_0-1},r_{i_0+1},\ldots, r_n)$. 
\end{Rem}

\begin{Lem}
    For each $D \in \ParDat_{s,\ul n}$, $\Bun_D$ is a smooth, locally of finite type algebraic stack over $k$.
\end{Lem}
\begin{proof}
    It is well known that the assertions are true for $\Bun_\emptyset$. By Remark \ref{Remark about isomorphisms between stacks having different data} we may hence assume that $D \in \ParDat_{s,\ul n}^\reg$. Write $D = (D_1,\ldots, D_s)$ and let $r := \rk(D)$. Let $F_i$ be the standard flag of type $D_i$ in $k^{\oplus r}$ (meaning that the stabilizer of $F_i$ in $\GL_r$ contains the Borel subgroup of upper triangular matrices). Let $P_i$ be the stabilizer of $F_i$ inside $\GL_r$. Then $P_i$ is a parabolic subgroup of $\GL_r$. Let $\mc O_i$ resp. $k_i$ be the completed local ring resp. residue field of $\mc O_X$ at $x_i$ and let $\mc P_i$ be the preimage of $P_i$ under $\GL_{r,\mc O_i} \to \GL_{r,k_i} \simeq \GL_r$. Let $\mc G$ be the group scheme constructed by modifying the constant group scheme $\GL_{r,X}$ over $X$ along $\Spec(\mc O_i)$ by $\mc P_i$. Then there is a natural isomorphism between $\Bun_{S^\oro, D}$ and $\Bun_{\mc G}$, the stack of $\mc G$-torsors. Now the result follows from \cite[Proposition~1]{Heinloth}. 
\end{proof}

Let us consider the case $T = \Spec(k)$ and $X = \mb P^1_k$. Then by the Birkhoff-Grothendieck Theorem for every vector bundle $\mc E$ of rank $r \in \mb Z_{\ge 1}$ on $X$ there exist line bundles $\mc L_1, \ldots, \mc L_r$ unique up to isomorphism and order such that $\mc E \simeq \bigoplus_{i=1}^r \mc L_i$. Since $\Pic(\mb P^1) \simeq \mb Z$, there are $z_1, \ldots, z_m \in \mb Z$ with $z_1 < \ldots < z_m$ and $t_1, \ldots, t_m \in \mb Z_{\ge 1}$ such that
\begin{align*}
    \mc E \simeq \bigoplus_{i=1}^m \mc O(z_i)^{\oplus t_i}.
\end{align*}
Such $z_i$ and $t_i$ are uniquely determined by $\mc E$ and only depend on the isomorphism class of $\mc E$. We set
\begin{align*}
    \ULD(\mc E) := \l z_1^{[t_1]}, \ldots, z_m^{[t_m]}\r \in \GlobDat.
\end{align*}
If $\mc F\b$ is an $(S^\oro,\ul n)$-parabolic vector bundle on $\mb P^1_k$, we set
\begin{align*}
    \ULD(\mc F\b) := \ULD(\UL(\mc F\b)).
\end{align*}

\begin{Def}
    For $D \in \ParDat_{s,\ul n}$ and $\ul d \in \GlobDat$ we denote by $\Bun_D^{\ul d}$ the locally closed reduced substack of $\Bun_D$ consisting of bundles $\mc E\b$ with $\ULD(\mc E\b) = \ul d$.
\end{Def}

\subsection{Structure of an abelian category}
For any $k$-scheme $T$, the category $\ParMod_{S^\oro,\ul n}(T)$ has all direct sums (objects in $\ParMod_{S^\oro,\ul n}(T)$ are certain diagrams in $\Mod_{\mc O_{X_T}}$ and we may take direct sums vertex-wise; for this one uses compatibility of tensor product and direct sums of $\mc O_{X_T}$-modules) and the structure of an abelian category (again we may form kernels and cokernels vertex-wise and use that $\mc O(-x)$ is locally free). In particular, we obtain a functor
\begin{align*}
    \ParMod_{S^\oro,\ul n}(T) \times \ParMod_{S^\oro,\ul n}(T) \to \ParMod_{S^\oro,\ul n}(T),
\end{align*}
which preserves parabolic vector bundles and thus induces a map
\begin{align*}
    \SUM_T : \Bun_{S^\oro,\ul n}(T) \times \Bun_{S^\oro,\ul n}(T)  \to \Bun_{S^\oro,\ul n}(T).
\end{align*}
These maps are compatible with pullback, hence give rise to a morphism of stacks
\begin{align*}
    \SUM : \Bun_{S^\oro,\ul n} \times \Bun_{S^\oro,\ul n} \to \Bun_{S^\oro,\ul n}.
\end{align*}
Since taking Euler characteristic is additive in direct sums of sheaves, $\SUM$ induces for $D_1, D_2 \in \ParDat_{s,\ul{n}}$ a morphism
\begin{align*}
    \SUM : \Bun_{S^\oro,D_1} \times \Bun_{S^\oro,D_2} \to \Bun_{S^\oro, D_1+D_2}.
\end{align*}

\subsection{Parabolic bundles via flags and subbundles}
\begin{Def}
    Let $D \in \ParDat_{s, \ul n}$. For a $k$-scheme $T$ we denote by $\FlagBun_{S^\oro, D}(T)$ the groupoid of tuples $(\mc E, (\mc F_i)_{i=1,\ldots, s})$ where $\mc E$ is a vector bundle on $X_T$ and $\mc F_i$ is a flag of vector bundles (meaning that successive quotients are again vector bundles) on $\{x_i\} \times T$
    \begin{align*}
        0 = \mc F_{i,0} \subset \mc F_{i,1} \subset \mc F_{i,2} \subset \ldots \subset \mc F_{i,n_i} 
    \end{align*}
    such that $\mc F_{i,n_i} = \mc E \vert_{\{x_i\} \times T}$ and such that $\rk(\mc F_{ij}) - \rk(\mc F_{i,j-1}) = d_{ij}$ for all $i \in \{1,\ldots, s\}$ and $j \in \{1,\ldots, n_i\}$.
\end{Def}

We define a functor
\begin{align*}
    \FlagBun_{S^\oro, D}(T) \to \Bun_{S^\oro, D}(T)
\end{align*}
as follows: Let $(\mc E, (\mc F_i)_{i=1,\ldots, s}) \in \FlagBun_{S^\oro, D}(T)$ and let $i \in \{1,\ldots, s\}$. Denote by $\iota_i$ the closed embedding $\{x_i\} \times T \hookrightarrow X_T$. Let $\rho_i : \mc E \to (\iota_i)_*(\mc E \vert_{\{x_i\} \times T})$ be the surjective map of sheaves obtained by tensoring with $\mc E$ the natural map $\mc O_{X_T} \to (\iota_i)_*\mc O_{\{x_i\} \times T}$. Then for $l \in \{0,\ldots, n_i\}$ set
\begin{align*}
    \mc E_i^{(-n_i+l)} := \rho_i^{-1}(\mc F_{i,l}).
\end{align*}
By construction, for $j \in \{1,\ldots, n_i\}$ there is a natural map $\phi_i^{(-j)} : \mc E_i^{(-j)} \hookrightarrow \mc E_i^{(-j+1)}$. Then $(\mc E_i\b, \phi_i\b) \in \Bun_{S^\oro, D}(T)$. \\

Conversely, given $(\mc E_i\b, \phi_i\b) \in \Bun_{S^\oro, D}(T)$, we define $(\mc E, (\mc F_i)_{i=1,\ldots, s}) \in \FlagBun_{S^\oro, D}(T)$ by setting
\begin{align*}
    \mc E &:= \UL(\mc E_i\b, \phi_i\b), \\
    \mc F_{i,l} &:= \im\l(\phi_i^{(-1)} \circ \ldots \circ \phi_i^{(-n_i+l)})\vert_{\{x_i\} \times T}\r.
\end{align*}

\begin{Prop}
    The two constructions above are inverse to each other and can be upgraded to an isomorphism of stacks
\begin{align*}
    \FlagBun_{S^\oro, D} \simeq \Bun_{S^\oro, D}.
\end{align*}
\end{Prop}

\subsection{Inherited parabolic structure}
 Let $(\mc E\b,\phi\b) \in \ParQCoh_{S^\oro,\ul n}(T)$ and consider a pair $(\mc F,\psi)$, where $\mc F \in \QCoh(X_T)$ and $\psi : \mc F \to \UL(\mc E\b)$ is a morphism of quasi-coherent sheaves. Assume the following:
 \begin{itemize}
     \item[i)] For all $i \in \{1,\ldots, s\}$ the morphism $\psi \vert_{\{x_i\} \times T}$ is injective.
     \item[ii)]  For all $i \in \{1,\ldots, s\}$ the morphism $\mc F \otimes \mc O(-(\{x_i\} \times T)) \to \mc F$ is injective.
     \end{itemize} 
 Under these assumptions we may endow $\mc F$ with the structure of a parabolic quasi-coherent sheaf inherited via $\psi$ from the parabolic structure on $\UL(\mc E\b)$ given by $\mc E\b$. More precisely, we construct a parabolic quasi-coherent sheaf $\Inh_{\mc E\b}(\psi) = \Inh_{\mc E\b}(\mc F, \psi) = (\mc F\b,\tau\b) \in \ParQCoh_{X,S^\oro,\ul n}(T)$ together with a morphism $\psi_{\Inh,\mc E\b} : \Inh_{\mc E\b}(\psi) \to \mc E\b$ such that $\UL(\Inh_{\mc E\b}(\psi)) = \mc F$ and $\UL(\psi_{\Inh,\mc E\b}) = \psi$. 

 We proceed as follows: For $i \in \{1,\ldots, s\}$ we define  
 \begin{align*}
     \mc F_i^{(-l)} \xrightarrow{\tau_i^{(-l)}} \mc F_i^{(-l+1)} 
 \end{align*}
 by induction on $l \in \{1,\ldots, n_i\}$ such that $\mc F_i^{(0)} = \mc F$ and such that each of the squares  \\
 \begin{tikzcd}
\mathcal{F}_i^{(-l)}  \arrow[rr, "\tau_i^{(-l)}"] \arrow[d] &  & \mathcal{F}_i^{(-l+1)} \arrow[d] \\
\mc E_i^{(-l)} \arrow[rr, "\phi_i^{(-l)}"']                   &  & \mc E_i^{(-l+1)}                
\end{tikzcd} \\
 is a pullback square, where the right vertical map is $\psi$ in the first step of the induction and in each further induction step the pullback of the right vertical map constructed in the previous step.

 For this construction to make sense, we must check that \\
 \begin{tikzcd}
\mc F \otimes \mc O(-x_i) \arrow[rr] \arrow[d, "\psi \otimes \mc O(-x_i)"'] &  & \mc F  \arrow[d, "\psi"] \\
\UL(\mc E\b) \otimes \mc O(-x_i) \arrow[rr]                    &  & \UL(\mc E\b)                                                 
\end{tikzcd} \\
is cartesian, where the horizontal maps are induced by the natural map $\mc O(-x_i) \to \mc O$. Extending this diagram to the diagram with exact rows \\
\begin{tikzcd}
\mc F \otimes \mc O(-x_i) \arrow[rr] \arrow[d, "\psi \otimes \mc O(-x_i)"'] &  & \mc F \arrow[d, "\psi"] \arrow[rr] &  & \mc F \vert_{\{x_i\} \times T} \arrow[d] \arrow[r] & 0 \\
\UL(\mc E\b) \otimes \mc O(-x_i) \arrow[rr]                                 &  & \UL(\mc E\b) \arrow[rr]            &  & \UL(\mc E\b) \vert_{\{x_i\} \times T} \arrow[r]    & 0
\end{tikzcd} \\
we see that assumptions i) and ii) imply the existence part of the universal property of pullbacks and assumption ii) implies the uniqueness part. \\

\begin{Rem}
    Note that the above construction in particular applies to the situation where $\mc F$ is a vector bundle, $\mc E\b$ is a parabolic vector bundle and $\psi$ is an injection of vector bundles\footnote{Meaning a monomorphism of sheaves whose cokernel is a vector bundle.}. In this case, $\Inh_{\mc E\b}(\psi)$ is a parabolic vector bundle.
\end{Rem}

\subsection{Irrelevant parabolic bundles}
We define the notion of an \textit{irrelevant parabolic bundle}. This notion depends on the choice of a twisted parabolic datum and is inspired by the following situation: Let $G$ be an abstract group acting on a set $A$, let $\chi : G \to \mb C^\times$ be a character and let $f : A \to \mb C$ be a function, which is $\chi$-twisted equivariant, i.e. $f(ga) = \chi(g)f(a)$ holds for all $g \in G$ and $a \in A$. If for $a \in A$ the restriction of $\chi$ to the stabilizer of $a$ in $G$ is non-trivial, then $f(a) = 0$. In this sense the element $a$ is "irrelevant" for $f$.

Let $K$ be a field extension of $k$.

\begin{Def} \label{Main Definition of irrelevance recognizing}
    Let $D^\tw = (D,\Lambda) \in \TwParDat_{s,\ul n}$ and let $\mc E\b \in \Bun_{S^\oro,D}(K)$.
    \begin{itemize}
        \item[1)] Let $\mc F\b \in \Bun_{S^\oro,\ul n}(K)$ and $\iota : \mc F\b \to \mc E\b$ be an injection of parabolic vector bundles. $\iota$ is called $D^\tw$-\textit{irrelevance recognizing} iff $\iota$ admits a complement and $D(\mc F\b) \in \ParDat_{s,\ul n, < D}^{D^\tw-\irec}$.
        \item[2)] Let $\mc F \in \Bun_\emptyset$ and let $\iota : \mc F \to \UL(\mc E\b)$ be an injection of vector bundles. $\iota$ is called $D^\tw$-\textit{irrelevance recognizing} iff the induced morphism $\iota_{\Inh,\mc E\b} : \Inh_{\mc E\b}(\iota) \to \mc E\b$ is $D^\tw$-irrelevance recognizing.
        \item[3)] $\mc E\b$ is called $D^\tw$-\textit{irrelevant} (or $\Lambda$-\textit{irrelevant}) iff there exists a $D^\tw$-irrelevance recognizing injection of parabolic vector bundles $\iota : \mc F\b \to \mc E\b$.
    \end{itemize}
\end{Def}

Let $D^\tw = (D,\Lambda) \in \TwParDat_{s,\ul n}$. We denote by 
\begin{align*}
    \vert \Bun_D \vert^{\Lambda-\irrel}
\end{align*}
the subspace of the topological space $\vert \Bun_D \vert$ consisting of isomorphism classes of all $\mc E\b \in \Bun_D(K)$, $K$ an algebraically closed field extension of $k$, which are $D^\tw$-irrelevant.

\subsection{Ext groups and noncuspidality}
Let $K$ be a field extension of $k$. Since $\ParMod_{S^\oro,\ul n}(K)$ is a $K$-linear abelian category, we can form for $\mc F\b, \mc H\b \in \ParMod_{S^\oro,\ul n}(K)$  the group $\Ext^1(\mc H\b,\mc F\b)$ of extensions of $\mc H\b$ by $\mc F\b$ up to the usual notion of equivalence. Forgetting the parabolic structure defines a natural map
\begin{align*}
    \eta_{\mc H\b,\mc F\b} : \Ext^1(\mc H\b,\mc F\b) \to \Ext^1(\UL(\mc H\b),\UL(\mc F\b)).
\end{align*}

\begin{Lem}
    If $\mc H\b,\mc F\b$ are parabolic vector bundles, the morphism $\eta_{\mc H\b,\mc F\b}$ is surjective.
\end{Lem}
\begin{proof}
    Consider a representative    \begin{align*}
        0 \to \UL(\mc F\b) \xrightarrow{j} \mc E \xrightarrow{p} \UL(\mc H\b) \to 0
    \end{align*}
    of an element in $\Ext^1(\UL(\mc H\b),\UL(\mc F\b))$. Assume the parabolic structure $\mc F\b$ resp. $\mc H\b$ is given at the point $x_i$ by flags
    \begin{align*}
        0 &= F_0 \subset  F_1 \subset \ldots \subset  F_{n_i} = \UL(\mc F\b) \otimes k(x_i) \\
        0 &=  H_0 \subset H_1 \subset \ldots \subset H_{n_i} = \UL(\mc H\b) \otimes k(x_i).
    \end{align*}
    For each $i \in \{1,\ldots, s\}$ pick a splitting $f_i : \UL(\mc H\b) \otimes k(x_i) \to \mc E \otimes k(x_i)$ of the surjective morpism of $k(x_i)$-vector spaces $p \otimes k(x_i)$. Then define for $l \in \{1,\ldots, n_i\}$
    \begin{align*}
        \mc E_i := (j \otimes k(x_i))(F_l) + f_i(H_l).
    \end{align*}
     Doing this for each point $x_i$, we obtain a parabolic bundle $\mc E\b$ with underlying bundle $\mc E$ and it is easy to check that $\mc E\b$ sits in an exact sequence of parabolic bundles
     \begin{align*}
        0 \to \mc F\b \xrightarrow{j\p} \mc E\b \xrightarrow{p\p} \mc H\b \to 0
    \end{align*}
    such that $\UL(j\p) = j$ and $\UL(p\p) = p$.
\end{proof}

\begin{Def} \label{Main Definition noncupsidality recognizing}
Let $D^\tw = (D,\Lambda) \in \TwParDat_{s,\ul n}$ and let $\mc E\b \in \Bun_{S^\oro,\ul n, D}(K)$. 
 \begin{itemize}
        \item[1)] Let $\mc F\b \in \Bun_{S^\oro,\ul n}(K)$ and $\iota : \mc F\b \to \mc E\b$ be an injection of parabolic vector bundles. $\iota$ is called $D^\tw$-\textit{noncuspidality recognizing} iff
        \begin{itemize}
            \item[i)] $\mf R(\Res_{D^\tw}(D( \mc F\b))) = 1$,
            \item[ii)] $\iota$ admits a complement, i.e. there is a parabolic vector bundle $\mc H\b$ such that
            \begin{align*}
                0 \to \mc F\b \xrightarrow{\iota} \mc E\b \to \mc H\b \to 0
            \end{align*}
            is split,
            \item[iii)] $\Ext^1(\mc H\b,\mc F\b) = 0$.
        \end{itemize}
        \item[2)] Let $\mc F \in \Bun_\emptyset(K)$ and let $\iota : \mc F \to \UL(\mc E\b)$ be an injection of vector bundles. $\iota$ is called  $D^\tw$-\textit{noncuspidality recognizing} iff the induced morphism $\iota_{\Inh,\mc E\b} : \Inh_{\mc E\b}(\iota) \to \mc E\b$ is $D^\tw$-noncuspidality recognizing.
        \item[3)] $\mc E\b$ is called $D^\tw$-\textit{noncuspidal} (or $\Lambda$-\textit{noncuspidal}) iff there exists a $D^\tw$-noncuspidality recognizing injection of parabolic vector bundles $\iota : \mc F\b \to \mc E\b$.
    \end{itemize}
\end{Def}

Let $D^\tw = (D,\Lambda) \in \TwParDat_{s,\ul n}$. We denote by 
\begin{align*}
    \vert \Bun_D \vert^{\Lambda-\noncusp}
\end{align*}
the subspace of the topological space $\vert \Bun_D \vert$ consisting of isomorphism classes of all $\mc E\b \in \Bun_D(K)$, $K$ an algebraically closed field extension of $k$, which are $D^\tw$-noncuspidal.

\subsection{Indecomposable bundles and the minimal locus}
Let $D \in \ParDat_{s, \ul n}$. 

\begin{Def}
    Let $K$ be an algebraically closed field extension of $k$. A parabolic vector bundle $\mc E\b \in \Bun_D(K)$ is called \textit{indecomposable} iff it is an indecomposable object in the abelian category $\ParMod_{X,S^\oro,\ul n}(K)$. We denote the full subcategory of $\Bun_D(K)$ consisting of indecomposable bundles by $\Bun_D(K)^\inde$.
\end{Def}

\begin{Lem} \label{Lemma characterization indecomposable}
    For $\mc E\b \in \Bun_D(K)$ the following are equivalent:
    \begin{itemize}
        \item[1)] $\mc E\b \in \Bun_D(K)^\inde$
        \item[2)] The only idempotents in $\End(\mc E\b)$ are $0$ and $1$.
    \end{itemize}
\end{Lem}
\begin{proof}
    Such a statement holds in any abelian category.
\end{proof}

\begin{Def}
    We let $\vert \Bun_D \vert^\inde$ be the set consisting of those points $P \in \vert \Bun_D \vert$ such that there exists an algebraically closed field extension $K$ of $k$ and a representative $\mc E\b \in \Bun_D(K)$ of $P$ such that $\mc E\b \in \Bun_D(K)^\inde$.  
\end{Def}

The set $\vert \Bun_D \vert^\inde$ need not be closed inside $\vert \Bun_D \vert$ (see Proposition \ref{Long Proposition about intereseting locus in degree one compoenent} for an explicit example). 

One can characterize indecomposability of a parabolic vector bundle via the existence of non-scalar semisimple elements in its automorphism group:

\begin{Lem} \label{Lemma criterion semisimple}
    Let $K$ be an algebraically closed field extension of $k$ and $\mc E\b \in \Bun_D(K)$. The following are equivalent:
    \begin{itemize}
        \item[1)] $\mc E\b \in \Bun_D(K)^\inde$.
        \item[2)] $\End(\mc E\b)$ does not contain non-scalar semisimple elements.
        \item[3)] $\Aut(\mc E\b)$ does not contain non-scalar semisimple elements.
        \item[4)] There does not exist a commutative diagram \\
        \begin{tikzcd}
{\mb G_{m,K}^{\times 2}} \arrow[rr, "f", hook] &                                                                      & \ul{\Aut}(\mc E\b) \\
                                               & {\mb G_{m,K}} \arrow[lu, "\Delta", hook] \arrow[ru, "\mf z"', hook'] &                   
\end{tikzcd} \\
where $\Delta$ is the diagonal, $\mf z$ is the natural inclusion (given by scalar multiplication) and $f$ is a monomorphism.
\end{itemize}
\end{Lem}
\begin{proof}
    Clearly, the negation of 1) implies the negation of 4). Also the negation of 4) implies the negation of 3), since by commutativity of the diagram and injectivity of $f$, $f(\mb G_m^{\times 2}(K))$ contains non-scalar elements (here we also use $\#K > 2$) and $K$-points of tori are semisimple. The negation of 3) trivially implies the negation of 2). \\
    To see that the negation of 2) implies the negation of 1), let $\phi \in \End(\mc E\b)$ be a non-scalar semisimple element. This means that its minimal polynomial $p_\phi$ can be written as 
    \begin{align*}
        p_\phi = \prod_{i=1}^n (X - \lambda_i)
    \end{align*}
    for $\lambda_1, \ldots, \lambda_n \in K$ pairwise distinct. Since $\phi$ is non-scalar, $n \ge 2$. Hence $q := \prod_{i=2}^n (X - \lambda_i)$ has degree $\ge 1$. Since the $\lambda_i$ are pairwise distinct, we find $\alpha, \beta \in K[X]$ such that $\alpha (X-\lambda_1) + \beta q = 1$. Set $e_1 := \alpha(\phi)(\phi - \lambda_1)$ and $e_2 := \beta(\phi) q(\phi)$. Then $e_1 + e_2 = 1$ and, since $\alpha(\phi)$, $(\phi - \lambda_1)$, $\beta(\phi)$ and $q(\phi)$ lie in the commutative subring $K[\phi]$ of $\End(\mc E\b)$, also $e_1e_2 = e_2e_1 = \alpha(\phi) \beta(\phi)(\phi - \lambda_1)q(\phi) = \alpha(\phi) \beta(\phi)p_\phi(\phi) = 0$. We obtain $e_1^2 - e_1 = (1-e_2)e_1 - e_1 = - e_1e_2 = 0$. Hence $e_1$ is idempotent and $\neq 1$ by construction. The same holds for $e_2$. Hence $e_1 \neq 0,1$ and so $\mc E\b \notin \Bun_D(K)^\inde$ by Lemma \ref{Lemma characterization indecomposable}.
\end{proof}

Recall that for any algebraic stack $\mf X$ locally of finite type over $k$ and for any $n \in \mb Z_{\ge 0} \cup \{\infty\}$ the set 
\begin{align*}
    \{x \in \vert \mf X \vert \ \vert \  \dim(\ul{\Aut}(x)) \le n\}
\end{align*}
is an open subset of $\vert \mf X \vert$. For $n = \mindim(\mf X)$, where
\begin{align*}
    \mindim(\mf X) := \min\{n\p \in \mb Z_{\ge 0} \cup \{\infty\} \ \vert \ \exists x \in \mf X \ \vert \text{ such that } \dim(\ul{\Aut}(x)) = n\p\} \in \mb Z_{\ge 0} \cup \{\infty\}
\end{align*}
we denote the corresponding open substack of $\mf X$ by $\mf X^\mini$.

\begin{Cor} \label{minimal locus inside interesting locus}
The following hold:
\begin{itemize}
    \item[i)] $\mindim(\Bun_D) \ge 1$.
    \item[ii)] If $\mindim(\Bun_D) = 1$, then $\vert \Bun_D^\mini \vert \subset \vert \Bun_D \vert^\inde$.
\end{itemize}
\end{Cor}
\begin{proof}
    $\mb G_m$ is contained in the automorphism group of every point of $\Bun_D$, hence i) is true. ii) follows from Lemma \ref{Lemma criterion semisimple}, part 4).
\end{proof}

We will see an example where the inclusion in ii) is not an equality (see Theorem \ref{Classification Theorem for degree zero}).

\subsection{The interesting locus}
Let $D^\tw = (D,\Lambda) \in \TwParDat_{s,\ul n}$. 

\begin{Def} \label{Definition of interesting locus}
    The $D^\tw$\textit{-interesting locus} is
\begin{align*}
    \vert \Bun_D \vert^{\Lambda-\inter} := \vert \Bun_D \vert \setminus \l \vert \Bun_D \vert^{\Lambda-\noncusp} \cup \vert \Bun_D \vert^{\Lambda-\irrel} \r.
\end{align*}
\end{Def}

If $D$ or even $D^\tw$ is clear from the context, we will call $\vert \Bun_D \vert^{\Lambda-\inter}$ the $\Lambda$\textit{-interesting} or just \textit{interesting locus}. Note that it always contains $\vert \Bun_D \vert^\inde$, hence if $\mindim(\Bun_D) = 1$ holds, also $\vert \Bun_D^\mini \vert$. In general, it will not be locally closed inside $\Bun_D$ (see Proposition \ref{Long Proposition about intereseting locus in degree one compoenent} for an explicit example).

We find this locus interesting for the following reason: With $D$ one can naturally associate an algebraic group $L$ and an $L$-torsor $\pi_D^+ : \Bun_D^+ \to \Bun_D$. With $\Lambda$, one can associate a character sheaf $\chi_\Lambda$ on $L$ and consider the derived category of $\chi_\Lambda$-equivariant sheaves on the total space of the torsor $\pi_D^+$. Then the Langlands philosophy predicts that the image under $\pi_D^+$ of the support\footnote{By \textit{support} we mean literally the set of points $P^+ \in \vert \Bun_D^+ \vert$ such that the stalk of the given object in the derived category at $P^+$ vanishes (and not the closure of this set).} of cuspidal sheaves lies inside $\vert \Bun_D \vert^\inter$.

\subsection{Elementary Hecke operator}
Let $T$ be a $k$-scheme and $X$, $s$, $S^\oro$, $\ul n$ as before. We now construct an additive functor
\begin{align*}
    \mb T^{(l), \ParBun,T}_s : \ParBun_{S^\oro,\ul n}(T) \to \ParBun_{S^\oro,\ul n}(T), 
\end{align*}
which should be thought of as the \textit{elementary degree-one lower Hecke operator} at $x_s$. Let $(\mc F\b,\phi\b) \in \ParBun_{S^\oro,\ul n}(T)$. Write $\ul n = (n_1, \ldots, n_s)$. We set
\begin{align*}
    \mc H_s^{(-j)} &:= \mc F_s^{(-j-1)}, \ j \in \{0,\ldots, n_s-1\}, \\
    \psi_s^{(-j)} &:= \phi_s^{(-j-1)}, \ j \in \{1,\ldots, n_s-1\}
\end{align*}
and define
\begin{align*}
    \psi_s^{(-n_s)} := \phi_s^{(-1)} \otimes \mc O(-x_s) : \mc F_s^{(-1)} \otimes \mc O(-x_s) \to \mc F_s^{(0)} \otimes \mc O(-x_s).
\end{align*}
For $i \in \{1,\ldots, s-1\}$ we set 
\begin{align*}
    \mc H_i^{(0)} &:= \mc F_s^{(-1)} \\
    \Psi_i^{(0)} &:= \l\mc H_i^{(0)} = \mc F_s^{(-1)} \xrightarrow{\phi_s^{(-1)}} F_s^{(0)}\r.
\end{align*}
Then we define 
\begin{align*}
    \mc H_i^{(-j)}, \\
    \psi_i^{(-j)} : \mc H_i^{(-j)} \to \mc H_i^{(-j+1)}, \\
    \Psi_i^{(-j)} : \mc H_i^{(-j)} \to \mc F_i^{(-j)}
\end{align*}
for $j \in \{1,\ldots, n_i\}$ inductively such that the squares 
\[
\begin{tikzcd}
\mc H_i^{(-j)} \arrow[d, "\Psi_i^{(-j)}"'] \arrow[rr, "\psi_i^{(-j)}"] &  & \mc H_i^{(-j+1)} \arrow[d, "\Psi_i^{(-j+1)}"] \\
\mc F_i^{(-j)} \arrow[rr, "\phi_i^{(-j)}"']                             &  & \mc F_i^{(-j+1)}                            
\end{tikzcd} 
\]
are cartesian. Since concatenations of cartesian diagrams are cartesian, $\mc H_i^{(0)}$ fits into a cartesian square 
\[
\begin{tikzcd}
\mc H_i^{(-n_i)} \arrow[d, "\Psi_i^{(-n_i)}"'] \arrow[rr, "\psi_i^{(-1)} \circ \ldots \circ \psi_i^{(-n_i)}"] &  & \mc H_i^{(0)} \arrow[d, "\Psi_i^{(0)}"] \\
\mc F_i^{(-n_i)} \arrow[rr, "\phi_i^{(-1)} \circ \ldots \circ \phi_i^{(-n_i)}"']                           &  & \mc F_i^{(0)}.                
\end{tikzcd}
\]
From the definition of parabolic bundles we obtain
\begin{align*}
    \phi_i^{(-1)} \circ \ldots \circ \phi_i^{(-n_i)} = \id \otimes \iota : \mc F_i^{(0)}(-x_i) \to \mc F_i^{(0)}.
\end{align*}
Then Lemma \ref{Lemma about cartesian diagram} below implies that the above diagram is naturally isomorphic to 
\[
\begin{tikzcd}
\mc H_i^{(0)}(-x_i) \arrow[d, "\Psi_i^{(0)} \otimes \id"'] \arrow[rr, "\id \otimes \iota_{x_i}"] &  & \mc H_i^{(0)} \arrow[d, "\Psi_i^{(0)}"] \\
\mc F_i^{(0)}(-x_i) \arrow[rr, "\id \otimes \iota_{x_i}"']                                         &  & \mc F_i^{(0)}.                            
\end{tikzcd}
\]
Hence, 
\begin{align*}
    \mc H\b &:= (\mc H_i^{(j)}, (\psi_i^{(j)})_{ij}) \in \ParBun_{S^\oro,\ul n}(T)  
\end{align*}
and we set
\begin{align*}
    \mb T^{(l), \ParBun,T}_s(\mc F\b) := \mc H\b.
\end{align*}

\begin{Lem} \label{Lemma about cartesian diagram}
    Let $Y$ be a scheme and $Z \subset Y$ an effective Cartier divisor. Let $\iota : \mc O_Y(-Z) \to \mc O_Y$ be the natural inclusion. Let $\mc H, \mc F$ be locally free sheaves of finite rank on $Y$. For any morphism of $\mc O_Y$-modules $\phi : \mc H \to \mc F$ such that the restriction $\phi \vert_Z : \mc H \vert_Z \to \mc F \vert_Z$ is injective the diagram 
    \[
    \begin{tikzcd}
\mc H(-Z) \arrow[d, "\phi \otimes \id"'] \arrow[r, "\id \otimes \iota"] & \mc H \arrow[d, "\phi"] \\
\mc F(-Z) \arrow[r, "\id \otimes \iota"']                               & \mc F                  
\end{tikzcd} 
\]
is cartesian.    
\end{Lem}
\begin{proof}
   Consider the following diagram with exact rows, where the vertical morphisms are induced by $\phi$: \\
   \begin{tikzcd}
0 \arrow[r] & \mc F(-Z) \arrow[r, "\ell"] \arrow[d, "\phi \otimes \mc O(-Z)"'] & \mc F \arrow[d, "\phi"'] \arrow[r, "\pi"] & (\iota_Z)_*(\mc F \vert_Z) \arrow[d, "(\iota_Z)_*(\phi \vert_Z)"] \arrow[r] & 0 \\
0 \arrow[r] & \mc H(-Z) \arrow[r, "\rho"']                                     & \mc H \arrow[r]                           & (\iota_Z)_*(\mc H \vert_Z) \arrow[r]                                        & 0
\end{tikzcd} \\
We want to show that the left square is cartesian. So let $\mc J$ be any $\mc O_Y$-module $Y$ together with morphisms $\alpha_1 : \mc J \to \mc F$ and $\alpha_2 : \mc J \to \mc H(-Z)$ such that $\phi \alpha_1 = \rho \alpha_2$. Then $\pi \alpha_1 = 0$. Indeed, by assumption $(\iota_Z)_*(\phi \vert_Z)$ is injective, hence $\pi \alpha_1 = 0$ is implied by $(\iota_Z)_*(\phi \vert_Z) \pi \alpha_1 = 0$, which follows from commutativity of the diagram and exactness of the lower row. Hence, $\alpha_1$ factors through $\alpha : \mc J \to \mc F(-Z)$.

It is left to show that $(\phi \otimes \mc O(-Z)) \circ \alpha = \alpha_2$. Since $\rho$ is a monomorphism, it suffices to show $\rho \circ (\phi \otimes \mc O(-Z)) \circ \alpha = \rho \circ \alpha_2$. This follows from $\rho \circ (\phi \otimes \mc O(-Z)) = \phi \circ \ell$, $\ell \alpha = \alpha_1$ and $\phi \alpha_1 = \rho \alpha_2$.
\end{proof}

Now we define $\mb T^{(l), \ParBun,T}_s$ on morphisms: Let $f = (f_i\b)_i : \mc F\b \to \tilde{\mc F}\b$ a morphism in $\ParBun_{S^\oro,\ul n}(T)$ and denote $\mb T^{(l), \ParBun,T}_s(\mc F\b)$ resp. $\mb T^{(l), \ParBun,T}_s(\tilde{\mc F}\b)$ by $\mc H\b$ resp. $\tilde{\mc H}\b$. Define
\begin{align*}
    g_s^{(-j)} &:= f_s^{(-j-1)}, j \in \{1,\ldots, n_s-2\}, \\
    g_s^{(-n_s-1)} &:= f_s^{(-1)} \otimes \mc O(-x_s).
\end{align*}
For $i \in \{1,\ldots, s\}$ define
\begin{align*}
    g_i^{(0)} := f_s^{(-1)} : \mc H_i^{(0)} = \mc F_s^{(-1)} \to \tilde{\mc F}_s^{(-1)} = \tilde{\mc H}_i^{(0)}.
\end{align*}
Then we define for fixed $i$ 
\begin{align*}
    g_i^{(-j)} : \mc H_i^{(-j)} \to \tilde{\mc H}_i^{(-j)}
\end{align*}
for $j \in \{1,\ldots, n_i-1\}$ inductively as the dotted arrow making the whole diagram 
\[
\begin{tikzcd}[row sep=large, column sep=large]
& \tilde{\mc H}_i^{(-j)}
    \arrow[rr]
    \arrow[dd]
& &
\tilde{\mc H}_i^{(-j+1)}
    \arrow[dd]
\\
\mc H_i^{(-j)}
    \arrow[rr]
    \arrow[dd]
    \arrow[ur, dotted]
& &
\mc H_i^{(-j+1)}
    \arrow[dd]
    \arrow[ur]
\\
& \tilde{\mc F}_i^{(-j)}
    \arrow[rr]
& &
\tilde{\mc F}_i^{(-j+1)}
\\
\mc F_i^{(-j)}
    \arrow[rr]
    \arrow[ur]
& &
\mc F_i^{(-j+1)}
    \arrow[ur]
\end{tikzcd} 
\]
commutative. Note that there is a unique such dotted arrow. Indeed, the back face of the cube is cartesian by construction of $\mb T^{(l), \ParBun,T}_s(\tilde{\mc F}\b)$. Concatenating all cubes as above for $j \in \{1,\ldots, n_i-1\}$ and using the fact that in each cube the back face is cartesian, one sees that the collection $(g_i^{(j)})_{ij}$ defines a morphism of parabolic $\mc O$-modules. \\

Using the universal property of pullbacks one shows that the above definitions are functorial, i.e. we obtain the desired functor $\mb T^{(l), \ParBun,T}_s$, and that this functor is additive. \\

Using the index shift $j \to j+1$ instead $j \to j-1$ at $x_s$ one defines the \textit{elementary  degree one upper Hecke operator} 
\begin{align*}
    \mb T^{(u), \ParBun,T}_s : \ParBun_{S^\oro,\ul n}(T) \to \ParBun_{S^\oro,\ul n}(T), 
\end{align*}
at $x_s$ analogously.

\begin{Lem} \label{Lemma Hecke operators inverse}
    $\mb T^{(l), \ParBun,T}_s$ and $\mb T^{(u), \ParBun,T}_s$ are mutually inverse additive equivalences of categories.
\end{Lem}
\begin{proof}
    By the above discussion we only have to show that  $\mb T_s^{(l)}$ and $\mb T_s^{(u)}$ are mutually inverse to each other. This follows from Lemma \ref{Technical Lemma about cartesian vs cocartesian} below.
\end{proof}

\begin{Lem} \label{Technical Lemma about cartesian vs cocartesian}
    Let $Y$ be any scheme. Consider a commutative diagram \\
    \begin{tikzcd}
\mc F \arrow[r, "f_2"] \arrow[d, "f_1"'] & \mc H_1 \arrow[d, "k_1"] \\
\mc H_2 \arrow[r, "k_2"']                & \mc K                   
\end{tikzcd} \\
in $\QCoh(Y)$ with $\mc H_i$ locally free. The following are equivalent:
\begin{itemize}
    \item[1)] Each of the following hold:
            \begin{itemize}
                \item[i)] $\mc K$ is locally free.
                 \item[ii)]  $k_1$ and $k_2$ are injective.
                \item[iii)] $\supp(\coker(k_1)) \cap \supp(\coker(k_2)) = \emptyset$
                \item[iv)] The diagram is cartesian.
            \end{itemize}
    \item[2)] Each of the following hold:
            \begin{itemize}
                \item[i)] $\mc F$ is locally free.
                 \item[ii)]  $f_1$ and $f_2$ are injective.
                \item[iii)] $\supp(\coker(f_1)) \cap \supp(\coker(f_2)) = \emptyset$
                \item[iv)] The diagram is cocartesian.
            \end{itemize}
\end{itemize}
If 1) or 2) hold, then all sheaves in the diagram are locally free, all morphisms injective and we have $\supp(\coker(f_i)) = \supp(\coker(k_i))$, $i = 1,2$.
\end{Lem}
\begin{proof}
    Note that all statements are local. By conditions i) and iii) we may reduce to a situation in which one of the maps $k_i$ in situation 1) resp. $f_i$ in situation 2) is an isomorphism. Then the statement is clear.
\end{proof}

Fix $\ast \in \{u,l\}$. The collection $\{\mb T^{(\ast), \ParBun,T}_s\}_T$ is functorial in $T$. Hence it induces a morphism of stacks
\begin{align*}
    \mb T^{(\ast)}_s : \Bun_{S^\oro,\ul n} \to \Bun_{S^\oro,\ul n}.
\end{align*}
Lemma \ref{Lemma Hecke operators inverse} implies:

\begin{Cor}
    $\mb T^{(l)}_s$ and $\mb T^{(u)}_s$ are mutually inverse isomorphisms of stacks. 
\end{Cor}

For $D = (D_1,\ldots, D_s) \in \ParDat_{s,\ul n}$ write $D_s = (d_1,\ldots, d_{n_s})$ and set
\begin{align*}
    D_s^{(l)} &:= (d_{n_s}, d_1, \ldots, d_{n_s-1}) \\
    D_s^{(u)} &:= (d_2, \ldots, d_{n_s}, d_1).
\end{align*}
Then define for $\ast \in \{u,l\}$
\begin{align*}
    D^{(\ast)} := (D_1,\ldots, D_{s-1}, D_s^{(\ast)}).
\end{align*}

\begin{Lem} \label{Lemma behavior of Hecke operators}
    Let $d \in \mb Z$ and $D \in \ParDat_{s,\ul n}$. \\
    $\mb T^{(l)}_s$ maps the open substack $\Bun_{S^\oro,\ul n}^d$ of $\Bun_{S^\oro,\ul n}$ to $\Bun_{S^\oro,\ul n}^{d-d_{n_s}}$ and $\Bun_{S^\oro,D}$ to $\Bun_{S^\oro,D^{(l)}}$. \\
    Similarly, $\mb T^{(u)}_s$ maps $\Bun_{S^\oro,\ul n}^d$ to $\Bun_{S^\oro,\ul n}^{d+d_1}$ and $\Bun_{S^\oro,D}$ to $\Bun_{S^\oro,D^{(u)}}$. \\
    In particular, $\mb T^{(l)}_s$ induces an isomorphism
    \begin{align*}
        \Bun_{S^\oro,D}^d \simeq \Bun_{S^\oro,D^{(l)}}^{d-d_{n_s}}.
    \end{align*}
\end{Lem}
\begin{proof}
    This follows immediately from the definitions. 
\end{proof}

\begin{Cor} \label{Corollary isomorphisms different degree components for special D}
    Let $d \in \mb Z$ and $D = (D_1,\ldots, D_s) \in \ParDat_{s,\ul n}$. Assume that $D_s = (1,\ldots, 1) \in \mb Z^{n_s}$. Then $\mb T^{(l)}_s$ induces an isomorphism
    \begin{align*}
        \Bun_{S^\oro,D}^d \simeq \Bun_{S^\oro,D}^{d-1}.
    \end{align*}
\end{Cor}
\begin{proof}
    Use that $D_s = (1,\ldots, 1)$ implies $D^{(l)} = D$ and Lemma \ref{Lemma behavior of Hecke operators}.
\end{proof}

We obtain:

\begin{Prop}
     Let $d \in \mb Z$ and $D = (D_1,\ldots, D_s) \in \ParDat_{s,\ul n}$. Assume that $D_s = (1,\ldots, 1) \in \mb Z^{n_s}$. Then $\mb T^{(l)}_s$ induces an isomorphism
    \begin{align*}
        \Bun_{S^\oro,D}^{d,\mini} \simeq \Bun_{S^\oro,D}^{d-1,\mini}
    \end{align*}
    and a homeomorphism
    \begin{align*}
         \vert \Bun_{S^\oro,\ul n, D}^d \vert^\inde \simeq \vert \Bun_{S^\oro,\ul n, D}^{d-1}\vert^\inde.
    \end{align*}
\end{Prop}
\begin{proof}
    The isomorphism follows from Corollary \ref{Corollary isomorphisms different degree components for special D} and the fact that an isomorphism of stack induces an isomorphism between the automorphism group of a point in the source and the automorphism group of its image in the target. \\
    The homeomorphism follows from Corollary \ref{Corollary isomorphisms different degree components for special D} and the additivity of the functor $\mb T^{(\ast), \ParBun,T}_s$, $\ast \in \{u,l\}$, for $T$ the spectrum of an algebraically closed field over $k$.
\end{proof}

\begin{Theo} \label{Proposition interesting locus and hecke operators}
    Let $D^\tw =  (D,\Lambda) = (D_1^\tw,\ldots, D_s^\tw) \in \TwParDat_{s,\ul n}$ and $\lambda \in \mb C^{\times}$. Assume that \\
    $D_s^\tw = ((1,\ldots, 1),(\lambda,\ldots, \lambda))$. Let $K \supset k$ be an algebraically closed field and let $\mc E\b \in \Bun_{S^\oro,\ul n, D}(K)$. 
    \begin{itemize}
        \item[1)] If $\mc E\b$ is $D^\tw$-irrelevant, so are $\mb T^{(l)}_s(\mc E\b)$ and $\mb T^{(u)}_s(\mc E\b)$.
        \item[2)] If $\mc E\b$ is $D^\tw$-noncuspidal, so are $\mb T^{(l)}_s(\mc E\b)$ and $\mb T^{(u)}_s(\mc E\b)$.
    \end{itemize}
    In particular $\mb T^{(l)}_s$ induces equivalences of categories 
    \begin{align*}
        \Bun_{S^\oro,D}^d(K)^\irrel &\simeq \Bun_{S^\oro,D}^{d-1}(K)^\irrel \\
         \Bun_{S^\oro,D}^d(K)^\noncusp &\simeq \Bun_{S^\oro,D}^{d-1}(K)^\noncusp \\
         \Bun_{S^\oro,D}^d(K)^\inter &\simeq \Bun_{S^\oro,D}^{d-1}(K)^\inter
    \end{align*}
    and a homeomorphism
    \begin{align*}
        \vert \Bun_{S^\oro,D}^d \vert^\inter &\simeq \vert \Bun_{S^\oro,D}^{d-1} \vert^\inter.
    \end{align*}
\end{Theo}
\begin{proof}
    Write $D = (D_1,\ldots, D_s)$. Note that $D_s = (1,\ldots, 1)$ implies $D(\mb T^{(l)}_s(\mc E\b)) = D(\mc E\b)$, so the statements make sense. \\
    Assume $\mc E\b$ is $D^\tw$-irrelevant. Let $\iota : \mc F \hookrightarrow \mc E\b$ be a $D^\tw$-irrelevance recognizing injection of parabolic vector bundles (defined over $K$). We claim that $\iota\p := \mb T^{(l),\ParBun,\Spec(K)}(\iota)$ is again a $D^\tw$-irrelevance-recognizing injection of parabolic vector bundles. Lemma \ref{Lemma Hecke operators inverse} implies that $\iota\p$ admits a complement. Hence, it suffices to check that $\mf R(\Res_{D^\tw}(D(T^{(l)}(\mc F\b)))) \neq 1$. We claim that $\mf R(\Res_{D^\tw}(D(T^{(l)}(\mc F\b)))) = \mf R(\Res_{D^\tw}(D(\mc F\b)))$, which is $\neq 1$ by assumption. \\
    To prove the claim, write $D(\mc F) = (\tilde D_1,\ldots, \tilde D_s)$, $\tilde r = \rk(\mc F)$, $r = \rk(\mc E)$. Write $\tilde D_i = (d_{ij})_j$ and $\Lambda = (\lambda_{ij})_{ij}$, so that $\lambda_{s,j} = \lambda$ for all $j \in \{1,\ldots, n_s\}$. Then $D(\mc F) < D$ implies $\tilde D_s = (l_1,\ldots, l_r)$ with $l_i \in \{0,1\}$ and $\#\{i \ \vert \ l_i = 1\} = \tilde r$. Using $\tilde D_s^{(l)} = (l_r,l_1,\ldots, l_{r-1}) =: (l_1\p,\ldots, l_r\p)$ we calculate
    \begin{align*}
        \mf R(\Res_{D^\tw}(D(T^{(l)}(\mc F\b)))) &= \mf R(\Res_{D^\tw}(\tilde D_1,\ldots, \tilde D_{s-1}, \tilde D_s^{(l)})) \\
        &= \l \prod_{i=1}^{s-1} \prod_{j=1}^{n_i} \lambda_{ij}^{d_{ij}} \r \l \prod_{j=1}^{n_s} \lambda_{sj}^{l_j\p} \r \\
         &= \l \prod_{i=1}^{s-1} \prod_{j=1}^{n_i} \lambda_{ij}^{d_{ij}} \r \l \prod_{j=1}^{n_s} \lambda^{l_j\p} \r \\
        &= \l \prod_{i=1}^{s-1} \prod_{j=1}^{n_i} \lambda_{ij}^{d_{ij}} \r \lambda^{\tilde r} \\
        &= \l \prod_{i=1}^{s-1} \prod_{j=1}^{n_i} \lambda_{ij}^{d_{ij}} \r \l \prod_{j=1}^{n_s} \lambda_{sj}^{l_j} \r \\
        &= \mf R(\Res_{D^\tw}(\tilde D_1,\ldots, \tilde D_{s-1}, \tilde D_s)) \\
        &= \mf R(\Res_{D^\tw}(D(\mc F\b))).
    \end{align*}
    This proves 1) for $\mb T^{(l)}$. The case of $\mb T^{(u)}$ is similar. \\
    To show 2), assume $\mc E\b$ is $D^\tw$-noncuspidal and let $\iota : \mc F \hookrightarrow \mc E\b$ be a $D^\tw$-noncuspidality recognizing injection of parabolic vector bundles (defined over $K$). We claim that $\iota\p := \mb T^{(l),\ParBun,\Spec(K)}(\iota)$ is again a $D^\tw$-noncuspidality recognizing injection of parabolic vector bundles. Let $\mc H\b$ denote the complement of $\iota$, i.e. $\iota$ fits into a split exact sequence
     \begin{align*}
                0 \to \mc F\b \xrightarrow{\iota} \mc E\b \to \mc H\b \to 0.
    \end{align*}
    By the same arguments as in the proof of 1) it suffices to show that
    \begin{align*}
        \Ext^1(\mb T^{(l)}(\mc H\b), \mb T^{(l)}(\mc F\b)) = 0. 
    \end{align*} 
    Indeed, exactness of the additive equivalence $\mb T^{(l), \ParBun,T}_s$ (see Lemma \ref{Lemma Hecke operators inverse}) implies
    \begin{align*}
        \Ext^1(\mb T^{(l)}(\mc H\b), \mb T^{(l)}(\mc F\b)) \simeq \Ext^1(\mc H\b,\mc F\b),
    \end{align*}
    which vanishes, since $\iota$ is $D^\tw$-noncuspidality recognizing. \\
    This shows 2) for $\mb T^{(l)}$. The case of $\mb T^{(u)}$ is similar.

\end{proof}

%% file: 4_Interesting_Locus_Classification.tex
\section{Classification results: The example of \texorpdfstring{$\mb P^1$}{} minus 4 points} \label{section on classification}
In this section we study the following example: $X = \mb P^1$, $S^\oro = (0,1,t,\infty)$, where $t \in \mb P^1(k) \setminus \{0,1,\infty\}$, $D^\tw = (D,\Lambda)$ with 
\begin{align} \label{Definition of Data in Example}
    D &:= ((1,2),(1,2),(1,2),(1,1,1)), \\
    \Lambda &:= ((\lambda,1),(\lambda,1),(\lambda, 1), (\lambda^{-1},\lambda^{-1},\lambda^{-1})),
\end{align}
where $\lambda \in \mb C^\times \setminus \{\pm 1\}$.
We will analyze the interesting locus $\Bun_D(k)^\inter$ (note that this notion depends also on $\Lambda$) in the sense that we give a complete classification of points in $\Bun_D(k)^\inter$ and $\Bun_D^\mini(k)$. 

In the whole chapter, $D$ and $\Lambda$ are the concrete data in (\ref{Definition of Data in Example}). We classify points in points in $\Bun_D(k)^\inter$ by explicitly exhibiting all points which lie in $\vert \Bun_D(k) \vert \setminus \vert \Bun_D(k) \vert^\inter$. A point in this set is either $\Lambda$-irrelevant or $\Lambda$-noncuspidal. To verify that a point is $\Lambda$-irrelevant is easier than to check that it is $\Lambda$-noncuspidal, so we first exclude as many $\Lambda$-irrelevant points as we can.

The main result of this section is the classification result Theorem \ref{Classification Theorem for degree zero} of the degree $0$ interesting locus $\vert \Bun_D^0(k) \vert^\inter$ over $k$. We also provide an analogous theorem for the degree $1$ component in Theorem \ref{Classification Theorem for degree one}. The reason for doing this is explained at the end of subsection \ref{subsection Preparation}.

\subsection{Preparations} \label{subsection Preparation} In this subsection we introduce notation that will be used in carrying out the above steps.

Set $\Pi_i := \mb P^2$ for $i \in \{1,2,3\}$ and $\Pi_4 := \FL_3$. For $i \in \{1,2,3,4\}$ set
\begin{align*}
    \Pi_{\le i} &:= \prod_{j=1}^i \Pi_j \\
    \Pi_{> i} &:= \prod_{j=i+1}^4 \Pi_j.
\end{align*}
To shorten notation we set $\Pi := \Pi_{\le 4}$. 

We choose a line bundle $\mc O(1)$ on $\mb P^1$ of degree $1$ and set $\mc O(d) := \mc O(1)^{\otimes d}$ for $d \in \mb Z$. We identify $\GlobDat^3$ with the set
 \begin{align*}
     \mb D^3 = \{(d_1,d_2,d_3) \in \mb Z^{\oplus 3} \ \vert \ d_1 \le d_2 \le d_3\}
 \end{align*}
 and write for $\ul d = (d_1,d_2,d_3) \in \mb D^3$
 \begin{align*}
     \mc O(\ul{d}) = \mc O(d_3) \oplus \mc O(d_2) \oplus \mc O(d_1).
 \end{align*}
Note that the direct summands are ordered by decreasing degree. For any ordered subset $\{i_1 < \ldots < i_q\}$ of $\{1,2,3\}$ we denote by $\iota_{i_1\ldots i_q}$ the inclusion into $\mc O(\ul d)$ of the direct sum of the $i_1$-th factor up to the $i_q$-th factor.

For each $\ul d \in \GlobDat^3$ and $i \in \{1,2,3,4\}$ we fix a basis $\f{e}_{\ul d}^{(i)}$ of $\mc O(\ul d) \otimes k(x_i)$ as follows: First let $e_0^{(i)}$ be the image of the constant global section $1$ of $\mc O$ inside $\mc O \otimes k(x_i)$. Consider the global sections $f_\infty := X$ and $g_\infty := Y$ of $\Gamma(\mb P^1,\mc O(1))$, viewed as elements in $\Hom(\mc O,\mc O(1))$. Set $e_1^{(i)} := Y(e_0^{(i)}) \in \mc O(1) \otimes k(x_i)$ for $i = 1,2,3$ and $e_1^{(4)} := X(e_0^{(i)})$. Then for $d \in \mb Z \setminus \{0\}$ set $e_d^{(i)} := (e_1^{(i)})^{\otimes d}$. Finally, for any $\ul d = (d_1,d_2,d_3) \in \mb D^3$ let $\f{e}_{\ul d}^{(i)}$ be the basis 
\begin{align*}
    \f{e}_{\ul d}^{(i)} := \l (\iota_1 \otimes k(x_i))(e_{z_3}^{(i)}), (\iota_2 \otimes k(x_i))(e_{z_2}^{(i)}), (\iota_3 \otimes k(x_i))(e_{z_1}^{(i)})\r
\end{align*}
of $\mc O(\ul d) = \mc O(z_3) \oplus \mc O(z_2) \oplus \mc O(z_1)$.

The basis $\f{e}_{\ul d}^{(i)}$ yields an identification
\begin{align} \label{formula identification flag variety}
    \prod_{i=1}^4 \mc O(\ul d) \otimes k(x_i) \simeq (k^{\oplus 3})^{\times 4}.
\end{align} 
The automorphism group $\ul{\Aut}(\mc O(\ul d))$ acts naturally on $\prod_{i=1}^4 \mc O(\ul d) \otimes k(x_i)$, hence via the above identification on $(k^{\oplus 3})^{\times 4}$ and on $\Pi$. Summarizing, for each $\ul d \in \GlobDat^3$ we obtain an action of $\ul{\Aut}(\mc O(\ul d))$ on $\Pi$, which we denote by $\act_{\ul d}$.

For each $\ul d \in \GlobDat^3$ the morphism there is a natural isomorphism
\begin{align} \label{iso between different isomorphism groups}
    \ul{\Aut}(\mc O(\ul d)) \simeq \ul{\Aut}(\mc O(\ul d) \otimes \mc O(1)) = \ul{\Aut}(\mc O(\ul d^{+1})),
\end{align}
where $\ul d^{+1} = (d_1+1,d_2+1,d_3+1)$ if $\ul d = (d_1,d_2,d_3)$. Now from the concrete construction of the bases $\f{e}_{\ul d}^{(i)}$ it follows that for each $\ul d \in \GlobDat^3$ the isomorphism (\ref{iso between different isomorphism groups}) is compatible with the actions $\act_{\ul d}$ and $\act_{\ul d^{+1}}$. In particular, if $\Omega \subset \Pi(k)$ is a complete system of representatives of the orbits of the action $\act_{\ul d}$, then $\Omega$ is also a complete system of representatives of the orbits of the action $\act_{\ul d^{+1}}$.

If $P \in \Pi(k)$ and $\ul d \in\GlobDat^3$ we denote by $(\mc O(\ul d), P)$ the parabolic bundle with underlying bundle $\mc O(\ul d)$ and parabolic structure given by $P$ using the identification (\ref{formula identification flag variety}).

For $i \in \{1,2,3\}$, $P_1 \in \Pi_{\le i}(k)$ and $P_2 \in \Pi_{>i}(k)$ we denote the point $(P_1,P_2) \in \Pi_{\le i}(k) \times \Pi_{>i}(k) = \Pi(k)$ by $P_1 \lrcorner P_2$. Similarly, for subsets $\mc P_1 \subset \Pi_{\le i}(k)$ and $\mc P_2 \subset \Pi_{>i}(k)$ we set
\begin{align*}
    \mc P_1 \lrcorner \mc P_2 := \{P_1 \lrcorner P_2 \ \vert \ P_i  \in \mc P_i\}.
\end{align*}

Note that the symmetric group on three letters $\mf S_3$ acts on $\Pi_{\le 3}$ by permuting the factors $\Pi_i$. For a subset $\Omega \subset \Pi_{\le 3}(k)$ we set
\begin{align*}
    \Omega^{\cup \langle \mf S_3 \rangle} := \left\{\sigma.P \ \vert \ \sigma \in \mf S_3, \ P \in \Omega \right\}.
\end{align*}

To formulate results in a concise way we introduce the following notion. Recall that $D^\tw$ denotes the concrete datum defined in (\ref{Definition of Data in Example}).

\begin{Def}
    Let $\mc D \subset \GlobDat$ and $i \in \{1,2,3,4\}$.
    \begin{itemize}
        \item[1)] A subset $\mc P \subset \Pi_{\le i}(k)$ is said to be $\mc D$\textit{-exhaustive} iff  for all $\ul d \in \mc D$ we have
        \begin{align*}
            \bigcup_{P \in \mc P} \ul{\Aut}(\mc O(\ul d))(k).P = \Pi_{\le i}(k).
        \end{align*}
        \item[2)] Assume $i < 4$. Let $\mc P_0 \subset \Pi_{\le i}(k)$, $\mc P_1 \subset \Pi_{> i}(k)$. Then $\mc P_1$ is said to be a $\mc D$\textit{-exhaustion of} $\mc P_0$ iff for all $\ul d \in \mc D$ and all $P \in \mc P_0$ we have
        \begin{align*}
            \bigcup_{P\p \in \mc P_1} \Stab_{\ul{\Aut}(\mc O(\ul d))}(P)(k).P\p = \Pi_{> i}(k). 
        \end{align*}
        If $\mc P_0 = \{P\}$, we say $\mc P_1$ is a $\mc D$-exhaustion of $P$.
    \item[3)] A subset $\mc P \subset \Pi_{\le i}(k)$ is said to be $\mc D$-irrelevant resp. $\mc D$-noncuspidal resp. $\mc D$-decomposable iff there exists a $\mc D$-exhaustion $P_1$ of $\mc P$ such that for every point $P \in \mc P \lrcorner \mc P_1$ and every $\ul d \in \mc D$, the parabolic bundle $(\mc O(\ul d),P)$ is $D^\tw$-irrelevant resp. $D^\tw$-noncuspidal resp. decomposable.
    \end{itemize}
\end{Def}

The purpose of this definition is to efficiently exclude many $D^\tw$-irrelevant resp. $D^\tw$-noncuspidal points at once.

Since we will have to write down points in coordinates, we explain our notation of points in projective space and flag varieties: We will denote points of $\mb P^2(k)$ by column vectors
\begin{align*}
    \begin{bmatrix}
        a_1 \\
        a_2 \\
        a_3
    \end{bmatrix}
\end{align*}
and elements of $\mb \FL_3(k)$ by 
\begin{align*}
    \lb \begin{pmatrix}
        a_1 \\
        a_2 \\
        a_3
    \end{pmatrix}, \begin{pmatrix}
        b_1 \\
        b_2 \\
        b_3
    \end{pmatrix} \rb, 
\end{align*}
where the two vectors 
\begin{align*}
   \f{a} := \begin{pmatrix}
        a_1 \\
        a_2 \\
        a_3
    \end{pmatrix}, \f{b} := \begin{pmatrix}
        b_1 \\
        b_2 \\
        b_3
    \end{pmatrix} \in k^{\oplus 3}
\end{align*}
are non-zero and linearly independent and for linearly independent $\f{a}\p, \f{b}\p \in k^{\oplus 3}$ we have
\begin{align*}
    \lb \f{a}, \f{b} \rb = \lb \f{a}\p, \f{b}\p \rb
\end{align*}
iff $\f{a}$ is a multiple of $\f{a}\p$ and $\f{a}, \f{b}$ span the same vector space as $ \f{a}\p, \f{b}\p$.

Finally, let us explain how we will perform the calculation: First, we note that we may write $\Bun_D(k)^\inter$ as a disjoint union (of groupoids)
\begin{align*}
    \Bun_D(k)^\inter = \sqcup_{d \in \mb Z} \Bun_D^d(k)^\inter.
\end{align*}
By Theorem \ref{Proposition interesting locus and hecke operators} the groupoid $\Bun_D^d(k)^\inter$ is independent of $d$ and an explicit equivalence 
\begin{align*}
    \Bun_D^{d+1}(k)^\inter \simeq \Bun_D^d(k)^\inter
\end{align*}
is given by the Hecke operator $\mb T_4^{(l)}$. Also, tensorization by a line bundle of degree $1$ induces an isomorphism 
\begin{align*}
    \Bun_D^{d}(k)^\inter \simeq \Bun_D^{d+3}(k)^\inter.
\end{align*}
These operations are related as follows: The threefold selfcomposition of $\mb T_4^{(l)}$ is equivalent to tensorization with $\mc O(-x_4)$. 

Now we take a closer look at coordinates: We can decompose each groupoid $\Bun_D^d(k)^\inter$ further by fixing the underlying global datum,
\begin{align*}
    \Bun_D^d(k)^\inter = \bigsqcup_{\ul d \in \GlobDat_d^3} \Bun_D^{\ul d}(k)^\inter,
\end{align*}
where $\GlobDat_d^3$ denotes the set of global data of rank $3$ and degree $d$. Each stack $\Bun_D^{\ul d}$ is, as a reduced locally closed substack of $\Bun_D$, isomorphic to a quotient stack,
\begin{align*}
    \Bun_D^{\ul d} \simeq \lb \Pi /\ul{\Aut}(\mc O(\ul d)) \rb.
\end{align*}
Now the discussion after (\ref{iso between different isomorphism groups}) together with Theorem \ref{Proposition interesting locus and hecke operators} imply the following: Given a subset $\Omega \subset \Pi(k)$ such that the image of $\Omega$ inside $\lb \Pi /\ul{\Aut}(\mc O(\ul d)) \rb(k)$ is a complete system of representatives of isomorphism classes in $\Bun_D^{\ul d}(k)^\inter$, then the image of $\Omega$ inside $\lb \Pi /\ul{\Aut}(\mc O(\ul d^{+1})) \rb(k)$ is a complete system of representatives of isomorphism classes in $\Bun_D^{\ul d^{+1}}(k)^\inter$ (recall $\ul d^{+m} = (d_1+m,d_2+m,d_3+m)$ if $\ul d = (d_1,d_2,d_3)$).

If $\ul d, \ul d\p \in \GlobDat^3$ are not related by $\ul d\p = \ul d^{+m}$ for some $m \in \mb Z$, then there is no such relation between coordinates of points in $\Bun_D^{\ul d}(k)^\inter$ and $\Bun_D^{\ul d\p}(k)^\inter$. We will calculate a list of representatives of isomorphism classes in $\Bun_D^d(k)^\inter$ for one value of $d = 0$ and $d=1$. The reason is that there is a nice chart of $\Bun_D^{0,\min}$, which we use for many calculations, but this chart can't immediately be extended to a neighborhood of the slightly larger locus $\vert \Bun_D^0 \vert^\inter$. However, when passing to the degree $1$ component, this is possible, since each point in $\vert \Bun_D^1 \vert^\inter$ has the same underlying bundle $\mc O(0,0,1)$, which is not true for $\vert \Bun_D^0 \vert^\inter$ (see Proposition \ref{First Proposition about Geometry}).

To obtain the list for the degree $1$ component we apply the elementary Hecke operator $\mb T_4^{(u)}$ to the points in the list for the degree $0$ component and use Theorem \ref{Proposition interesting locus and hecke operators}.

In the calculations, we will distinguish four cases. For $D \in \GlobDat^3$ of degree $0$, its sequence of multiplicities (cf. subsection \ref{subsection Global Data}) has to be one of the ordered sequences:
\begin{align*}
    (3), (1,2), (2,1), (1,1,1).
\end{align*}
Each of the subsections of this chapter deals with one of the four cases. We emphasize again that in each case it is our aim to exhibit (efficiently) all points in $\vert \Bun_D^0(k) \vert$ which have the given sequence of multiplicities and are $D^\tw$-irrelevant or $D^\tw$-noncuspidal and then write down a list of the remaining ones, which constitute the interesting locus.

\subsection{The case $(3)$} \label{Subsection cas (3)}
Let $\mc D_{(3)} := \{\ul d \in\GlobDat \ \vert \ \deg(\ul d) = 0, \ \secmul(\ul d) = (3)\}$. This set only contains the element $\ul 0 = (0,0,0)$ with $\mc O(\ul 0) = \mc O^{\oplus 3}$. Its automorphism group is $\ul{\Aut}(\mc O(\ul 0)) = \GL_3$ and the action on $\Pi$ is the standard diagonal action.

Consider the following four points in $\Pi_{\le 3}(k)$: 
\begin{align*}
    P_{(3)}^{(4)} := \l  \begin{bmatrix}
                       1 \\ 0 \\ 0 \\
                      \end{bmatrix},  \begin{bmatrix}
                         1 \\ 0 \\ 0 \\
                        \end{bmatrix}, \begin{bmatrix}
                         1 \\ 0 \\ 0 \\
                          \end{bmatrix} \r \\
    P_{(3)}^{(3)} := \l  \begin{bmatrix}
                       1 \\ 0 \\ 0 \\
                      \end{bmatrix},  \begin{bmatrix}
                         1 \\ 0 \\ 0 \\
                        \end{bmatrix}, \begin{bmatrix}
                         0 \\ 1 \\ 0 \\
                          \end{bmatrix} \r \\
    P_{(3)}^{(2)} := \l  \begin{bmatrix}
                       1 \\ 0 \\ 0 \\
                      \end{bmatrix},  \begin{bmatrix}
                         0 \\ 1 \\ 0 \\
                        \end{bmatrix}, \begin{bmatrix}
                         1 \\ 1 \\ 0 \\
                          \end{bmatrix} \r \\
    P_{(3)}^{(1)} := \l  \begin{bmatrix}
                       1 \\ 0 \\ 0\\
                      \end{bmatrix},  \begin{bmatrix}
                         0 \\ 1 \\ 0 \\
                        \end{bmatrix}, \begin{bmatrix}
                         0 \\ 0 \\ 1 \\
                          \end{bmatrix} \r \\
\end{align*}
Set $\mc P_{(3)} := \{P_{(3)}^{(1)}, P_{(3)}^{(2)}, P_{(3)}^{(4)}\} \cup \{ P_{(3)}^{(3)} \}^{\cup \langle \mf S_3 \rangle}$ (which has $3+3 = 6$ elements).

\begin{Lem}
    $\mc P_{(3)}$ is $\mc D_{(3)}$-exhaustive.
\end{Lem}

\begin{Lem}
    The set $\{P_{(3)}^{(2)},\ P_{(3)}^{(4)} \} \cup \{ P_{(3)}^{(3)} \}^{\cup \langle \mf S_3 \rangle} = \mc P_{(3)}  \setminus \{P_{(3)}^{(1)}\}$ is $\mc D_{(3)}$-irrelevant.
\end{Lem}
\begin{proof}
        Let $\mc Q_{(3)}^{(4)}$ be the set consisting of the following six points in $\FL_3(k)$:
        \[
\setlength{\arraycolsep}{8pt}
\begin{array}{ccc}
\left(
\begin{bmatrix}0\\0\\1\end{bmatrix},
\begin{bmatrix}0\\1\\0\end{bmatrix}
\right)
&
\left(
\begin{bmatrix}0\\0\\1\end{bmatrix},
\begin{bmatrix}1\\0\\0\end{bmatrix}
\right)
&
\left(
\begin{bmatrix}0\\1\\0\end{bmatrix},
\begin{bmatrix}0\\0\\1\end{bmatrix}
\right)
\\[4pt]
\left(
\begin{bmatrix}0\\1\\0\end{bmatrix},
\begin{bmatrix}1\\0\\0\end{bmatrix}
\right)
&
\left(
\begin{bmatrix}1\\0\\0\end{bmatrix},
\begin{bmatrix}0\\0\\1\end{bmatrix}
\right)
&
\left(
\begin{bmatrix}1\\0\\0\end{bmatrix},
\begin{bmatrix}0\\1\\0\end{bmatrix}
\right)
\end{array}
\]
Then $\mc Q_{(3)}^{(4)}$ is a $\mc D_{(3)}$-exhaustion of $P_{(3)}^{(4)}$ and for any $P \in \{P_{(3)}^{(4)}\} \lrcorner \mc Q_{(3)}^{(4)}$ the inclusion $\iota_1$ is $\mc D_{(3)}$-$D^\tw$-irrelevance recognizing for the parabolic bundle $\mc E\b := (\mc O(\ul d), P)$: Indeed, direct inspection shows that that $\Inh_{\mc E\b}(\iota_1) \hookrightarrow \mc E\b$ admits a complement and that $D(\Inh_{\mc E\b}(\iota))$ is of the form $((1,0),(1,0),(1,0),e)$ for $e \in \{(1,0,0),(0,1,0),(0,0,1)\}$. Hence, we obtain
\begin{align*}
    \mf R\l \Res_{D^\tw}\l D(\Inh_{\mc E\b} \r(\iota_1)) \r = \lambda^{3} \cdot \lambda^{-1} = \lambda^2,
\end{align*}
which is $\neq 1$ in view of $\lambda \notin \{\pm 1\}$. This handles the case of $P_{(3)}^{(4)}$. \\
Let $\mc Q_{(3)}^{(3)}$ be the set consisting of the following nine points in $\FL_3(k)$:
\[
\setlength{\arraycolsep}{8pt}
\begin{array}{ccc}
\left(
\begin{bmatrix}0\\0\\1\end{bmatrix},
\begin{bmatrix}1\\1\\0\end{bmatrix}
\right)
&
\left(
\begin{bmatrix}0\\0\\1\end{bmatrix},
\begin{bmatrix}0\\1\\0\end{bmatrix}
\right)
&
\left(
\begin{bmatrix}0\\0\\1\end{bmatrix},
\begin{bmatrix}1\\0\\0\end{bmatrix}
\right)
\\[4pt]
\left(
\begin{bmatrix}1\\1\\0\end{bmatrix},
\begin{bmatrix}0\\0\\1\end{bmatrix}
\right)
&
\left(
\begin{bmatrix}1\\1\\0\end{bmatrix},
\begin{bmatrix}1\\0\\0\end{bmatrix}
\right)
&
\left(
\begin{bmatrix}0\\1\\0\end{bmatrix},
\begin{bmatrix}0\\0\\1\end{bmatrix}
\right)
\\[4pt]
\left(
\begin{bmatrix}0\\1\\0\end{bmatrix},
\begin{bmatrix}1\\0\\0\end{bmatrix}
\right)
&
\left(
\begin{bmatrix}1\\0\\0\end{bmatrix},
\begin{bmatrix}0\\0\\1\end{bmatrix}
\right)
&
\left(
\begin{bmatrix}1\\0\\0\end{bmatrix},
\begin{bmatrix}0\\1\\0\end{bmatrix}
\right)
\end{array}
\]
Then $\mc Q_{(3)}^{(3)}$ is a $\mc D_{(3)}$-exhaustion of $P_{(3)}^{(3)}$ and for any $P \in \{P_{(3)}^{(3)}\} \lrcorner \mc Q_{(3)}^{(3)}$ the inclusion $\iota_3$ is $\mc D_{(3)}$-$D^\tw$-irrelevance recognizing. To see this one argues as before, but uses 
\begin{align*}
    \lambda^{-1} \neq 1
\end{align*}
in the end. Also note that now the inclusions $\iota_1$ or $\iota_2$ do not always work, namely they do not work for the first, fourth, fifth and sixth point in the list defining $\mc Q_{(3)}^{(3)}$.\\
The argument for the points in $\{ P_{(3)}^{(3)} \}^{\cup \langle \mf S_3 \rangle} \setminus \{P_{(3)}^{(3)}\}$ as well as for the point $P_{(3)}^{(2)}$ is analogous to the one for $P_{(3)}^{(3)}$. 
\end{proof}

Define for $i \in \{1,2,3\}$ points in $\mb A^3(k)$ as follows ($\delta_{ij}$ denotes the Kronecker delta):
\begin{align*}
    M^{(1)}_i := \begin{pmatrix}
        \delta_{i1} \\
        \delta_{i2} \\
        \delta_{i3}
    \end{pmatrix},\   M^{(2)}_i := \begin{pmatrix}
        1-\delta_{i1} \\
        1-\delta_{i2} \\
        1-\delta_{i3} \\
        \end{pmatrix}, \  M^{(3)} := \begin{pmatrix}
        1 \\
        1 \\
        1
    \end{pmatrix}.
\end{align*}
We fix an element $\sigma \in \mf S_3$ with no fixed points. Also consider the inclusion
\begin{align*}
    \iota_{\mb P^1,\mb P^2} : \mb P^1 &\hookrightarrow \mb P^2 \\
    v = [v_1:v_2] &\mapsto P_v := \begin{pmatrix}
        v_1 \\
        v_2 \\
        0
    \end{pmatrix}.
\end{align*}
Define subsets of $\FL_3(k)$ as follows:
\begin{align*}
    \mc M^{(1)} &:= \left\{ \lb M^{(1)}_i, M^{(1)}_j \rb \ \vert \ i \in \{1,2,3\}, \ j \in \{1,2,3\} \setminus \{i\} \right\} \cup \left\{ \lb M^{(1)}_i, M^{(2)}_i \rb \ \vert \ i \in \{1,2,3\} \right\} \\
    \mc M^{(2)}_A &:= \left\{ \lb M^{(2)}_i, M^{(1)}_{\sigma(i)} \rb \ \vert \ i \in \{1,2,3\}\right\} \\
    \mc M^{(2)}_B &:= \left\{ \lb M^{(2)}_i, M^{(1)}_i \rb \ \vert \ i \in \{1,2,3\} \right\} \\
    \mc M^{(2)}_C &:= \left\{ \lb M^{(2)}_i, M^{(2)}_{\sigma(i)} \rb \ \vert \ i \in \{1,2,3\} \right\}. \\
    \mc M^{(3)} &:= \left\{ \lb M^{(3)}, P_v \rb \ \vert \ v \in \mb P^1(k) \right\} 
\end{align*}

\begin{Lem}
The following hold:
\begin{itemize}
    \item[1)] The union of the sets $\mc M^{(1)}$, $\mc M^{(2)}_A$, $\mc M^{(2)}_B$, $\mc M^{(2)}_C$ and $\mc M^{(3)}$ is a $\mc D_{(3)}$-exhaustion for $P_{(3)}^{(1)}$. 
    \item[2)] $\{P_{(3)}^{(1)}\} \lrcorner \mc M^{(1)}$ and $\{P_{(3)}^{(1)}\} \lrcorner \mc M^{(2)}_A$ are $\mc D_{(3)}$-noncuspidal.
    \item[3)] The three points in $\{P_{(3)}^{(1)}\} \lrcorner \mc M^{(2)}_B$ are decomposable, but neither $\mc D_{(3)}$-noncuspidal nor $\mc D_{(3)}$-irrelevant.
    \item[4)] All points in $\{P_{(3)}^{(1)}\} \lrcorner \l\mc M_C^{(2)} \cup \mc M^{(3)}\r$ are $\mc D_{(3)}$-indecomposable.
\end{itemize}
\end{Lem}
\begin{proof}
    1) is easily verified by direct calculation. \\
    Consider the case of $\{P_{(3)}^{(1)}\} \lrcorner \mc M^{(2)}_A$. This set consists of three points of the form $\lb M^{(2)}_i, M^{(1)}_{\sigma(i)} \rb$ for some $i \in \{1,2,3\}$. Write $\{1,2,3\} = \{i,l,m\}$ with $l < m$. We claim that for such a point $\iota_{lm}$ is $\mc D_{(3)}$-$D^\tw$-noncuspidality recognizing. We prove this in the case $i=1$, and $\sigma(i) = 2$, the proof in the other cases (i.e. different $i \in \{1,2,3\}$ and other choices of $\sigma$) is the same up to permuting the coordinates. We have
    \begin{align*}
        P_{(3)}^{(1)} \lrcorner \lb M^{(2)}_1, M^{(1)}_2 \rb = \l \begin{bmatrix}
                       1 \\ 0 \\ 0\\
                      \end{bmatrix},  \begin{bmatrix}
                         0 \\ 1 \\ 0 \\
                        \end{bmatrix}, \begin{bmatrix}
                         0 \\ 0 \\ 1 \\
                          \end{bmatrix}, \lb \begin{pmatrix}
                              0 \\
                              1 \\
                              1
                          \end{pmatrix}, \begin{pmatrix}
                              0 \\
                              1 \\
                              0
                          \end{pmatrix} \rb \r.
    \end{align*}
Set
\begin{align*}
    \mc E\b := \l \mc O^{\oplus 3} , P_{(3)}^{(1)} \lrcorner \lb M^{(2)}_1, M^{(1)}_2 \rb \r.
\end{align*}
We obtain
\begin{align*}
    \Inh_{\mc E\b}(\iota_{23}) = \l \mc O^{\oplus 2}, \l \begin{bmatrix}
                        0 \\ 0\\
                      \end{bmatrix},  \begin{bmatrix}
                          1 \\ 0 \\
                        \end{bmatrix}, \begin{bmatrix}
                          0 \\ 1 \\
                          \end{bmatrix}, \lb \begin{pmatrix}
                              1 \\
                              1
                          \end{pmatrix}, \begin{pmatrix}
                              1 \\
                              0
                          \end{pmatrix} \rb \r \r.
\end{align*}
and the natural morphism $\Inh_{\mc E\b}(\iota_{23})  \hookrightarrow \mc E\b$ admits a complement, which we denote by $\mc Q\b$. Then we have
\begin{align*}
    \mc Q\b = \l \mc O, \l 1, 0, 0, \lb 0,0 \rb \r \r.
\end{align*}
Since $\Ext^1(\mc O, \mc O^{\oplus 2}) = 0$, any extension of $\mc Q\b$ by $\Inh_{\mc E\b}(\iota_{23})$ is equivalent to an extension of the form
\begin{align} \label{general extension}
    0 \to \Inh_{\mc E\b}(\iota_{23}) \xrightarrow{\iota_{23}} \l \mc O^{\oplus 3}, \l \begin{bmatrix}
                       1 \\ a_2 \\ a_3\\
                      \end{bmatrix},  \begin{bmatrix}
                         0 \\ 1 \\ 0 \\
                        \end{bmatrix}, \begin{bmatrix}
                         0 \\ 0 \\ 1 \\
                          \end{bmatrix}, \lb \begin{pmatrix}
                              0 \\
                              1 \\
                              1
                          \end{pmatrix}, \begin{pmatrix}
                              0 \\
                              1 \\
                              0
                          \end{pmatrix} \rb \r \r \xrightarrow{\pr_1} \mc Q\b \to 0,
\end{align}
where the first maps is induced by the inclusion of the second and third coordinate, the second map is induced by the projection to the first coordinate and $a_2,a_3 \in k$. Now the automorphism
\begin{align*}
    \begin{pmatrix}
        1 & & \\
        -a_2 & 1 & \\
        -a_3 & & 1
    \end{pmatrix} \in \GL_3(k) = \ul{\Aut}(\mc O^{\oplus 3})
\end{align*}
exhibits the extension \ref{general extension} as equivalent to the trivial extension of $\mc Q\b$ by $\Inh_{\mc E\b}(\iota_{23})$. \\
The argument for points in $\{P_{(3)}^{(1)}\} \lrcorner \mc M^{(1)}$ is similar. This finishes the proof of 2). \\

Concerning 3), the set $\{P_{(3)}^{(1)}\} \lrcorner \mc M^{(2)}_B$ contains three points, namely $P_{(3)}^{(1)} \lrcorner \lb M^{(2)}_i, M^{(1)}_i \rb$ for $i \in \{1,2,3\}$. They are decomposable. We argue by contradiction that they are neither $\mc D_{(3)}$-irrelevant nor $\mc D_{(3)}$-noncuspidal. In doing so we restrict to the case $i = 1$, the other cases are handled similarly. We have
    \begin{align*}
        P_{(3)}^{(1)} \lrcorner \lb M^{(2)}_i, M^{(1)}_i \rb = \l \begin{bmatrix}
                       1 \\ 0 \\ 0\\
                      \end{bmatrix},  \begin{bmatrix}
                         0 \\ 1 \\ 0 \\
                        \end{bmatrix}, \begin{bmatrix}
                         0 \\ 0 \\ 1 \\
                          \end{bmatrix}, \lb \begin{pmatrix}
                              0 \\
                              1 \\
                              1
                          \end{pmatrix}, \begin{pmatrix}
                              1 \\
                              0 \\
                              0
                          \end{pmatrix} \rb \r.
    \end{align*}
Set
\begin{align*}
    \mc E\b := \l \mc O^{\oplus 3} , P_{(3)}^{(1)} \lrcorner \lb M^{(2)}_i, M^{(1)}_i \rb \r.
\end{align*}
We first prove the following claim: If $\iota : \mc F\b \hookrightarrow \mc E$ is a monomorphism of parabolic vector bundles that admits a complement with $\rk(\mc F\b) \in \{1,2\}$, then $\iota$ is isomorphic either to $\iota_1$ or to $\iota_{23}$. 

We first assume $\rk(\mc F\b) = 1$. Since $\iota$ admits a complement, we may assume $\UL(\mc F\b) = \mc O$. Then we can write $\iota$ as 
\begin{align*}
    \iota = (a_1,a_2,a_3) : \mc O \hookrightarrow \mc O^{\oplus 3}
\end{align*}
(the map on parabolic vector bundles is uniquely determined by the map on underlying bundles). Write
\begin{align*}
    D(\mc F\b) = (q_1,q_2,q_3,q_4)
\end{align*}
with $q_i \in \{(1,0),(0,1)\}$ for $i = 1,2,3$ and $q_4 \in \{(1,0,0),(0,1,0),(0,0,1)\}$. Since the first three entries of $ P_{(3)}^{(1)} \lrcorner \lb M^{(2)}_i, M^{(1)}_i \rb$ are pairwise different, the quantity $\#\{i \in \{1,2,3\} \ \vert \ q_i = (0,1)\}$ is either $2$ or $3$. On the other hand, since $\rk(\mc F\b) = 1$, its complement can at most contain two of the three lines
\begin{align*}
    \begin{bmatrix}
                       1 \\ 0 \\ 0\\
                      \end{bmatrix},  \begin{bmatrix}
                         0 \\ 1 \\ 0 \\
                        \end{bmatrix}, \begin{bmatrix}
                         0 \\ 0 \\ 1 \\
                          \end{bmatrix}.
\end{align*}
Hence, we must have $\#\{i \in \{2,3,4\} \ \vert \ q_i = (0,1)\} = 2$, i.e. $\iota$ has to be a coordinate inclusion: $\iota = \iota_j$ for some $j \in \{1,2,3\}$. If $\iota = \iota_2$, then its complement must contain all three lines 
\begin{align*}
     \begin{bmatrix}
                       1 \\ 0 \\ 0\\
                      \end{bmatrix},  \begin{bmatrix}
                         0 \\ 0 \\ 1 \\
                        \end{bmatrix}, \begin{bmatrix}
                         0 \\ 1 \\ 1 \\
                          \end{bmatrix}
\end{align*}
and thus have rank $3$, a contradiction. A similar contradiction occurs in the case $\iota_3$. Hence, we conclude $\iota = \iota_1$ and have shown the claim in the case $\rk(\mc F\b) = 1$. The case $\rk(\mc F\b) = 2$ is similar. 

Now assume that $\mc E\b$ is irrelevant. Then there exists a $D^\tw$-irrelevance recognizing monomorphism of parabolic vector bundles $\iota : \mc F\b \hookrightarrow \mc E$ and by symmetry for direct sum decompositions of the notion of irrelevance we may assume $\rk(\mc F\b) = 1$. Using the claim, we may assume $\iota = \iota_1$. But then 
\begin{align*}
    \mf R(\Res_{D^\tw}(D(\mc F\b))) = \lambda \cdot \lambda^{-1} = 1,
\end{align*}
so $\iota$ is not $D^\tw$-irrelevance recognizing. 

Now assume that $\mc E\b$ is noncuspidal and let $\iota : \mc F\b \hookrightarrow \mc E\b$ be $D^\tw$-noncuspidality recognizing. The claim implies that $\iota = \iota_1$ or $\iota = \iota_{23}$. Thus it suffices to show
\begin{align*}
    \Ext^1(\Inh_{\mc E\b}(\iota_{23}),\Inh_{\mc E\b}(\iota_1)) \neq 0, \\
    \Ext^1(\Inh_{\mc E\b}(\iota_1),\Inh_{\mc E\b}(\iota_{23})) \neq 0.
\end{align*}
In the first case consider the extension
\begin{align} 
    0 \to \Inh_{\mc E\b}(\iota_1) \xrightarrow{\iota_1} \l \mc O^{\oplus 3}, \l \begin{bmatrix}
                       1 \\ 0 \\ 0\\
                      \end{bmatrix},  \begin{bmatrix}
                         a \\ 1 \\ 0 \\
                        \end{bmatrix}, \begin{bmatrix}
                         b \\ 0 \\ 1 \\
                          \end{bmatrix}, \lb \begin{pmatrix}
                              c \\
                              1 \\
                              1
                          \end{pmatrix}, \begin{pmatrix}
                              1 \\
                              0 \\
                              0
                          \end{pmatrix} \rb \r \r \xrightarrow{\iota_{23}} \Inh_{\mc E\b}(\iota_{23}) \to 0
\end{align}
for $a,b,c \in k$. This is equivalent to the trivial extension iff there is an element $A$ in the subgroup
\begin{align*}
    \left\{\begin{pmatrix}
        \tilde a_1 & \tilde c_1 & \tilde c_2 \\
        & \tilde a_2 & \\
        & & \tilde a_2
    \end{pmatrix} \ \vert \ \tilde a_i \in k^{\times}, \ \tilde c_i \in k   \right\}
\end{align*}
of $\GL_3$ such that $A$ maps
\begin{align*}
    A : \begin{pmatrix}
                         a \\ 1 \\ 0 \\
                        \end{pmatrix} \mapsto  \begin{pmatrix}
                         0 \\ 1 \\ 0 \\
                        \end{pmatrix}, \ \begin{pmatrix}
                         b \\ 0 \\ 1 \\
                          \end{pmatrix} \mapsto \begin{pmatrix}
                         0 \\ 0 \\ 1 \\
                          \end{pmatrix}, \  \begin{pmatrix}
                              c \\
                              1 \\
                              1
                          \end{pmatrix} \mapsto \begin{pmatrix}
                              0 \\
                              1 \\
                              1
                          \end{pmatrix}.
\end{align*}
Such an element $A$ cannot exist as soon $a,b,c$ are chosen such that $a+b \neq c$. This shows $\Ext^1(\Inh_{\mc E\b}(\iota_{23}),\Inh_{\mc E\b}(\iota_1)) \neq 0$. \\

For the second case, consider the extension
\begin{align} 
    0 \to \Inh_{\mc E\b}(\iota_{23}) \xrightarrow{\iota_{23}} \l \mc O^{\oplus 3}, \l \begin{bmatrix}
                       1 \\ a_2 \\ a_3\\
                      \end{bmatrix},  \begin{bmatrix}
                         0 \\ 1 \\ 0 \\
                        \end{bmatrix}, \begin{bmatrix}
                         0 \\ 0 \\ 1 \\
                          \end{bmatrix}, \lb \begin{pmatrix}
                              0 \\
                              1 \\
                              1
                          \end{pmatrix}, \begin{pmatrix}
                              1 \\
                              b_2 \\
                              b_3
                          \end{pmatrix} \rb \r \r \xrightarrow{\iota_1} \Inh_{\mc E\b}(\iota_1) \to 0.
\end{align}
 This is equivalent to the trivial extension iff there is an element $A$ in the subgroup
\begin{align*}
    \left\{\begin{pmatrix}
        \tilde a_1 &  & \\
        \tilde c_1 & \tilde a_2 & \\
        \tilde c_2 & & \tilde a_2
    \end{pmatrix} \ \vert \ \tilde a_i \in k^{\times}, \ \tilde c_i \in k   \right\}
\end{align*}
of $\GL_3$ such that there exists $\eta \in k^{\times }$ such that $A$ maps
\begin{align*}
    A : \begin{pmatrix}
                         1 \\ a_2 \\ a_3 \\
                        \end{pmatrix} \mapsto  \begin{pmatrix}
                         1 \\ 0 \\ 0 \\
                        \end{pmatrix}, \ \begin{pmatrix}
                              1 \\
                              b_2 \\
                              b_3
                          \end{pmatrix} \mapsto \eta\begin{pmatrix}
                              0 \\
                              1 \\
                              1
                          \end{pmatrix} + \begin{pmatrix}
                              1 \\
                              0 \\
                              0
                          \end{pmatrix} = \begin{pmatrix}
                              1 \\
                              \eta \\
                              \eta
                          \end{pmatrix} 
\end{align*}
Such an element $A$ cannot exist as soon as $a_i,b_i$ are chosen such that $a_3-a_2 \neq b_3 - b_2$. This shows $\Ext^1(\Inh_{\mc E\b}(\iota_1),\Inh_{\mc E\b}(\iota_{23})) \neq 0$. 

Finally, 4) follows easily by direct examination.
\end{proof}

\subsection{The cases $(1,2)$ and $(2,1)$} \label{subsection cases (1,2) and (2,1)}
Let $\mc D_{(1,2)} := \{\ul d \in\GlobDat \ \vert \ \secmul(\ul d) = (1,2), \ \deg(\ul d) = 0\}$. This means that for $\ul d \in \mc D_{(1,2)}$ there is $d \in \mb Z_{\ge 1}$ such that $\ul d = (-2d,d,d) \in \mb D^3 = \GlobDat^3$, hence $\mc O(\ul d) = \mc O(d)^{\oplus 2} \oplus \mc O(-2d)$. \\

\begin{Lem} \label{Deomposability in the (1,2)-case}
    There is a set $\Omega \subset \Pi(k)$ such that $\Omega$ is a $\mc D_{(1,2)}$-exhaustion and such that $\Omega$ is $\mc D$-decomposable w.r.t. $\iota_3$.
\end{Lem}
\begin{proof}
    Let $P \in \Pi(k)$ be an arbitrary point. We need to show that for any $\ul d \in \mc D_{(1,2)}$ there is an automorphism $\phi \in \Aut(\mc O(\ul d))$ such that $\phi(P)$ is decomposable w.r.t. $\iota_3$. Write 
    \begin{align*}
   P = \l \begin{bmatrix}
        a_1^{(1)} \\
        a_2^{(1)} \\
        a_3^{(1)}
    \end{bmatrix}, \begin{bmatrix}
        a_1^{(2)} \\
        a_2^{(2)} \\
        a_3^{(2)}
    \end{bmatrix},  \begin{bmatrix}
        a_1^{(3)} \\
        a_2^{(3)} \\
        a_3^{(3)}
    \end{bmatrix}, \lb \begin{pmatrix}
        a_1^{(4)} \\
        a_2^{(4)} \\
        a_3^{(4)}
    \end{pmatrix}, \begin{pmatrix}
        b_1 \\
        b_2 \\
        b_3
    \end{pmatrix} \rb  \r
\end{align*}
and define
    \begin{align*}
        I := \{i \in \{1,2,3,4\} \ \vert \ a_3^{(i)} \neq 0 \}.
    \end{align*}
Since $d \in \mb Z_{\ge 1}$, we have $d-(-2d) \ge 3$, hence $\dim_k(\Hom(\mc O(-2d),\mc O(d)) \ge 4$ and the map
\begin{align*}
    \Hom(\mc O(-2d),\mc O(d)) \to \prod_{i=1}^4 \Hom(\mc O(-2d) \otimes k(x_i), \mc O(d) \otimes k(x_i))
\end{align*}
is surjective. Thus, we find $\psi_1, \psi_2 \in  \Hom(\mc O(-2d),\mc O(d))$ such that for all $i \in I$, $j \in \{1,2\}$ we have
\begin{align*}
    \l\psi_j \otimes k(x_i)\r(1) = \frac{a_j^{(i)}}{a_3^{(i)}}.
\end{align*}
Then if $a_3^{(4)} \neq 0$, for
\begin{align*}
    \phi := \begin{pmatrix}
        1 & & -\psi_1 \\
        & 1 & -\psi_2 \\
        & & 1
    \end{pmatrix}
\end{align*}
the point $\phi(P)$ is decomposable w.r.t. $\iota_3$. If $a_3^{(4)} = 0$, we may consider instead of $I$ the set 
\begin{align*}
    I\p := 
        \{i \in \{1,2,3\} \ \vert \ a_3^{(i)} \neq 0\} \cup \{4\},
\end{align*}
pick by the same argument as before $\psi_1, \psi_2 \in  \Hom(\mc O(-2d),\mc O(d))$ such that for all $i \in I\p \cap \{1,2,3\}$, $j \in \{1,2\}$ we have
\begin{align*}
    \l\psi_j \otimes k(x_i)\r(1) = \frac{a_j^{(i)}}{a_3^{(i)}}.
\end{align*} 
and we have
\begin{align*}
    \l\psi_j \otimes k(x_4)\r(1) = \frac{b_j}{b_3}.
\end{align*} 
Then for $\phi$ as above the point $\phi(P)$ is decomposable w.r.t. $\iota_3$.
\end{proof}

Consider the stratification
\begin{align*}
    \mb P^2 = T_3 \cup T_{12},
\end{align*}
where $T_3$ is the open subscheme of points 
\begin{align*}
    \begin{bmatrix}
        a_1 \\
        a_2 \\
        a_3
    \end{bmatrix}
\end{align*}
with $a_3 \neq 0$ and $T_{12}$ is its complement. For a point
\begin{align*}
    P = (P_1,P_2,P_3,P_4) \in \prod_{i=1}^4 \Pi_i
\end{align*}
we set
\begin{align*}
    w_1(P) := \#\{i \in \{1,2,3\} \ \vert \ P_i \in T_3\} \in \{0,1,2,3\}. 
\end{align*}
Let $\Omega$ be as in \ref{Deomposability in the (1,2)-case}.

\begin{Lem}
Let $P \in \Omega$.
\begin{itemize}
    \item[1)] If $w_1(P) \neq 1$, then $P$ is $\mc D_{(1,2)}$-irrelevant.
    \item[2)] If $w_1(P) = 1$, then $P$ is $\mc D_{(1,2)}$-noncuspidal.
\end{itemize}
\end{Lem}
\begin{proof}
    1) is easily verified. \\
    For 2) we claim that $\iota_{12}$ is $\mc D_{(1,2)}$-$D^\tw$-noncuspidality recognizing. Indeed, first note
    \begin{align*}
        \Ext^1(\mc O(-2d),\mc O(d)^{\oplus 2}) \simeq H^1(\mb P^1,\mc O(3d))^{\oplus 2} = 0
    \end{align*}
    in view of $d \ge 1$. Now the claim follows from the fact that because of $d-(-2d) = 3d \ge 3$, for any two distinct points $y_1,y_2 \in S$ the natural map
    \begin{align*}
        \Hom(\mc O(-2d),\mc O(d)^{\oplus 2}) \to \prod_{i=1}^2 \Hom(\mc O(-2d) \otimes k(y_i), \mc O(d)^{\oplus 2} \otimes k(y_i)) 
    \end{align*}
    is surjective.
\end{proof}

Now we consider the case $(2,1)$. Let $\mc D_{(2,1)} := \{\ul d \in\GlobDat \ \vert \ \secmul(\ul d) = (2,1), \ \deg(\ul d) = 0\}$. This means that for $\ul d \in \mc D_{(2,1)}$ there is $d \in \mb Z_{\ge 1}$ such that $\ul d = (-d,-d,2d)$, hence $\mc O(\ul d) = \mc O(2d) \oplus \mc O(-d)^{\oplus 2}$.  A similar analysis as in the $(1,2)$-case yields:

\begin{Lem}
     There is a set $\Omega \subset \Pi(k)$ such that $\Omega$ is a $\mc D_{(2,1)}$-exhaustion and such that $\Omega$ is $\mc D$-decomposable w.r.t. $\iota_1$. Each point in any such choice of $\Omega$ is $\mc D_{(2,1)}$-irrelevant or $\mc D_{(2,1)}$-noncuspidal.
\end{Lem}

\subsection{The case $(1,1,1)$} \label{subsection case (1,1,1)}
Let $\mc D_{\ul 1} := \{\ul d \in\GlobDat \ \vert \ \secmul(\ul d) = (1,1,1), \ \deg(\ul d) = 0\}$. This means that for $\ul d \in \mc D_{\ul 1}$ there are $d_1, d_2 \in \mb Z$ with $-(d_1+d_2) < d_1 < d_2$ such that $\ul d = (-(d_1+d_2), d_1, d_2)$, hence $\mc O(\ul d) = \mc O(d_2) \oplus \mc O(d_1) \oplus \mc O(-d_1-d_2)$. \\

Consider the following points in $\Pi_{\le 3}(k)$:
\[
\setlength{\arraycolsep}{5pt}
\begin{array}{c@{\;}c@{\qquad}c@{\;}c}
P_{\underline{1}}^{(1)} &=
\begin{bmatrix}0&0&0\\0&0&1\\1&1&1\end{bmatrix}
&
P_{\underline{1}}^{(7)} &=
\begin{bmatrix}0&1&1\\0&0&0\\1&0&0\end{bmatrix}
\\[5pt]
P_{\underline{1}}^{(2)} &=
\begin{bmatrix}0&0&0\\0&0&0\\1&1&1\end{bmatrix}
&
P_{\underline{1}}^{(8)} &=
\begin{bmatrix}0&0&1\\1&1&1\\0&0&0\end{bmatrix}
\\[5pt]
P_{\underline{1}}^{(3)} &=
\begin{bmatrix}0&0&0\\0&0&1\\1&1&0\end{bmatrix}
&
P_{\underline{1}}^{(9)} &=
\begin{bmatrix}0&0&0\\1&1&1\\0&0&0\end{bmatrix}
\\[5pt]
P_{\underline{1}}^{(4)} &=
\begin{bmatrix}0&0&1\\0&0&0\\1&1&0\end{bmatrix}
&
P_{\underline{1}}^{(10)} &=
\begin{bmatrix}0&0&1\\1&1&0\\0&0&0\end{bmatrix}
\\[5pt]
P_{\underline{1}}^{(5)} &=
\begin{bmatrix}0&0&0\\0&1&1\\1&0&0\end{bmatrix}
&
P_{\underline{1}}^{(11)} &=
\begin{bmatrix}0&1&1\\1&0&0\\0&0&0\end{bmatrix}
\\[5pt]
P_{\underline{1}}^{(6)} &=
\begin{bmatrix}0&0&1\\0&1&0\\1&0&0\end{bmatrix}
&
P_{\underline{1}}^{(12)} &=
\begin{bmatrix}1&1&1\\0&0&0\\0&0&0\end{bmatrix}
\end{array}
\]
We define:
\begin{align*}
    \mc P_{\ul 1, A} &:= \{P_{\ul 1}^{(1)}, P_{\ul 1}^{(3)},P_{\ul 1}^{(7)},P_{\ul 1}^{(8)}, P_{\ul 1}^{(9)}, P_{\ul 1}^{(10)}, P_{\ul 1}^{(11)}, P_{\ul 1}^{(12)}\} \\
    \mc P_{\ul 1, B} &:= \{P_{\ul 1}^{(2)}, P_{\ul 1}^{(4)}, P_{\ul 1}^{(5)}, P_{\ul 1}^{(6)}\} \\
    \mc P_{\ul 1} &:= \mc P_{\ul 1, A} \cup \mc P_{\ul 1, B}
\end{align*}

\begin{Lem}
    $\mc P_{\ul 1}^{\cup \langle \mf S_3 \rangle}$ is $\mc D_{\ul 1}$-exhaustive.
\end{Lem}

\begin{Lem}[The cases 1,3,7,8,9,10,11,12]
    $\mc P_{\ul 1, A}^{\cup \langle \mf S_3 \rangle}$ is $\mc D_{\ul 1}$-irrelevant.
\end{Lem}
\begin{proof}
    We provide the argument for the points in $\mc P_{\ul 1, A}$; the general follows from this by permutation. We claim that there exist subsets $\mc Q_{\ul 1, A}, \mc Q_{\ul 1, B} \subset \FL_3(k)$ such that the following hold:
    \begin{itemize}
        \item[i)] $\mc Q_{\ul 1, A}$ is a $\mc D_{\ul 1}$-exhaustion of $P_{\ul 1}^{(i)}$ for all $i \in \{1,3,7\}$
        \item[ii)] $\mc Q_{\ul 1, B}$ is a $\mc D_{\ul 1}$-exhaustion of $P_{\ul 1}^{(i)}$ for all $i \in \{9,10,11,12\}$
        \item[iii)] For each $i \in \{1,3,7\}$ and $P \in \{P_{\ul 1}^{(i)}\} \lrcorner \mc Q_{\ul 1, A}$ the inclusion $\iota_1$ is $\mc D_{\ul 1}$-$D^\tw$-irrelevance recognizing.
        \item[iv)] For each $i \in \{9,10,11,12\}$ and $P \in \{P_{\ul 1}^{(i)}\} \lrcorner \mc Q_{\ul 1, B}$ the inclusion $\iota_3$ is $\mc D_{\ul 1}$-$D^\tw$-irrelevance recognizing.
    \end{itemize}
    The Lemma then follows except in the case of $P_{\ul 1}^{(8)}$, which one easily verifies separately. \\
    To prove the claim, we let $\mc Q_{\ul 1, A} \subset \FL_3(k)$ be the set consisting of the following points:
\[
\setlength{\arraycolsep}{8pt}
\begin{array}{ccc}
\lb
\begin{pmatrix}0\\\mu_1\\1\end{pmatrix},
\begin{pmatrix}0\\1\\0\end{pmatrix}
\rb
&
\lb
\begin{pmatrix}0\\\mu_2\\1\end{pmatrix},
\begin{pmatrix}1\\0\\0\end{pmatrix}
\rb
&
\lb
\begin{pmatrix}0\\1\\0\end{pmatrix},
\begin{pmatrix}0\\0\\1\end{pmatrix}
\rb
\\[4pt]
\lb
\begin{pmatrix}0\\1\\0\end{pmatrix},
\begin{pmatrix}1\\0\\0\end{pmatrix}
\rb
&
\lb
\begin{pmatrix}1\\0\\0\end{pmatrix},
\begin{pmatrix}0\\\mu_3\\\mu_4\end{pmatrix}
\rb
&
\end{array}
\]
for arbitrary $\mu_1,\mu_2 \in k$ and $[\mu_3:\mu_4] \in \mb P^1(k)$. Then iii) is clearly true. To show i), note that for every $\ul d \in \mc D_{\ul 1}$ the group 
\begin{align*}
    H_A := \left\{\begin{pmatrix}
                a & b & c \\
                & d & 0 \\
                & & d
            \end{pmatrix} \ \vert \ a,d \in k^\times, \ b,c \in k \right\} \subset \GL_3
\end{align*}
is contained in image under $\ul{\Aut}(\mc O(\ul d)) \to \ul{\Aut}(\mc O(\ul d) \otimes k(x_4)) \simeq \GL_3$ of $\Stab_{\ul{\Aut}(\mc O(\ul d))}(P_{\ul 1}^{(1)})$ as well as $\Stab_{\ul{\Aut}(\mc O(\ul d))}(P_{\ul 1}^{(3)})$. Now to prove i), it suffices to check that $Q_{\ul 1, A}$ is a system of representatives of the orbits of the action of $H_A$ on $\FL_3(k)$. We omit the easy verification. \\
Regarding ii) and iv), let $\mc Q_{\ul 1, B} \subset \FL_3(k)$ be the set consisting of the following points:
\[
\setlength{\arraycolsep}{8pt}
\begin{array}{cc}
\lb
\begin{pmatrix}0\\0\\1\end{pmatrix},
\begin{pmatrix}\mu_1\\\mu_2\\0\end{pmatrix}
\rb
&
\lb
\begin{pmatrix}0\\0\\1\end{pmatrix},
\begin{pmatrix}0\\1\\0\end{pmatrix}
\rb
\\[5pt]
\lb
\begin{pmatrix}0\\0\\1\end{pmatrix},
\begin{pmatrix}1\\0\\0\end{pmatrix}
\rb
&
\lb
\begin{pmatrix}\mu_3\\\mu_4\\0\end{pmatrix},
\begin{pmatrix}0\\0\\1\end{pmatrix}
\rb
\\[5pt]
\lb
\begin{pmatrix}\mu_5\\\mu_6\\0\end{pmatrix},
\begin{pmatrix}1\\0\\0\end{pmatrix}
\rb
&
\lb
\begin{pmatrix}0\\1\\0\end{pmatrix},
\begin{pmatrix}0\\0\\1\end{pmatrix}
\rb
\\[5pt]
\lb
\begin{pmatrix}0\\1\\0\end{pmatrix},
\begin{pmatrix}1\\0\\0\end{pmatrix}
\rb
&
\lb
\begin{pmatrix}1\\0\\0\end{pmatrix},
\begin{pmatrix}0\\0\\1\end{pmatrix}
\rb
\\[5pt]
\lb
\begin{pmatrix}1\\0\\0\end{pmatrix},
\begin{pmatrix}0\\1\\0\end{pmatrix}
\rb
&
\end{array}
\]
for arbitrary $[\mu_1:\mu_2], [\mu_3:\mu_4], [\mu_5:\mu_6] \in \mb P^1(k)$. Then iv) is clearly true. To show ii), note that the group 
\begin{align*}
    H_B := \left\{\begin{pmatrix}
                a & 0 & c \\
                & d & e \\
                & & f
            \end{pmatrix} \ \vert \ a,d,f \in k^\times, \ c,e \in k\right\} \subset \GL_3
\end{align*}
is contained in the image under $\ul{\Aut}(\mc O(\ul d)) \to \ul{\Aut}(\mc O(\ul d) \otimes k(\infty)) \simeq \GL_3$ of $\Stab_{\ul{\Aut}(\mc O(\ul d))}(P_{\ul 1}^{(i)})$ for all $i \in \{9,10,11,12\}$. Now to prove ii), it suffices to check that $Q_{\ul 1, B}$ is a system of representatives of the orbits of the action of $H_B$ on $\FL_3(k)$. We omit the easy verification. 
\end{proof}

\begin{Lem}[The case $4$] \label{Case 4}
Let $P \in \{P_{\ul 1}^{(4)}\}^{\cup \langle \mf S_3 \rangle}$. There are sets $\mc Q^{(4)}_A, \mc Q^{(4)}_B \subset \FL_3(k)$ such that
\begin{itemize}
    \item[i)] $\mc Q^{(4)}_A \cup \mc Q^{(4)}_B$ is a $\mc D_{\ul 1}$-exhaustion for $P$,
    \item[ii)] $\{P\} \lrcorner \mc Q^{(4)}_A$ is $\mc D_{\ul 1}$-irrelevant,
    \item[iii)] $\{P\} \lrcorner \mc Q^{(4)}_B$ is $\mc D_{\ul 1}$-noncuspidal.
\end{itemize}
\end{Lem}
\begin{proof}
    We only do the case $P = P_{\ul 1}^{(4)}$, the other cases follow by symmetry. We define $\mc Q^{(4)}_A$ to be the set consisting of the following points in $\FL_3(k)$:
\[
\setlength{\arraycolsep}{8pt}
\begin{array}{cc}
\lb
\begin{pmatrix}0\\0\\1\end{pmatrix},
\begin{pmatrix}0\\1\\0\end{pmatrix}
\rb
&
\lb
\begin{pmatrix}0\\0\\1\end{pmatrix},
\begin{pmatrix}1\\0\\0\end{pmatrix}
\rb
\\[4pt]
\lb
\begin{pmatrix}0\\1\\0\end{pmatrix},
\begin{pmatrix}0\\0\\1\end{pmatrix}
\rb
&
\lb
\begin{pmatrix}0\\1\\0\end{pmatrix},
\begin{pmatrix}1\\0\\0\end{pmatrix}
\rb
\\[4pt]
\lb
\begin{pmatrix}1\\0\\0\end{pmatrix},
\begin{pmatrix}0\\0\\1\end{pmatrix}
\rb
&
\lb
\begin{pmatrix}1\\0\\0\end{pmatrix},
\begin{pmatrix}0\\1\\0\end{pmatrix}
\rb
\end{array}
\]
and $\mc Q^{(4)}_B$ to be the set consisting of the following
three points in $\FL_3(k)$:
\[
\setlength{\arraycolsep}{8pt}
\begin{array}{cc}
Q_1^{(4)} :=
\lb
\begin{pmatrix}0\\1\\1\end{pmatrix},
\begin{pmatrix}0\\1\\0\end{pmatrix}
\rb
&
Q_2^{(4)} :=
\lb
\begin{pmatrix}0\\1\\1\end{pmatrix},
\begin{pmatrix}1\\0\\0\end{pmatrix}
\rb
\\[5pt]
Q_3^{(4)} :=
\lb
\begin{pmatrix}1\\0\\0\end{pmatrix},
\begin{pmatrix}0\\1\\1\end{pmatrix}
\rb
&
\end{array}
\]
Then i) and ii) are easily checked. Regarding iii), we claim that $\iota_1$ is $\mc D_{\ul 1}$-$D^\tw$-noncuspidality recognizing. Indeed, write $\supp(\ul d) = \{d_1 < d_2 < d_3\}$ and set $\mc E_j\b := (\mc O(\ul d), P_{\ul 1}^{(4)} \lrcorner Q_j^{(4)})$. Then
\begin{align*}
   \Ext^1(\mc O(d_1) \oplus \mc O(d_2), \mc O(d_3)) \simeq H^1(\mc O(d_3-d_2)) \oplus H^1(\mc O(d_3-d_1)) = 0,
\end{align*}
in view of $d_1 < d_2 < d_3$. Thus, any class in $\Ext^1(\Inh_{\mc E_j\b}(\iota_{23}),\Inh_{\mc E_j\b}(\iota_1))$ can be represented by an extension of the form
\begin{align} \label{possible extensions in case 4}
\begin{cases}
     0 \to \Inh_{\mc E_j\b}(\iota_{1}) \xrightarrow{\iota_{1}} \l \mc O(\ul d), \l \begin{bmatrix}
                       a \\ 0 \\ 1\\
                      \end{bmatrix},  \begin{bmatrix}
                        b \\ 0 \\ 1 \\
                        \end{bmatrix}, \begin{bmatrix}
                         1 \\ 0 \\ 0 \\
                          \end{bmatrix}, \lb \begin{pmatrix}
                              c \\
                              1 \\
                              1
                          \end{pmatrix}, \begin{pmatrix}
                              d \\
                              1 \\
                              0
                          \end{pmatrix} \rb \r \r \xrightarrow{\iota_{23}} \Inh_{\mc E_j\b}(\iota_{23}) \to 0, & j = 1 \\
                                0 \to \Inh_{\mc E_j\b}(\iota_{1}) \xrightarrow{\iota_{1}} \l \mc O(\ul d), \l \begin{bmatrix}
                       a \\ 0 \\ 1\\
                      \end{bmatrix},  \begin{bmatrix}
                        b \\ 0 \\ 1 \\
                        \end{bmatrix}, \begin{bmatrix}
                         1 \\ 0 \\ 0 \\
                          \end{bmatrix}, \lb \begin{pmatrix}
                              c \\
                              1 \\
                              1
                          \end{pmatrix}, \begin{pmatrix}
                              1 \\
                              0 \\
                              0
                          \end{pmatrix} \rb \r \r \xrightarrow{\iota_{23}} \Inh_{\mc E_j\b}(\iota_{23}) \to 0, & j = 2 \\
                               0 \to \Inh_{\mc E_j\b}(\iota_{1}) \xrightarrow{\iota_{1}} \l \mc O(\ul d), \l \begin{bmatrix}
                       a \\ 0 \\ 1\\
                      \end{bmatrix},  \begin{bmatrix}
                        b \\ 0 \\ 1 \\
                        \end{bmatrix}, \begin{bmatrix}
                         1 \\ 0 \\ 0 \\
                          \end{bmatrix}, \lb \begin{pmatrix}
                              1 \\
                              0 \\
                              0
                          \end{pmatrix}, \begin{pmatrix}
                              d \\
                              1 \\
                              1
                          \end{pmatrix} \rb \r \r \xrightarrow{\iota_{23}} \Inh_{\mc E_j\b}(\iota_{23}) \to 0, & j = 3 
\end{cases}
\end{align}
for some $a,b,c,d \in k$ and it suffices to show that each such extension is isomorphic to the trivial one. We only treat the case $j = 1$, the same argument works in the other two cases. In the case $j = 1$ we first pick a morphism $\phi_1 \in \Hom(\mc O(d_2), \mc O(d_3))$ such that for the action of $\ul{\Aut}(\mc O(\ul d))$ on $k^{\oplus 3}$ through $\ul{\Aut}(\mc O(\ul d) \otimes k(\infty))$ we have
\begin{align*}
    \begin{pmatrix}
        1 & \phi_1 & \\
        & 1 & \\
        & & 1
    \end{pmatrix}.\begin{pmatrix}
        d \\
        1 \\
        0
    \end{pmatrix} = \begin{pmatrix}
        0 \\
        1 \\
        0
    \end{pmatrix}.
\end{align*}
Then we use that in view of $d_3-d_1 \ge 2$ the map 
\begin{align*}
    \Hom(\mc O(d_1),\mc O(d_3)) \to \prod_{i = 1,2,4} \Hom(\mc O(d_1) \otimes k(x_i), \mc O(d_3) \otimes k(x_i)) 
\end{align*}
is surjective. Pick a preimage $\phi_2$ of $(-a,-b,-(c-d))$. Then
\begin{align*}
     \begin{pmatrix}
        1 & \phi_1 & \phi_2\\
        & 1 & \\
        & & 1
    \end{pmatrix} \in \Aut(\mc O(\ul d))
\end{align*}
extends to the desired isomorphism between \ref{possible extensions in case 4} and the trivial extension.  
\end{proof}

\begin{Lem}[The case $5$]
Let $P \in \{P_{\ul 1}^{(5)}\}^{\cup \langle \mf S_3 \rangle}$. There are sets $\mc Q^{(5)}_A, \mc Q^{(5)}_B \subset \FL_3(k)$ such that
\begin{itemize}
    \item[i)] $\mc Q^{(5)}_A \cup \mc Q^{(5)}_B$ is a $\mc D_{\ul 1}$-exhaustion for $P$,
    \item[ii)] $\{P\} \lrcorner \mc Q^{(5)}_A$ is $\mc D_{\ul 1}$-irrelevant,
    \item[iii)] $\{P\} \lrcorner \mc Q^{(5)}_B$ is $\mc D_{\ul 1}$-noncuspidal.
\end{itemize}
\end{Lem}
\begin{proof}
    The proof is very similar to the one of \ref{Case 4}. We omit it.
\end{proof}

\begin{Lem}[The case $6$]
The set $\{P_{\ul 1}^{(6)}\}^{\cup \langle \mf S_3 \rangle}$ is $\mc D_{\ul 1}$-noncuspidal.
\end{Lem}
\begin{proof}
    We show this for $P = P_{\ul 1}^{(6)}$, the general case follows by symmetry. Let $\mc Q^{(6)}$ be the set consisting of the following points in $\FL_3(k)$:
\[
\setlength{\arraycolsep}{8pt}
\begin{array}{cc}
Q_1^{(6)} =
\lb
\begin{pmatrix}0\\0\\1\end{pmatrix},
\begin{pmatrix}0\\1\\0\end{pmatrix}
\rb
&
Q_2^{(6)} =
\lb
\begin{pmatrix}0\\0\\1\end{pmatrix},
\begin{pmatrix}1\\0\\0\end{pmatrix}
\rb
\\[5pt]
Q_3^{(6)} =
\lb
\begin{pmatrix}0\\1\\0\end{pmatrix},
\begin{pmatrix}0\\0\\1\end{pmatrix}
\rb
&
Q_4^{(6)} =
\lb
\begin{pmatrix}0\\1\\0\end{pmatrix},
\begin{pmatrix}1\\0\\0\end{pmatrix}
\rb
\\[5pt]
Q_5^{(6)} =
\lb
\begin{pmatrix}1\\0\\0\end{pmatrix},
\begin{pmatrix}0\\0\\1\end{pmatrix}
\rb
&
Q_6^{(6)} =
\lb
\begin{pmatrix}1\\0\\0\end{pmatrix},
\begin{pmatrix}0\\1\\0\end{pmatrix}
\rb
\end{array}
\]
Clearly, $\mc Q^{(6)}$ is a $\mc D_{\ul 1}$-exhaustion of $P$. We claim that $\iota_{1}$ is $\mc D_{\ul 1}$-$D^\tw$-noncuspidality recognizing for each point in $\{P_{\ul 1}^{(6)}\} \lrcorner \mc Q^{(6)}$. Indeed, write $\supp(\ul d) = \{d_1 < d_2 < d_3\}$ and set $\mc E_j\b := (\mc O(\ul d), P_{\ul 1}^{(6)} \lrcorner Q_j^{(6)})$. Then
\begin{align*}
   \Ext^1(\mc O(d_1) \oplus \mc O(d_2), \mc O(d_3)) \simeq H^1(\mc O(d_3-d_2)) \oplus H^1(\mc O(d_3-d_1)) = 0,
\end{align*}
in view of $d_1 < d_2 < d_3$. Thus, any class in $\Ext^1(\Inh_{\mc E_j\b}(\iota_{23}),\Inh_{\mc E_j\b}(\iota_1))$ can be represented by an extension of the form
\begin{align} \label{possible extensions in case 6}
\begin{cases}
     0 \to \Inh_{\mc E_j\b}(\iota_{1}) \xrightarrow{\iota_{1}} \l \mc O(\ul d), \l \begin{bmatrix}
                       a \\ 0 \\ 1\\
                      \end{bmatrix},  \begin{bmatrix}
                        b \\ 1 \\ 0 \\
                        \end{bmatrix}, \begin{bmatrix}
                         1 \\ 0 \\ 0 \\
                          \end{bmatrix}, \lb\begin{pmatrix}
        c \\ 
        0 \\
        1 
    \end{pmatrix}, \begin{pmatrix}
        d \\
        1 \\
        0
    \end{pmatrix} \rb \r \r \xrightarrow{\iota_{23}} \Inh_{\mc E_j\b}(\iota_{23}) \to 0, & j = 1 \\
                                0 \to \Inh_{\mc E_j\b}(\iota_{1}) \xrightarrow{\iota_{1}} \l \mc O(\ul d), \l \begin{bmatrix}
                       a \\ 0 \\ 1\\
                      \end{bmatrix},  \begin{bmatrix}
                        b \\ 1 \\ 0 \\
                        \end{bmatrix}, \begin{bmatrix}
                         1 \\ 0 \\ 0 \\
                          \end{bmatrix}, \lb\begin{pmatrix}
        c \\ 
        0 \\
        1 
    \end{pmatrix}, \begin{pmatrix}
        1 \\
        0 \\
        0
    \end{pmatrix} \rb \r \r \xrightarrow{\iota_{23}} \Inh_{\mc E_j\b}(\iota_{23}) \to 0, & j = 2 \\
                               0 \to \Inh_{\mc E_j\b}(\iota_{1}) \xrightarrow{\iota_{1}} \l \mc O(\ul d), \l \begin{bmatrix}
                       a \\ 0 \\ 1\\
                      \end{bmatrix},  \begin{bmatrix}
                        b \\ 1 \\ 0 \\
                        \end{bmatrix}, \begin{bmatrix}
                         1 \\ 0 \\ 0 \\
                          \end{bmatrix}, \lb\begin{pmatrix}
        c \\ 
        1 \\
        0 
    \end{pmatrix}, \begin{pmatrix}
        d \\
        0 \\
        1
    \end{pmatrix} \rb \r \r \xrightarrow{\iota_{23}} \Inh_{\mc E_j\b}(\iota_{23}) \to 0, & j = 3 \\
   0 \to \Inh_{\mc E_j\b}(\iota_{1}) \xrightarrow{\iota_{1}} \l \mc O(\ul d), \l \begin{bmatrix}
                       a \\ 0 \\ 1\\
                      \end{bmatrix},  \begin{bmatrix}
                        b \\ 1 \\ 0 \\
                        \end{bmatrix}, \begin{bmatrix}
                         1 \\ 0 \\ 0 \\
                          \end{bmatrix}, \lb\begin{pmatrix}
        c \\ 
        1 \\
        0 
    \end{pmatrix}, \begin{pmatrix}
        1 \\
        0 \\
        0
    \end{pmatrix} \rb \r \r \xrightarrow{\iota_{23}} \Inh_{\mc E_j\b}(\iota_{23}) \to 0, & j = 4 \\
                                0 \to \Inh_{\mc E_j\b}(\iota_{1}) \xrightarrow{\iota_{1}} \l \mc O(\ul d), \l \begin{bmatrix}
                       a \\ 0 \\ 1\\
                      \end{bmatrix},  \begin{bmatrix}
                        b \\ 1 \\ 0 \\
                        \end{bmatrix}, \begin{bmatrix}
                         1 \\ 0 \\ 0 \\
                          \end{bmatrix}, \lb\begin{pmatrix}
        1 \\ 
        0 \\
        0 
    \end{pmatrix}, \begin{pmatrix}
        d \\
        0 \\
        1
    \end{pmatrix} \rb \r \r \xrightarrow{\iota_{23}} \Inh_{\mc E_j\b}(\iota_{23}) \to 0, & j = 5 \\
                               0 \to \Inh_{\mc E_j\b}(\iota_{1}) \xrightarrow{\iota_{1}} \l \mc O(\ul d), \l \begin{bmatrix}
                       a \\ 0 \\ 1\\
                      \end{bmatrix},  \begin{bmatrix}
                        b \\ 1 \\ 0 \\
                        \end{bmatrix}, \begin{bmatrix}
                         1 \\ 0 \\ 0 \\
                          \end{bmatrix}, \lb\begin{pmatrix}
        1\\ 
        0 \\
        0 
    \end{pmatrix}, \begin{pmatrix}
        d \\
        1 \\
        0
    \end{pmatrix} \rb \r \r \xrightarrow{\iota_{23}} \Inh_{\mc E_j\b}(\iota_{23}) \to 0, & j = 6
\end{cases}
\end{align}
for some $a,b,c,d \in k$ and it suffices to show that each such extension is isomorphic to the trivial one. This follows in each case from the fact that for two different points $y_1,y_2 \in S$ in view of $d_3 > d_2 > d_1$ the maps
\begin{align*}
    \Hom(\mc O(d_2),\mc O(d_3)) &\to \prod_{i=1}^2\Hom(\mc O(d_2) \otimes k(y_i),\mc O(d_3) \otimes k(y_i)) \\
    \Hom(\mc O(d_1),\mc O(d_3)) &\to \prod_{i=1}^2\Hom(\mc O(d_1) \otimes k(y_i),\mc O(d_3) \otimes k(y_i))
\end{align*}
are both surjective.
\end{proof}

It remains to discuss the case of $P_{\ul 1}^{(2)}$. Let $\mc Q^{(2)}$ be the subset of $\FL_3(k)$ consisting of the following points:
\[
\setlength{\arraycolsep}{8pt}
\begin{array}{cc}
Q_1^{(2)} =
\lb
\begin{pmatrix}0\\1\\1\end{pmatrix},
\begin{pmatrix}1\\1\\0\end{pmatrix}
\rb
&
Q_2^{(2)} =
\lb
\begin{pmatrix}0\\1\\1\end{pmatrix},
\begin{pmatrix}0\\1\\0\end{pmatrix}
\rb
\\[5pt]
Q_3^{(2)} =
\lb
\begin{pmatrix}0\\1\\1\end{pmatrix},
\begin{pmatrix}1\\0\\0\end{pmatrix}
\rb
&
Q_4^{(2)} =
\lb
\begin{pmatrix}1\\0\\1\end{pmatrix},
\begin{pmatrix}0\\1\\0\end{pmatrix}
\rb
\\[5pt]
Q_5^{(2)} =
\lb
\begin{pmatrix}1\\0\\1\end{pmatrix},
\begin{pmatrix}1\\0\\0\end{pmatrix}
\rb
&
Q_6^{(2)} =
\lb
\begin{pmatrix}0\\0\\1\end{pmatrix},
\begin{pmatrix}0\\1\\0\end{pmatrix}
\rb
\\[5pt]
Q_7^{(2)} =
\lb
\begin{pmatrix}0\\0\\1\end{pmatrix},
\begin{pmatrix}1\\0\\0\end{pmatrix}
\rb
&
Q_8^{(2)} =
\lb
\begin{pmatrix}0\\1\\0\end{pmatrix},
\begin{pmatrix}1\\0\\1\end{pmatrix}
\rb
\\[5pt]
Q_9^{(2)} =
\lb
\begin{pmatrix}0\\1\\0\end{pmatrix},
\begin{pmatrix}0\\0\\1\end{pmatrix}
\rb
&
Q_{10}^{(2)} =
\lb
\begin{pmatrix}0\\1\\0\end{pmatrix},
\begin{pmatrix}1\\0\\0\end{pmatrix}
\rb
\\[5pt]
Q_{11}^{(2)} =
\lb
\begin{pmatrix}1\\0\\0\end{pmatrix},
\begin{pmatrix}0\\1\\1\end{pmatrix}
\rb
&
Q_{12}^{(2)} =
\lb
\begin{pmatrix}1\\0\\0\end{pmatrix},
\begin{pmatrix}0\\0\\1\end{pmatrix}
\rb
\\[5pt]
Q_{13}^{(2)} =
\lb
\begin{pmatrix}1\\0\\0\end{pmatrix},
\begin{pmatrix}0\\1\\0\end{pmatrix}
\rb
&
\end{array}
\]
Let 
\begin{align*}
    \mc D_{\ul 1}\p := \{(-1,0,1)\}.
\end{align*}
Then we have:

\begin{Lem}[The case 2]
    The following hold:
    \begin{itemize}
        \item[i)] $\mc Q^{(2)}$ is a $\mc D_{\ul 1}$-exhaustion for $P_{\ul 1}^{(2)}$.
        \item[ii)] $\{P_{\ul 1}^{(2)}\} \lrcorner \l \mc Q^{(2)} \setminus \{Q_1^{(2)}\} \r$ is $\mc D_{\ul 1}$-irrelevant.
        \item[iii)] $P_{\ul 1}^{(2)} \lrcorner Q_1^{(2)}$ is $\l \mc D_{\ul 1} \setminus \mc D_{\ul 1}\p\r$-irrelevant.
        \item[iv)] $P_{\ul 1}^{(2)} \lrcorner Q_1^{(2)}$ is $\mc D_{\ul 1}\p$-indecomposable and for $\ul d \in \mc D_{\ul 1}\p$ we have
        \begin{align*}
            \ul{\Aut}\l \mc O(\ul d), P_{\ul 1}^{(2)} \lrcorner Q_1^{(2)}\r) = \left\{ \begin{pmatrix}
                a & bY & \\
                & a & \\
                & & a
            \end{pmatrix} \ \vert \ a \in k^{\times}, \ b \in k\right\} \simeq \mb G_m \times \mb G_a.
        \end{align*}
    \end{itemize}
\end{Lem}
\begin{proof}
i) and ii) are easily checked. For iii) let $\ul d \in \mc D_{\ul 1}$ with $\supp(\ul d) = \{d_1 < d_2 < d_3\}$. Then $d_3-d_1 < 3$ implies $\ul d \in \mc D_{\ul 1}\p$, i.e. $(d_1,d_2,d_3) = (-1,0,1)$. If this is not the case, we may first pick $\phi_1 \in \Hom(\mc O(d_2), \mc O(d_3))$ such that for the action of $\Aut(\mc O(\ul d))$ on $k^{\oplus 3}$ through $\Aut(\mc O(\ul d) \otimes k(\infty))$ we have
\begin{align*}
    \begin{pmatrix}
        1 & \phi_1 & \\
        & 1 & \\
        & & 1
    \end{pmatrix}.\begin{pmatrix}
        1 \\
        1 \\
        0
    \end{pmatrix} = \begin{pmatrix}
        0 \\
        1 \\
        0
    \end{pmatrix}.
\end{align*}
Then we use that in view of $d_3-d_1 \ge 3$ the map 
\begin{align*}
    \Hom(\mc O(d_1),\mc O(d_3)) \to \prod_{i = 1,2,3,4} \Hom(\mc O(d_1) \otimes k(x_i), \mc O(d_3) \otimes k(x_i)) 
\end{align*}
is surjective. Pick a preimage $\phi_2$ of $(0,0,0,1)$. Then
\begin{align*}
     \begin{pmatrix}
        1 & \phi_1 & \phi_2\\
        & 1 & \\
        & & 1
    \end{pmatrix} \in \Aut(\mc O(\ul d))
\end{align*}
defines an isomorphism between $(\mc O(\ul d), P_{\ul 1}^{(2)} \lrcorner Q_1^{(2)})$ and $(\mc O(\ul d), P_{\ul 1}^{(2)} \lrcorner Q_2^{(2)})$, which is $\{\ul d\}$-irrelevant by ii). \\
Finally in the case $\ul d \in \mc D_{\ul 1}\p$ one easily checks the formula for the automorphism group. Its semisimple part is isomorphic to the central $\mb G_m$. This implies indecomposability of $\l \mc O(\ul d), P_{\ul 1}^{(2)} \lrcorner Q_1^{(2)}\r$ (\ref{Lemma criterion semisimple}).
\end{proof}

\subsection{The degree $0$ locus: Summary} \label{Subsection about degree zero locus summary}
We summarize the results of the last sections: Let $X = \mb P^1$, $S^\oro = (0,1,t,\infty)$, where $t \in \mb P^1(k) \setminus \{0,1,\infty\}$, $D^\tw = (D,\Lambda)$ with 
\begin{align*}
    D &:= ((1,2),(1,2),(1,2),(1,1,1)), \\
    \Lambda &:= ((\lambda,1),(\lambda,1),(\lambda, 1), (\lambda^{-1},\lambda^{-1},\lambda^{-1})).
\end{align*}
Set $\ul 0 := (0,0,0)$. Consider the following points in $\Bun_D^0(k)$:
\begin{align*}
    P^{0,\dec}_0 &:= \l \mc O(\ul 0), \l \begin{bmatrix}
                       0 \\ 0 \\ 1\\
                      \end{bmatrix},  \begin{bmatrix}
                        0 \\ 1 \\ 0 \\
                        \end{bmatrix}, \begin{bmatrix}
                         1 \\ 0 \\ 0 \\
                          \end{bmatrix}, \lb \begin{pmatrix}
                              1 \\
                              1 \\
                              0
                          \end{pmatrix}, \begin{pmatrix}
                              0 \\
                              0 \\
                              1
                          \end{pmatrix} \rb \r \r \\
    P^{0,\dec}_1 &:= \l \mc O(\ul 0), \l \begin{bmatrix}
                       0 \\ 0 \\ 1\\
                      \end{bmatrix},  \begin{bmatrix}
                        0 \\ 1 \\ 0 \\
                        \end{bmatrix}, \begin{bmatrix}
                         1 \\ 0 \\ 0 \\
                          \end{bmatrix}, \lb \begin{pmatrix}
                              0 \\
                              1 \\
                              1
                          \end{pmatrix}, \begin{pmatrix}
                              1 \\
                              0 \\
                              0
                          \end{pmatrix} \rb \r \r  \\
    P^{0,\dec}_\infty &:= \l \mc O(\ul 0), \l \begin{bmatrix}
                       0 \\ 0 \\ 1\\
                      \end{bmatrix},  \begin{bmatrix}
                        0 \\ 1 \\ 0 \\
                        \end{bmatrix}, \begin{bmatrix}
                         1 \\ 0 \\ 0 \\
                          \end{bmatrix}, \lb \begin{pmatrix}
                              1 \\
                              0 \\
                              1
                          \end{pmatrix}, \begin{pmatrix}
                              0 \\
                              1 \\
                              0
                          \end{pmatrix} \rb \r \r \\
\end{align*}
For $\alpha \in k \setminus \{0,1\}$ set
\begin{align*}
     P^0_\alpha :=  \l \mc O(\ul 0), \l \begin{bmatrix}
                       0 \\ 0 \\ 1\\
                      \end{bmatrix},  \begin{bmatrix}
                        0 \\ 1 \\ 0 \\
                        \end{bmatrix}, \begin{bmatrix}
                         1 \\ 0 \\ 0 \\
                          \end{bmatrix}, \lb \begin{pmatrix}
                              1 \\
                              1 \\
                              1
                          \end{pmatrix}, \begin{pmatrix}
                              0 \\
                              \alpha \\
                              1
                          \end{pmatrix} \rb \r \r  
\end{align*}
Also consider:
\begin{align*}
    P^{0,\vee}_0 &:= \l \mc O(\ul 0), \l \begin{bmatrix}
                       0 \\ 0 \\ 1\\
                      \end{bmatrix},  \begin{bmatrix}
                        0 \\ 1 \\ 0 \\
                        \end{bmatrix}, \begin{bmatrix}
                         1 \\ 0 \\ 0 \\
                          \end{bmatrix}, \lb \begin{pmatrix}
                              1 \\
                              1 \\
                              1
                          \end{pmatrix}, \begin{pmatrix}
                              0 \\
                              0 \\
                              1
                          \end{pmatrix} \rb \r \r \\
    P^{0,\vee}_1 &:= \l \mc O(\ul 0), \l \begin{bmatrix}
                       0 \\ 0 \\ 1\\
                      \end{bmatrix},  \begin{bmatrix}
                        0 \\ 1 \\ 0 \\
                        \end{bmatrix}, \begin{bmatrix}
                         1 \\ 0 \\ 0 \\
                          \end{bmatrix}, \lb \begin{pmatrix}
                              1 \\
                              1 \\
                              1
                          \end{pmatrix}, \begin{pmatrix}
                              1 \\
                              0 \\
                              0
                          \end{pmatrix} \rb \r \r \\
    P^{0,\vee}_\infty &:= \l \mc O(\ul 0), \l \begin{bmatrix}
                       0 \\ 0 \\ 1\\
                      \end{bmatrix},  \begin{bmatrix}
                        0 \\ 1 \\ 0 \\
                        \end{bmatrix}, \begin{bmatrix}
                         1 \\ 0 \\ 0 \\
                          \end{bmatrix}, \lb \begin{pmatrix}
                              1 \\
                              1 \\
                              1
                          \end{pmatrix}, \begin{pmatrix}
                              0 \\
                              1 \\
                              0
                          \end{pmatrix} \rb \r \r \\
     P^{0,\wedge}_0 &:= \l \mc O(\ul 0), \l \begin{bmatrix}
                       0 \\ 0 \\ 1\\
                      \end{bmatrix},  \begin{bmatrix}
                        0 \\ 1 \\ 0 \\
                        \end{bmatrix}, \begin{bmatrix}
                         1 \\ 0 \\ 0 \\
                          \end{bmatrix}, \lb \begin{pmatrix}
                              1 \\
                              1 \\
                              0
                          \end{pmatrix}, \begin{pmatrix}
                              0 \\
                              1 \\
                              1
                          \end{pmatrix} \rb \r \r \\
    P^{0,\wedge}_1 &:= \l \mc O(\ul 0), \l \begin{bmatrix}
                       0 \\ 0 \\ 1\\
                      \end{bmatrix},  \begin{bmatrix}
                        0 \\ 1 \\ 0 \\
                        \end{bmatrix}, \begin{bmatrix}
                         1 \\ 0 \\ 0 \\
                          \end{bmatrix}, \lb \begin{pmatrix}
                              0 \\
                              1 \\
                              1
                          \end{pmatrix}, \begin{pmatrix}
                              1 \\
                              0 \\
                              1
                          \end{pmatrix} \rb \r \r \\
    P^{0,\wedge}_\infty &:= \l \mc O(\ul 0), \l \begin{bmatrix}
                       0 \\ 0 \\ 1\\
                      \end{bmatrix},  \begin{bmatrix}
                        0 \\ 1 \\ 0 \\
                        \end{bmatrix}, \begin{bmatrix}
                         1 \\ 0 \\ 0 \\
                          \end{bmatrix}, \lb \begin{pmatrix}
                              1 \\
                              0 \\
                              1
                          \end{pmatrix}, \begin{pmatrix}
                              0 \\
                              1 \\
                              1
                          \end{pmatrix} \rb \r \r 
\end{align*}
Finally, set
\begin{align*}
    P_{\frac{t-1}{t}}^{0,\spc} := \l \mc O(-1,0,1), \l \begin{bmatrix}
                       0 \\ 0 \\ 1\\
                      \end{bmatrix},  \begin{bmatrix}
                        0 \\ 0 \\ 1 \\
                        \end{bmatrix}, \begin{bmatrix}
                         0 \\ 0 \\ 1 \\
                          \end{bmatrix}, \lb \begin{pmatrix}
                              0 \\
                              1 \\
                              1
                          \end{pmatrix}, \begin{pmatrix}
                              1 \\
                              1 \\
                              0
                          \end{pmatrix} \rb \r \r. 
\end{align*}
Summarizing the results of subsections \ref{Subsection cas (3)}, \ref{subsection cases (1,2) and (2,1)} and \ref{subsection case (1,1,1)} we obtain:

\begin{Theo} \label{Classification Theorem for degree zero}
    The degree $0$ interesting locus $\vert \Bun_D^0(k) \vert^\inter$ defined over $k$ consists of 
    \begin{itemize}
        \item[a)] the points $P_\alpha^0$, $\alpha \in k \setminus \{0,1\}$,
        \item[b)] the points $P_i^{0,\vee}$, $i \in \{0,1,\infty\}$,
        \item[c)] the points $P_i^{0,\wedge}$, $i \in \{0,1,\infty\}$,
        \item[d)] the points $P_i^{0,\dec}$, $i \in \{0,1,\infty\}$,
        \item[e)] the point $P_{\frac{t-1}{t}}^{0,\spc}$.
    \end{itemize}
    Each of these points except the points in d) and e) is contained in $\Bun_D^{0,\mini}(k)$ and has automorphism group isomorphic to $\mb G_m$. The points in d) have automorphism group isomorphic to $\mb G_m \times \mb G_m$.  The point $P_{\frac{t-1}{t}}^{0,\spc}$ has automorphism group isomorphic to $\mb G_m \times \mb G_a$.
\end{Theo}

\subsection{The degree $1$ component}
Consider the following points in $\Bun_D^1(k)$:
\begin{align*}
    P^{1,\dec}_0 &:= \l \mc O(0,0,1), \l \begin{bmatrix}
                       0 \\ 0 \\ 1\\
                      \end{bmatrix},  \begin{bmatrix}
                        0 \\ 1 \\ 0 \\
                        \end{bmatrix}, \begin{bmatrix}
                         0 \\ 1 \\ 0 \\
                          \end{bmatrix}, \lb \begin{pmatrix}
                              0 \\
                              0 \\
                              1
                          \end{pmatrix}, \begin{pmatrix}
                              1 \\
                              1 \\
                              0
                          \end{pmatrix} \rb \r \r \\
    P^{1,\dec}_1 &:= \l \mc O(0,0,1), \l \begin{bmatrix}
                       0 \\ 0 \\ 1\\
                      \end{bmatrix},  \begin{bmatrix}
                        0 \\ 0 \\ 1 \\
                        \end{bmatrix}, \begin{bmatrix}
                         0 \\ 1 \\ 0 \\
                          \end{bmatrix}, \lb \begin{pmatrix}
                              0 \\
                              1 \\
                              0
                          \end{pmatrix}, \begin{pmatrix}
                              1 \\
                              0 \\
                              1
                          \end{pmatrix} \rb \r \r \\
    P^{1,\dec}_\infty &:= \l \mc O(0,0,1), \l \begin{bmatrix}
                       0 \\ 0 \\ 1\\
                      \end{bmatrix},  \begin{bmatrix}
                        0 \\ 1 \\ 0 \\
                        \end{bmatrix}, \begin{bmatrix}
                         0 \\ 0 \\ 1 \\
                          \end{bmatrix}, \lb \begin{pmatrix}
                              0 \\
                              1 \\
                              0
                          \end{pmatrix}, \begin{pmatrix}
                              1 \\
                              0 \\
                              1
                          \end{pmatrix} \rb \r \r  
\end{align*}
For $\alpha \in k \setminus \{0,1,\frac{t-1}{t}\}$ set
\begin{align*}
    P^1_\alpha &:= \l \mc O(0,0,1), \l \begin{bmatrix}
                       0 \\ 0 \\ 1\\
                      \end{bmatrix},  \begin{bmatrix}
                        0 \\ 1 \\ 0 \\
                        \end{bmatrix}, \begin{bmatrix}
                         0 \\ 1 \\ 1 \\
                          \end{bmatrix}, \lb \begin{pmatrix}
                              0 \\
                              \alpha \\
                              1
                          \end{pmatrix}, \begin{pmatrix}
                              1 \\
                              1 \\
                              0
                          \end{pmatrix} \rb \r \r 
\end{align*}
Also consider:
\begin{align*}
    P^{1,\vee}_0 &:= \l \mc O(0,0,1), \l \begin{bmatrix}
                       0 \\ 0 \\ 1\\
                      \end{bmatrix},  \begin{bmatrix}
                        0 \\ 1 \\ 0 \\
                        \end{bmatrix}, \begin{bmatrix}
                         0 \\ 1 \\ 1 \\
                          \end{bmatrix}, \lb \begin{pmatrix}
                              0 \\
                              0 \\
                              1
                          \end{pmatrix}, \begin{pmatrix}
                              1 \\
                              1 \\
                              0
                          \end{pmatrix} \rb \r \r \\
    P^{1,\vee}_1 &:= \l \mc O(0,0,1), \l \begin{bmatrix}
                       0 \\ 0 \\ 1\\
                      \end{bmatrix},  \begin{bmatrix}
                        0 \\ 1 \\ 0 \\
                        \end{bmatrix}, \begin{bmatrix}
                         0 \\ 1 \\ 1 \\
                          \end{bmatrix}, \lb \begin{pmatrix}
                              0 \\
                              1 \\
                              1
                          \end{pmatrix}, \begin{pmatrix}
                              1 \\
                              1 \\
                              0
                          \end{pmatrix} \rb \r \r \\
    P^{1,\vee}_\infty &:= \l \mc O(0,0,1), \l \begin{bmatrix}
                       0 \\ 0 \\ 1\\
                      \end{bmatrix},  \begin{bmatrix}
                        0 \\ 1 \\ 0 \\
                        \end{bmatrix}, \begin{bmatrix}
                         0 \\ 1 \\ 1 \\
                          \end{bmatrix}, \lb \begin{pmatrix}
                              0 \\
                              1 \\
                              0
                          \end{pmatrix}, \begin{pmatrix}
                              1 \\
                              0 \\
                              1
                          \end{pmatrix} \rb \r \r \\
    P^{1,\wedge}_0 &:= \l \mc O(0,0,1), \l \begin{bmatrix}
                       0 \\ 0 \\ 1\\
                      \end{bmatrix},  \begin{bmatrix}
                        0 \\ 1 \\ 0 \\
                        \end{bmatrix}, \begin{bmatrix}
                         0 \\ 1 \\ 0 \\
                          \end{bmatrix}, \lb \begin{pmatrix}
                              0 \\
                              1 \\
                              1
                          \end{pmatrix}, \begin{pmatrix}
                              1 \\
                              1 \\
                              0
                          \end{pmatrix} \rb \r \r \\
    P^{1,\wedge}_1 &:= \l \mc O(0,0,1), \l \begin{bmatrix}
                       0 \\ 0 \\ 1\\
                      \end{bmatrix},  \begin{bmatrix}
                        0 \\ 0 \\ 1 \\
                        \end{bmatrix}, \begin{bmatrix}
                         0 \\ 1 \\ 0 \\
                          \end{bmatrix}, \lb \begin{pmatrix}
                              0 \\
                              1 \\
                              1
                          \end{pmatrix}, \begin{pmatrix}
                              1 \\
                              1 \\
                              0
                          \end{pmatrix} \rb \r \r \\
    P^{1,\wedge}_\infty &:= \l \mc O(0,0,1), \l \begin{bmatrix}
                       0 \\ 0 \\ 1\\
                      \end{bmatrix},  \begin{bmatrix}
                        0 \\ 1 \\ 0 \\
                        \end{bmatrix}, \begin{bmatrix}
                         0 \\ 0 \\ 1 \\
                          \end{bmatrix}, \lb \begin{pmatrix}
                              0 \\
                              1 \\
                              1
                          \end{pmatrix}, \begin{pmatrix}
                              1 \\
                              1 \\
                              0
                          \end{pmatrix} \rb \r \r 
\end{align*}
Finally, set
\begin{align*}
    P_{\frac{t-1}{t}}^1 &:= \l \mc O(0,0,1), \l \begin{bmatrix}
                       0 \\ 0 \\ 1\\
                      \end{bmatrix},  \begin{bmatrix}
                        0 \\ 1 \\ 0 \\
                        \end{bmatrix}, \begin{bmatrix}
                         0 \\ 1 \\ 1 \\
                          \end{bmatrix}, \lb \begin{pmatrix}
                              1 \\
                              \frac{t-1}{t} \\
                              1
                          \end{pmatrix}, \begin{pmatrix}
                              0 \\
                              1 \\
                              0
                          \end{pmatrix} \rb \r \r \\
    P_{\frac{t-1}{t}}^{1,\spc} &:= \l \mc O(0,0,1), \l \begin{bmatrix}
                       0 \\ 0 \\ 1\\
                      \end{bmatrix},  \begin{bmatrix}
                        0 \\ 1 \\ 0 \\
                        \end{bmatrix}, \begin{bmatrix}
                         0 \\ 1 \\ 1 \\
                          \end{bmatrix}, \lb \begin{pmatrix}
                              1 \\
                              \frac{t-1}{t} \\
                              1
                          \end{pmatrix}, \begin{pmatrix}
                              1 \\
                              0 \\
                              0
                          \end{pmatrix} \rb \r \r 
\end{align*}

We have the following classification result for points in the interesting locus of degree $1$:

\begin{Theo} \label{Classification Theorem for degree one}
    The degree $1$ interesting locus $\vert \Bun_D^1(k) \vert^\inter$ defined over $k$ consists of 
    \begin{itemize}
        \item[a)] the points $P_\alpha^1$, $\alpha \in k \setminus \{0,1\}$,
        \item[b)] the points $P_i^{1,\vee}$, $i \in \{0,1,\infty\}$,
        \item[c)] the points $P_i^{1,\wedge}$, $i \in \{0,1,\infty\}$,
        \item[d)] the points $P_i^{1,\dec}$, $i \in \{0,1,\infty\}$,
        \item[e)] the point $P_{\frac{t-1}{t}}^{1,\spc}$.
    \end{itemize}
    Each of these points except the points in d) and e) is contained in $\Bun_D^{1,\mini}(k)$ and has automorphism group isomorphic to $\mb G_m$. The points in d) have automorphism group isomorphic to $\mb G_m \times \mb G_m$. The point $P_{\frac{t-1}{t}}^{1,\spc}$ in e) has automorphism group isomorphic to $\mb G_m \times \mb G_a$.
\end{Theo}

Theorem \ref{Classification Theorem for degree one} is proved by applying the elementary Hecke operator at $\infty$ to the points listed above and checking that one obtains precisely the points in the list in Theorem \ref{Classification Theorem for degree zero}. The Theorem then follows from Theorem \ref{Proposition interesting locus and hecke operators}. We illustrate the calculation in the case of $P^{1,\dec}_i$ in Lemma \ref{Theorem Hecke operators of decomposable relevant points} below and leave all other cases to the reader.

Set $f_\infty := X, g_\infty := Y \in \Hom(\mc O,\mc O(1)) = \Gamma(\mb P^1,\mc O(1))$. Then
\begin{align*}
    f_\infty(\infty) &= 1, \\
    g_\infty(\infty) &= 0, \\
    g_\infty(0) &= g_\infty(1) = g_\infty(t) = 1, \\
    f_\infty(0) &= 0, \\
    f_\infty(1) &= 1, \\
    f_\infty(t) &= t.
\end{align*}

\begin{Lem} \label{Theorem Hecke operators of decomposable relevant points}
    For $i \in \{0,1,\infty\}$, $\mb T^{(l)}_4$ maps $P^{1,\dec}_i$ to $P^{0,\dec}_i$.
\end{Lem}
\begin{proof}
   Consider the diagram $\mc E_4\b$, defined as  \\
\[
\begin{tikzcd}[column sep=huge, ampersand replacement=\&]
\mc O \oplus \mc O(-1) \oplus \mc O(-1)
\arrow[r, "{\left(\begin{smallmatrix}
1 & 0 & -X \\
0 & Y & 0 \\
0 & 0 & 1
\end{smallmatrix}\right)}"]
\&
\mc O \oplus \mc O \oplus \mc O(-1)
\arrow[r, "{\left(\begin{smallmatrix}
1 & 0 & 0 \\
0 & 1 & 0 \\
0 & 0 & Y
\end{smallmatrix}\right)}"]
\&
\mc O \oplus \mc O \oplus \mc O
\arrow[r, "{\left(\begin{smallmatrix}
Y & 0 & X \\
0 & 1 & 0 \\
0 & 0 & 1
\end{smallmatrix}\right)}"]
\&
\mc O(1) \oplus \mc O \oplus \mc O.
\end{tikzcd}
\]
Then for $i \in \{1,\infty\}$ there are $\mc E_{i,1}\b, \mc E_{i,2}\b, \mc E_{i,3}\b$ such that 
\begin{align*}
    P^{1,\dec}_i = \l \mc E_{i,1}\b, \mc E_{i,2}\b, \mc E_{i,3}\b, \mc E_4\b \r.
\end{align*}
Since the third term in the diagram is $\mc O^{\oplus 3}$, we have $\UL(\mb T^{(l)}(P^{1,\dec}_i)) = \mc O^{\oplus 3}$. We now consider the case $i = 1$ and determine the parabolic structure on $\mc O^{\oplus 3}$ defined by $\mb T^{(l)}(P^{1,\dec}_1)$. To see what it is at the point $\infty \in S$, we calculate
\begin{align*}
    \begin{pmatrix}
            1 & &  \\
            & 1 & \\
            & & Y
        \end{pmatrix}  \begin{pmatrix}
            1 & & -X \\
            & Y & \\
            & & 1
        \end{pmatrix} = \begin{pmatrix}
            1 & & -X \\
            & Y & \\
            & & Y
        \end{pmatrix}.
\end{align*}
The induced morphism on the fiber at $\infty$ has image spanned by
\begin{align*}
    \begin{pmatrix}
        1 \\
        0 \\
        0
    \end{pmatrix}.
\end{align*}
Similarly, the morphism induced on the fiber at $\infty$ by
\begin{align*}
     \begin{pmatrix}
            1 & &  \\
            & 1 & \\
            & & Y
        \end{pmatrix}
\end{align*}
as image spanned by
\begin{align*}
    \begin{pmatrix}
        1 \\
        0 \\
        0
    \end{pmatrix}, \begin{pmatrix}
        0 \\
        1 \\
        0
    \end{pmatrix}.
\end{align*}
So we obtain the flag
\begin{align*}
    \lb \begin{pmatrix}
        1 \\
        0 \\
        0
    \end{pmatrix}, \begin{pmatrix}
        0 \\
        1 \\
        0
    \end{pmatrix} \rb \in \FL_3(k)
\end{align*}
at $\infty$. \\
To calculate the structure at the points in $S \setminus \{\infty\}$, we have to solve the equations
\begin{align*}
    A_1 v_1 =
        \begin{bmatrix}
            0 \\
            0 \\
            1
        \end{bmatrix},\  A_2 v_2 =
        \begin{bmatrix}
            0 \\
            0 \\
            1
        \end{bmatrix}, \ A_3v_3 = 
         \begin{bmatrix}
            0 \\
            1 \\
            0
        \end{bmatrix},
\end{align*}
where $A_i$ denotes the image of 
\begin{align*}
    \begin{pmatrix}
            Y & & X \\
            & 1 & \\
            & & 1 
    \end{pmatrix}
\end{align*}
in $\Hom(\mc O^{\oplus 3}, \mc O(0,0,1)) \otimes k(x_i)$. Note that this equation is solvable: Since we consider modifications at $\infty$, $A_i$ is bijective for $i \in \{1,2,3\}$. We obtain
\begin{align*}
    v_1 = 
    \begin{bmatrix}
        0 \\
        0 \\
        1
    \end{bmatrix},\ v_2 = \begin{bmatrix}
        -1 \\
        0 \\
        1
    \end{bmatrix},\ v_3 =  \begin{bmatrix}
        0 \\
        1 \\
        0
    \end{bmatrix}.
\end{align*}
Hence,
\begin{align*}
    \mb T^{(l)}_4(P_1^{1,\dec}) = \l\mc O^{\oplus 3}, \l \begin{bmatrix}
        0 \\
        0 \\
        1
    \end{bmatrix}, \begin{bmatrix}
        -1 \\
        0 \\
        1
    \end{bmatrix},  \begin{bmatrix}
        0 \\
        1 \\
        0
    \end{bmatrix},  \lb \begin{pmatrix}
        1 \\
        0 \\
        0
    \end{pmatrix}, \begin{pmatrix}
        0 \\
        1 \\
        0
    \end{pmatrix} \rb \r \r,
\end{align*}
which is isomorphic to $P_1^{0,\dec}$ via
\begin{align*}
    \begin{pmatrix}
        0 & 1 & 0 \\
        1 & 0 & 0 \\
        1 & 0 & 1
    \end{pmatrix} \in \Aut(\mc O(\ul 0)).
\end{align*}
For $i = \infty$, a similar calculation yields
\begin{align*}
    \mb T^{(l)}_4(P_\infty^{1,\dec}) = \l\mc O^{\oplus 3}, \l \begin{bmatrix}
        0 \\
        0 \\
        1
    \end{bmatrix}, \begin{bmatrix}
        0 \\
        1 \\
        0
    \end{bmatrix},  \begin{bmatrix}
        -t \\
        0 \\
        1
    \end{bmatrix},  \lb \begin{pmatrix}
        1 \\
        0 \\
        0
    \end{pmatrix}, \begin{pmatrix}
        0 \\
        1 \\
        0
    \end{pmatrix} \rb \r \r,
\end{align*}
which is isomorphic to $P_\infty^{0,\dec}$ via
\begin{align*}
    \begin{pmatrix}
        1 & & \\
        & t & \\
        1 & & t
    \end{pmatrix} \in \Aut(\mc O(\ul 0)).
\end{align*}
Finally, for $P_0^{1,\dec}$ we consider \\
\[
\begin{tikzcd}[column sep=huge, ampersand replacement=\&]
\mc O \oplus \mc O(-1) \oplus \mc O(-1)
\arrow[r, "{\left(\begin{smallmatrix}
1 & -X &  \\
0 & 1 & 0 \\
0 & 0 & Y
\end{smallmatrix}\right)}"]
\&
\mc O \oplus \mc O(-1) \oplus \mc O
\arrow[r, "{\left(\begin{smallmatrix}
1 & 0 & 0 \\
0 & Y & 0 \\
0 & 0 & 1
\end{smallmatrix}\right)}"]
\&
\mc O \oplus \mc O \oplus \mc O
\arrow[r, "{\left(\begin{smallmatrix}
Y & X & 0 \\
0 & 1 & 0 \\
0 & 0 & 1
\end{smallmatrix}\right)}"]
\&
\mc O(1) \oplus \mc O \oplus \mc O
\end{tikzcd}
\] \\
and obtain
\begin{align*}
     \mb T^{(l)}_4(P_0^{1,\dec}) = \l\mc O^{\oplus 3}, \l \begin{bmatrix}
        0 \\
        0 \\
        1
    \end{bmatrix}, \begin{bmatrix}
        -1 \\
        1 \\
        0
    \end{bmatrix},  \begin{bmatrix}
        -t \\
        1 \\
        0
    \end{bmatrix},  \lb \begin{pmatrix}
        1 \\
        0 \\
        0
    \end{pmatrix}, \begin{pmatrix}
        0 \\
        0 \\
        1
    \end{pmatrix} \rb \r \r,
\end{align*}
which is isomorphic to $P_0^{0,\dec}$ via
\begin{align*}
    \begin{pmatrix}
        1 & 1 &  \\
        1 & t &  \\
         &  & 1
    \end{pmatrix} \in \Aut(\mc O(\ul 0)).
\end{align*}
\end{proof}

%% file: 5_Geometry.tex
\section{Geometry: The example of \texorpdfstring{$\mb P^1$}{} minus 4 points}  \label{Last section}
In this section we continue the study of the example $X = \mb P^1$, $S^\oro = (0,1,t,\infty)$, where $t \in \mb P^1(k) \setminus \{0,1,\infty\}$, $D^\tw = (D,\Lambda)$ with 
\begin{align*}
    D &:= ((1,2),(1,2),(1,2),(1,1,1)), \\
    \Lambda &:= ((\lambda,1),(\lambda,1),(\lambda, 1), (\lambda^{-1},\lambda^{-1},\lambda^{-1})),
\end{align*}
which was also considered in the last section. We now focus on the geometry of $\vert \Bun_D^d \vert^\inter$. On the one hand, as we have seen in Theorem \ref{Proposition interesting locus and hecke operators}, this topological space is independent of $d$. In particular, in proving a statement about $\vert \Bun_D^d \vert^\inter$ we may choose $d$ adapted to the problem we study. \\
On the other hand, for $d_1 \neq d_2$ the following  is not always true (see Proposition \ref{First Proposition about Geometry}i) for an example): There exist open subsets $U_i \subset \vert \Bun_{\GL_3}^{d_i} \vert$, $i = 1,2$, and a homeomorphism $\rho : U_1 \simeq U_2$ s.t. the map $h_i : \vert \Bun_D^{d_i} \vert^\inter \to \vert \Bun_{\GL_3}^{d_i} \vert$ induced by the morphism of stacks $\Bun_D^{d_i} \to \Bun_{\GL_3}^{d_i}$ factors through $U_i$ and s.t. $\rho$ fits into a commutative diagram \\
\begin{tikzcd}
\vert \Bun_D^{d_1} \vert^\inter \arrow[rr, "\mb T_4^{(l)}"] \arrow[d, "h_1"'] &  & \vert \Bun_D^{d_2} \vert^\inter \arrow[d, "h_2"] \\
U_1 \arrow[rr, "\rho"']                                                       &  & U_2.                                             
\end{tikzcd} \\
Nevertheless, if $d_2 \equiv d_1 \mod 3$, then such a property holds. This is why in Proposition \ref{First Proposition about Geometry} we only consider $d \in \{0,1,2\}$. We denote by $\Bun_D^{\ul 0, (-1,0,1)}$ the open substack of $\Bun_D$ consisting of parabolic bundles with underlying bundle isomorphic to $\mc O(\ul 0)$ or $\mc O(-1,0,1)$.

\begin{Prop} \label{First Proposition about Geometry}
    The following are true:
    \begin{itemize}
        \item[i)] For $d \in \{0,1\}$ we have:
        \begin{align*}
            \vert \Bun_D^d \vert^\inter \subset \begin{cases}
                \vert \Bun_D^{\ul 0, (-1,0,1)} \vert, & d = 0 \\
                \vert \Bun_D^{(0,0,1)} \vert, & d = 1. 
            \end{cases}
        \end{align*}
        \item[ii)] For $d = 0,1$ the point $P_{\frac{t-1}{t}}^{d,\spc}$ is closed inside $\vert \Bun_D^d \vert^\inter$.
        \item[iii)] For $d = 0,1$ the set 
        \begin{align*}
            \vert \Bun_D^d \vert^\inter \setminus \{P_0^{d,\dec}, P_1^{d,\dec}, P_\infty^{d,\dec},P_{\frac{t-1}{t}}^{d,\spc}\}
        \end{align*}
        is open in $\vert \Bun_D^d \vert^\inter$.
    \end{itemize}
\end{Prop}
\begin{proof}
   Part i) follows immediately from the classification theorems \ref{Classification Theorem for degree zero} and \ref{Classification Theorem for degree one}. As explained before, Theorem \ref{Proposition interesting locus and hecke operators} implies that it suffices to show ii) for $d = 0$. In this case it follows immediately from i) and the continuity of the natural morphism $\Bun_D^0 \to \Bun_{\GL_3}^0$. Finally, according to theorems \ref{Classification Theorem for degree zero} and \ref{Classification Theorem for degree one} the set in 3) is the topological spaces of $\Bun_D^{d,\mini}$, which is open in $\Bun_D^d$.
\end{proof}

\subsection{Degree $0$ locus: Charts.}
In this section we introduce charts for the stack $\Bun_D^\internull$ corresponding to the open subset $\vert \Bun_D^0 \vert^\inter \setminus \{P_{\frac{t-1}{t}}^{0,\spc}\}$ of $\vert \Bun_D^0 \vert$ (cf. Proposition \ref{First Proposition about Geometry}). 

By Proposition \ref{First Proposition about Geometry}, $\vert \Bun_D^0 \vert^\inter \setminus \{P_{\frac{t-1}{t}}^{0,\spc}\}$ is contained in $\vert \Bun_D^{\ul 0} \vert$.  Using the map $\Bun_D \to \Bun_{\GL_3}$ we see that $\Bun_D^{\ul 0}$ is an open substack of $\Bun_D^0$. Since $\ul{\Aut}(\mc O(\ul 0)) \simeq \GL_3$,
\begin{align*}
    \Bun_D^{\ul 0} \simeq [\Pi/\GL_3],
\end{align*}
where $\GL_3$ acts diagonally on $\Pi = (\mb P^2)^{\times 3} \times \FL_3$. Let $\Pi^{\ul 0} \subset \Pi$ be the pullback under the projection $\Pi \to (\mb P^2)^{\times 3}$ of the open subscheme of $(\mb P^2)^{\times 3}$ consisting of those $(l_1,l_2,l_3)$, which are linearly independent. By \ref{Classification Theorem for degree zero}, $\vert \Bun_D^0 \vert^\inter \setminus \{P_{\frac{t-1}{t}}^{0,\spc}\}$ is contained in the open substack $\mf X_1^{\ul 0}$ of $\Bun_D^{\ul 0}$ corresponding to $\Pi^{\ul 0}$. The following easy Lemma will be used more than once:

\begin{Lem} \label{Lemma about uniformaization of quotient stacks}
    Let $X, Y$ be $k$-schemes and $G$ a $k$-group scheme acting on $X$ and $Y$. Endow $X \times Y$ with the diagonal $G$-action. Assume that $X$ is smooth and let $x_0 \in X(k)$. Assume that the orbit map 
    \begin{align*}
        G &\to X \\
        g &\mapsto g.x_0
    \end{align*}
    is smooth and surjective. Then:
    \begin{align*}
        [(X \times Y)/G] \simeq [Y/\Stab_{G}(x_0)]
    \end{align*}
\end{Lem}

Since any point in $\Pi^{\ul 0}(k)$ is isomorphic under $\GL_3(k)$ to a point of the form
\begin{align*}
    \l x_0, w \r
\end{align*}
for 
\begin{align*}
    x_0 := \l \lsb 0 \\ 0 \\ 1 \rsb, \lsb 0 \\ 1 \\ 0 \rsb, \lsb 1 \\ 0 \\ 0 \rsb \r \in (\mb P^2)^{\times 3}(k)
\end{align*}
and some $w \in \FL_3(k)$, the same is true after any field extension of $k$ and we have
\begin{align*}
    \Stab_{\GL_3}(x_0) \simeq \mb G_m^{\times 3}
\end{align*}
(the torus consisting of diagonal matrices in $\GL_3$), we obtain from Lemma \ref{Lemma about uniformaization of quotient stacks}
\begin{align*}
    \lb \Pi^{\ul 0}/\GL_3 \rb \simeq \lb \FL_3/ \mb G_m^{\times 3} \rb.
\end{align*}
For $i,j \in \{1,2,3\}$ with $i \neq j$ let $Y_{ij}$ resp. $Z_{ij}$ be the closed reduced subscheme of $\FL_3$ given by 
\begin{align*}
    \lb \begin{pmatrix}
        a_1 \\ a_2 \\ a_3 
    \end{pmatrix}, \begin{pmatrix}
        b_1 \\ b_2 \\ b_3
    \end{pmatrix} \rb
\end{align*}
s.t. $a_i = a_j = 0$ resp. $a_i = a_ib_j-a_jb_i = 0$. Let $\FL_3^\internull \subset \FL_3$ be the open complement of the union $\bigcup_{i \neq j} (Y_{ij} \cup Z_{ij})$. It is clearly $\mb G_m^{\times 3}$-invariant.

\begin{Prop} \label{Proposition theta zero}
    The map
    \begin{align*}
        \theta^0 : \FL_3^\internull &\to \Pi \\
         \lb \begin{pmatrix}
        a_1 \\ a_2 \\ a_3 
    \end{pmatrix}, \begin{pmatrix}
        b_1 \\ b_2 \\ b_3
    \end{pmatrix} \rb &\mapsto   \l \begin{bmatrix}
       0 \\
       0 \\
       1
    \end{bmatrix}, \begin{bmatrix}
        0 \\
        1 \\
        0
    \end{bmatrix},  \begin{bmatrix}
        1 \\
        0 \\
        0
    \end{bmatrix},  \lb \begin{pmatrix}
        a_1 \\
        a_2 \\
        a_3
    \end{pmatrix}, \begin{pmatrix}
        b_1 \\
        b_2 \\
        b_3
    \end{pmatrix} \rb \r
    \end{align*}
    induces an isomorphism 
    \begin{align*}
        \lb \FL_3^\internull / \mb G_m^{\times 3} \rb \simeq \Bun_D^\internull.
    \end{align*}
\end{Prop}
\begin{proof}
    This follows immediately from the above discussion and Theorem \ref{Classification Theorem for degree zero}.
\end{proof}

\subsection{Degree 0 locus: Symmetry} \label{Subsection on Symmetry}
Denote by $\tau_{12}, \tau_{13}, \tau_{23}$ the transpositions in the symmetric group $\mf S_3 = \Aut(\{1,2,3\})$. The assignment $\{\tau_{12}, \tau_{13}, \tau_{23}\} \to \PGL_2(k)$ given by
\begin{align*}
    \tau_{12} &\mapsto \begin{bmatrix}
        1 & 0 \\
        1 & -1 
    \end{bmatrix} \\
     \tau_{13} &\mapsto \begin{bmatrix}
        -1 & 1 \\
        0 & 1 
    \end{bmatrix} \\
     \tau_{23} &\mapsto \begin{bmatrix}
        0 & 1 \\
        1 & 0 
    \end{bmatrix}
\end{align*}
extends uniquely to a representation
\begin{align*}
    \rho : \mf S_3 \to \PGL_2(k).
\end{align*}
We let $\PGL_2(k)$ act on $\mb P^1$ via Möbius transformations. Then we let $\mf S_3$ act on $\mb P^1$ through $\rho$. 

Explicitly, the action is given in affine coordinates $x = [x:1]$ of $\mb P^1$ by
\begin{align*}
    \tau_{12} : x &\mapsto \frac{x}{x-1} \\
    \tau_{13} : x &\mapsto 1-x \\
    \tau_{23} : x &\mapsto \frac{1}{x}. 
\end{align*}
In particular, it preserves the subset $\{0,1,\infty\} \subset \mb P^1(k)$.

Also, permutation of coordinates, i.e.
\begin{align*}
    \sigma. \lb \begin{pmatrix}
        a_1 \\
        a_2 \\
        a_3
    \end{pmatrix}, \begin{pmatrix}
        b_1 \\
        b_2 \\
        b_3
    \end{pmatrix} \rb = \lb \begin{pmatrix}
        a_{\sigma(1)} \\
        a_{\sigma(2)} \\
        a_{\sigma(3)}
    \end{pmatrix}, \begin{pmatrix}
        b_{\sigma(1)} \\
        b_{\sigma(2)} \\
        b_{\sigma(3)}
    \end{pmatrix} \rb
\end{align*}
defines a right action of $\mf S_3$ on $\FL_3$ (for which we still use left action notation $\sigma.\ast$ to avoid awkwardness). Similarly, we let $\mf S_3$ act (from the right) on $\mb G_m^{\times 3}$ by permutation of coordinates. We denote these right actions by $\rho_{\FL_3}$ resp. $\rho_{\mb G_m^{\times 3}}$. A direct calculation shows:

\begin{Prop} \label{Action on Buninter}
    $\FL_3^\internull$ is $\mf S_3$-stable. Denoting the action of $\mb G_m^{\times 3}$ on $\FL_3$ by $\act$, there is for each $\sigma \in \mf S_3$ an equality
    \begin{align*}
        \rho_{\FL_3}(\sigma) \circ \act = \act \circ (\rho_{\mb G_m^{\times 3}}(\sigma), \rho_{\FL_3}(\sigma))
    \end{align*}
    of maps $\mb G_m^{\times 3} \times \FL_3 \to \FL_3$. In particular, $\rho_{\FL_3}$ induces an action of $\mf S_3$ on $\Bun_D^\internull$.
\end{Prop}

\subsection{Degree $0$ locus: Map to $\mb P^1$}
Consider the morphism
\begin{align*}
    d^0 : \FL_3^\internull &\rightarrow \mb P^1 \\
   \lb \begin{pmatrix}
        a_1 \\ a_2 \\ a_3 
    \end{pmatrix}, \begin{pmatrix}
        b_1 \\ b_2 \\ b_3
    \end{pmatrix} \rb &\mapsto [a_3(a_1b_2 - a_2b_1) : a_2(a_1b_3 - a_3b_1)].
\end{align*}
The definition of $\FL_3^\internull$ and the calculation
\begin{align*}
   \lb \begin{pmatrix}
        a_1 \\ a_2 \\ a_3 
    \end{pmatrix}, \begin{pmatrix}
        b_1 \\ b_2 \\ b_3
    \end{pmatrix} + \lambda \begin{pmatrix}
        a_1 \\ a_2 \\ a_3 
    \end{pmatrix} \rb &\mapsto [a_3(a_1(b_2+\lambda a_2) - a_2(b_1+\lambda a_1)) : a_2(a_1(b_3+ \lambda a_3) - a_3(b_1+ \lambda a_1)))] \\
    &= [a_3(a_1b_2 - a_2b_1) : a_2(a_1b_3 - a_3b_1)]
\end{align*}
imply that it is well-defined.

\begin{Lem} \label{Lemma smooth locus} 
    $d^0$ is invariant under the $\mb G_m^{\times 3}$-action on $\FL_3^{\internull}$, so induces a morphism
    \begin{align*}
        \mf d^0 : \Bun_D^\internull \underset{\ref{Proposition theta zero}}{\simeq} \lb \FL_3^\internull / \mb G_m^{\times 3} \rb \to \mb P^1.
    \end{align*}
    The singular locus of $d^0$ is the disjoint union (in $\FL_3^\internull$)
    \begin{align*}
         \{a_1 = 0, \ a_2b_3-a_3b_2 = 0\} \sqcup \{a_2 = 0, \ a_1b_3-a_3b_1 = 0\} \sqcup \{a_3 = 0, \ a_1b_2-a_2b_1 = 0\}.
    \end{align*}
    In particular, $\mf d^0$ is smooth away from the closed locus $\left\{P_i^{0,\dec} \ \vert \ i \in \{0,1,\infty\}\right\}$.
\end{Lem}
\begin{proof}
    The preimage of $\{[x:y] \ \vert \ y \neq 0\} \simeq \mb A^1$ under $d^0$ is $U := \{a_2 \neq 0, \ a_1b_3 - a_3b_1 \neq 0\}$. Then $f \vert_U$ is given by 
    \begin{align*}
        F : U &\to \mb A^1 \\
       \lb \begin{pmatrix}
        a_1 \\ a_2 \\ a_3 
    \end{pmatrix}, \begin{pmatrix}
        b_1 \\ b_2 \\ b_3
    \end{pmatrix} \rb &\mapsto \frac{a_3(a_1b_2 - a_2b_1)}{a_2(a_1b_3 - a_3b_1)}.
    \end{align*}
    The derivatives of $F$ in the variables $a_1,a_2,a_3,b_1,b_2, b_3$ are given by
\begin{align}
\frac{a_3(a_2b_1b_3-a_3b_1b_2)}{a_2(a_1b_3-a_3b_1)^2} \\
-\frac{a_1a_3b_2}{a_2^2(a_1b_3-a_3b_1)} \\
\frac{a_1b_3(a_1b_2-a_2b_1)}{a_2(a_1b_3-a_3b_1)^2} \\
-\frac{a_3(a_1a_2b_3-a_1a_3b_2)}{a_2(a_1b_3-a_3b_1)^2} \\
\frac{a_1a_3}{a_2(a_1b_3-a_3b_1)} \\
-\frac{a_1a_3(a_1b_2-a_2b_1)}{a_2(a_1b_3-a_3b_1)^2}
\end{align}
It is easy to see that these all vanish on 
\begin{align*}
         \{a_1 = 0, \ a_2b_3-a_3b_2 = 0\} \sqcup \{a_3 = 0, \ a_1b_2-a_2b_1 = 0\}.
    \end{align*}
Conversely, if they all vanish, we obtain from equation (5.5) that $a_1a_3 = 0$. Assume $a_1 = 0$. Then $a_3 \neq 0$ and $b_1 \neq 0$ by definition of $\FL_3^\internull$. Hence, by equation (5.1), $a_2b_3-a_3b_2 = 0$. Similarly, if $a_3 = 0$, then $a_1b_2-a_2b_1 = 0$ by equation (5.3). Since 
\begin{align*}
    \{a_2 = 0, \ a_1b_3-a_3b_1 = 0\} \cap U = \emptyset,
\end{align*}
this completes the calculation over $\{[x:y] \ \vert \ y \neq 0\}$. The calculation for $\{[x:y] \ \vert \ x \neq 0\}$ is similar.
\end{proof}

\begin{Prop} \label{Equivariance of map to Pone}
    $\mf d^0$ is equivariant for the action of $\mf S_3$ on $\Bun_D^\internull$ constructed in Proposition \ref{Action on Buninter} and the action on $\mb P^1$ via Möbius transformations (cf. subsection \ref{Subsection on Symmetry}).
\end{Prop}
\begin{proof}
    Direct calculation.
\end{proof}

\subsection{Minimal locus: Coarse Moduli} \label{subsection on coarse moduli}
Let $\FL_3^\minnull \subset \FL_3^\internull$ denote the smooth locus of $d^0$, i.e. $\FL_3^\minnull$ is the complement in $\FL_3^\internull$ of 
\begin{align} \label{disjoint union describing decomposable bundles}
    \{a_1 = a_2b_3-a_3b_2 = 0\} \sqcup \{a_2 =  a_1b_3-a_3b_1 = 0\} \sqcup \{a_3 = a_1b_2-a_2b_1 = 0\}.
\end{align}

\begin{Lem}
    The open substack of $\Bun_D^\internull$ corresponding to $\FL_3^\minnull$ is $\Bun_D^\minnull$.
\end{Lem}
\begin{proof}
    According to Theorem \ref{Classification Theorem for degree zero} there are only three points lying in $\vert \Bun_D^\internull \vert$, but not in $\Bun_D^\mini$, namely $P_i^{0,\dec}$, $i = 0,1,\infty$. One easily checks that each point corresponds to one of the loci in the disjoint union (\ref{disjoint union describing decomposable bundles}).
\end{proof}

We also need the following Lemma:

\begin{Lem} \label{Closedness of Points}
    Each point of $\vert \Bun_D^\minnull(k) \vert$ is closed inside $\Bun_D^\minnull$.
\end{Lem}
\begin{proof}
    One easily checks for each $p  \in \mb P^1 \setminus \{0,1,\infty\}$ the preimage $(d^0)^{-1}$ is a single $\mb G_m^{\times 3}$-orbit. Hence, over $\mb P^1 \setminus \{0,1,\infty\}$, $\mf d^0$ is a bijection on topological spaces, so that each point in $\l\mf d^0\r^{-1}\l\mb P^1(k) \setminus \{0,1,\infty\}\r$ is closed inside $\Bun_D^\internull$ by continuity. Each such point also lies inside $\Bun_D^\minnull$ and is thus closed in $\Bun_D^\minnull$. \\
    For $p = 0$ one checks that $\l d^0 \r^{-1}(0)$ is the closed subscheme in $\FL_3^\internull$ given by $a_3(a_1b_2-a_2b_1) = 0$. This locus consists of three $\mb G_m^{\times 3}$-orbits: 
    \begin{align} \label{loci in proof of closedness of points}
        \{a_3 &= a_1b_2-a_2b_1 = 0\} \\
        \{a_3 &= 0, \ a_1b_2-a_2b_1 \neq 0\}, \\
        \{a_3 &\neq 0, \ a_1b_2-a_2b_1 = 0\},
    \end{align} 
    which correspond (in this order) to the points $P_0^{0,\dec}$, $P_0^{0,\wedge}$ and $P_0^{0,\vee}$ respectively. Of these points only $P_0^{0,\wedge}$ and $P_0^{0,\vee}$ lie in $\Bun_D^\mini$. Since the loci (5.9) and (5.10) are disjoint, we obtain 
    \begin{align*}
        \l \mf d^0 \r^{-1}(0) \cap \Bun_D^\minnull \simeq B\mb G_m \sqcup B\mb G_m,
    \end{align*}
    which implies the claim for $p = 0$. The cases $p = 1,\infty$ follow from $\mf S_3$-symmetry.
\end{proof}

Lemma \ref{Closedness of Points} implies that there is an open substack $\Bun_D^\gennull$ of $\Bun_D^\minnull$ such that
\begin{align*}
    \vert \Bun_D^\gennull \vert = \vert \Bun_D^\minnull \vert \setminus \left\{P_i^{0,\vee}, P_i^{0,\wedge} \ \vert \ i \in \{0,1,\infty\}\right\}
\end{align*}
and further that for any subset $R \subset \{0,1,\infty\}$ there is a unique open substack $\Bun_D^\Rnull$ of $\Bun_{D}^\minnull$ such that $\Bun_D^\gennull \subset \Bun_D^\Rnull$ and such that for each $i \in \{0,1,\infty\}$ we have
\begin{align*}
\vert \Bun_D^\Rnull \vert \cap \l\mf d^0\r^{-1}(i) = \begin{cases}
    \left \{P_i^{0,\vee} \right\}, & i \in R  \\
    \left\{ P_i^{0,\wedge} \right\}, & i \notin R.
\end{cases}
\end{align*}
It follows from Theorem \ref{Classification Theorem for degree zero} that the collection of stacks $\Bun_D^\Rnull$, $R \subset \{0,1,\infty\}$, covers $\Bun_D^\minnull$. Recall the action of $\mf S_3$ on $\{0,1,\infty\}$ constructed in subsection \ref{Subsection on Symmetry}. It induces an action $\mf S_3 \to \Aut(\f{P}(\{0,1,\infty\}))$ on the power set.

\begin{Lem}
    Let $R \subset \{0,1,\infty\}$ and $\sigma \in \mf S_3$. Then the action of $\sigma$ on $\Bun_D^\minnull$ induces an isomorphism
    \begin{align*}
        \Bun_D^\Rnull \simeq \Bun_D^{0,\sigma(R)}.
    \end{align*}
\end{Lem}
\begin{proof}
    This follows from Theorem \ref{Classification Theorem for degree zero} and Proposition \ref{Equivariance of map to Pone}.
\end{proof}

In particular, the substacks $\Bun_D^{0,\{0,1,\infty\}}$ and $\Bun_D^{0,\emptyset}$ are $\mf S_3$-invariant.

Consider the projective line $\mc P\p$ with double points at $0$, $1$ and $\infty$ i.e.
\begin{align*}
    \mc P\p := \mb P^1 \sqcup_{\mb P^1 \setminus \{0,1,\infty\}} \mb P^1,
\end{align*}
where both maps $\mb P^1 \setminus \{0,1,\infty\} \to \mb P^1$ are the natural inclusion. Let $p\p : \mc P\p \to \mb P^1$ be the natural map induced by $\id_{\mb P^1} \sqcup \id_{\mb P^1}$.

\begin{Prop} \label{Proposition coarse moduli}
    For any section $s\p : \mb P^1 \to \mc P\p$ of $p\p$ there is a unique morphism
    \begin{align*}
        \mf c^0_{s\p} : \Bun_D^\minnull \to \mc P\p
    \end{align*}
    s.t.
    \begin{itemize}
        \item[i)] the diagram \\
        \begin{tikzcd}
\Bun_D^\minnull \arrow[r, "\mf c_{s\p}^0"] \arrow[rd, "\mf d^0 \vert"'] & \mc P\p \arrow[d, "p\p"] \\
                                                                                & \mb P^1                 
\end{tikzcd} \\
        commutes and
        \item[ii)] for each $i \in \{0,1,\infty \}$ the point $P_i^{0,\vee}$ is mapped to $s\p(i)$.
    \end{itemize}
    The morphism $\mf c^0_{s\p}$ is a neutral $\mb G_m$-gerbe. In particular, $\mc P\p$ is the coarse moduli space\footnote{In the sense of \cite[Definition~4.2.1]{Alper}.} of $\Bun_D^\minnull$. 
\end{Prop}
\begin{proof}
    Recall that we realized $\Bun_D^\minnull$ as the quotient stack $\lb \FL_3^\minnull / \mb G_m^{\times 3} \rb$. The diagonal $\mb G_m$ inside $\mb G_m^{\times 3}$ acts trivially on $\FL_3^\minnull$. Hence, we obtain a $\mb G_m$-gerbe
    \begin{align} \label{the gerbe obtained by modding out center}
        \lb \FL_3^\minnull / \mb G_m^{\times 3} \rb \to \lb \FL_3^\minnull / L \rb,
    \end{align}
    where $L = \mb G_m^{\times 3}/\mb G_m$ \cite[Lemma~\href{https://stacks.math.columbia.edu/tag/06PE}{06PE}]{StacksProject}. It is neutral, since $\mb G_m^{\times 3} \to L$ admits a splitting. The smooth surjective morphism $d^0 \vert : \FL_3^\minnull \to \mb P^1$ induces a smooth surjective morphism $\delta^0 : \lb \FL_3^\minnull / L \rb \to \mb P^1$. By the defining property of $\Bun_D^\minnull$, points in $\lb \FL_3^\minnull / L \rb$ have trivial automorphism groups, hence $\lb \FL_3^\minnull / L \rb$ is an algebraic space \cite[Proposition~\href{https://stacks.math.columbia.edu/tag/04SZ}{04SZ}]{StacksProject}. Further, since $\dim(\FL_3) = 3$ and $\dim(L) = 2$, it has dimension $1$. We see that $\delta^0$ is smooth of relative dimension $0$, hence ètale. The classification theorem \ref{Classification Theorem for degree zero} and the fact that (\ref{the gerbe obtained by modding out center}) is bijective on isomorphism classes of $k$-points show that $\delta^0$ is bijective on $k$-points over $\mb P^1 \setminus \{0,1,\infty\}$ and that the fiber of $\delta^0$ over each of the point $p \in \{0,1,\infty\}$ is isomorphic to $\Spec(k) \sqcup \Spec(k)$. Thus, choosing for each $p \in \{0,1,\infty\}$ any point $y_p$ of the two points in the fiber $\l\delta^0\r^{-1}(p)$, the restriction of $\delta^0$ to the open algebraic subspace $U := \lb \FL_3^\minnull / L \rb \setminus \{y_0,y_1,y_\infty\}$ of $\lb \FL_3^\minnull / L \rb$ is étale, surjective and bijective on $k$-points, hence an isomorphism. Now the claim is clear.
\end{proof}

\subsection{The degree $1$ locus: Charts} \label{Degree 1 locus charts}
The image of the smooth morphism $\theta^0 : \FL_3^\internull \to \Bun_D^0$ constructed in Proposition \ref{Proposition theta zero} does not contain the point $P^{0,\spc}_{\frac{t-1}{t}}$. To get a chart around this point it is easier to use the degree $1$ component (recall from Theorem \ref{Proposition interesting locus and hecke operators} that the isomorphism $\Bun_D^0 \simeq \Bun_D^1$ given via the elementary Hecke operator at $\infty$ induces a homeomorphism $\vert \Bun_D^0 \vert^\inter \simeq \vert \Bun_D^1 \vert^\inter$). So we will be interested in the neighbourhood of the pendant of $P_{\frac{t-1}{t}}^{0,\spc}$ in the degree $1$ component, which is $P_{\frac{t-1}{t}}^{1,\spc}$. 

Let $\Pi^{(1)}$ be the domain of definition of the rational map
\begin{align*}
    r\p : \Pi &\dashrightarrow (\mb P^1)^{\times 4} \\
    \l \begin{bmatrix}
        a_1^{(1)} \\
        a_2^{(1)} \\
        a_3^{(1)}
    \end{bmatrix}, \begin{bmatrix}
        a_1^{(2)} \\
        a_2^{(2)} \\
        a_3^{(2)}
    \end{bmatrix},  \begin{bmatrix}
        a_1^{(3)} \\
        a_2^{(3)} \\
        a_3^{(3)}
    \end{bmatrix},  \lb \begin{pmatrix}
        a_1 \\
        a_2 \\
        a_3
    \end{pmatrix}, \begin{pmatrix}
        b_1 \\
        b_2 \\
        b_3
    \end{pmatrix} \rb \r &\mapsto  \l \begin{bmatrix}
        a_2^{(1)} \\
        a_3^{(1)}
    \end{bmatrix}, \begin{bmatrix}
        a_2^{(2)} \\
        a_3^{(2)}
    \end{bmatrix},  \begin{bmatrix}
        a_2^{(3)} \\
        a_3^{(3)}
    \end{bmatrix}, \begin{bmatrix}
        a_2 \\
        a_3
    \end{bmatrix} \r 
\end{align*}
Let $\Pi^{1} \subset \Pi^{(1)}$ denote the preimage of the complement of all partial diagonals in $(\mb P^1)^{\times 4}$, i.e. of all points $(l_1,l_2,l_3,l_4)$ s.t. $l_i \neq l_j$ for $i \neq j$. We omit the direct calculation that proves the following Lemma.

\begin{Lem}
    The open subscheme $\Pi^1$ of $\Pi$ is $\ul{\Aut}(\mc O(0,0,1))$-invariant and hence corresponds to an open substack $\Bun_D^\onegentc$ of
    \begin{align*}
        \Bun_D^{(0,0,1)} \simeq [\Pi/ \ul{\Aut}(\mc O(0,0,1))].
    \end{align*}
    We have
    \begin{align*}
        \vert \Bun_D^\onegentc \vert \cap \vert \Bun_D^1 \vert^{\inter} = \vert \Bun_D^1 \vert^{\inter} \setminus \left\{P_i^{1,\dec}, P_i^{1,\wedge}, P_i^{1,\vee} \ \vert \ i \in \{0,1,\infty\}\right\},
    \end{align*}
    so in particular $P_{\frac{t-1}{t}}^{1,\spc} \in \vert \Bun_D^\onegentc \vert$.
\end{Lem}

Define
\begin{align*}
    \Upsilon^1 := \l(\mb G_m \setminus \{1\}\r) \times \mb A^1 \times \mb P^1.
\end{align*}
Consider the map
\begin{align*}
        \theta^1 :\Upsilon^1 &\to \Pi^{1} \\
        (u,v,[w_1:w_2]) &\mapsto \l \begin{bmatrix}
                       0 \\ 0 \\ 1\\
                      \end{bmatrix},  \begin{bmatrix}
                        0 \\ 1 \\ 0 \\
                        \end{bmatrix}, \begin{bmatrix}
                         0 \\ 1 \\ 1 \\
                          \end{bmatrix}, \lb \begin{pmatrix}
                              v \\
                              u \\
                              1
                          \end{pmatrix}, \begin{pmatrix}
                              w_1 \\
                              w_2 \\
                              0
                          \end{pmatrix} \rb \r.    
\end{align*}
Let $B \subset \GL_2$ be the standard upper Borel subgroup. Let $H$ be the semidirect product $\mb G_m \ltimes \mb G_a$ for the multiplication action of $\mb G_m$ on $\mb G_a$. We have a quotient map
\begin{align*}
    q : B &\to H \\
    \begin{pmatrix}
        a & c \\
        & b
    \end{pmatrix} &\mapsto \l \frac{a}{b}, \frac{c}{b} \r.
\end{align*}
Let $H$ act on $\Upsilon^1$ via
\begin{align*}
    (a,c). (u,v,[w_1:w_2]) = \l u,av+c\l u - \frac{t-1}{t}\r, [aw_1 + cw_2 : w_2]\r
\end{align*}
and let $B$ act on $\Upsilon^1$ through $q$.

\begin{Lem} \label{Lemma theta one}
    $\theta^1$ induces an isomorphism 
    \begin{align*}
         \lb \Upsilon^1/B \rb \simeq \Bun_D^\onegentc.
    \end{align*}
\end{Lem}
\begin{proof}
    Use Lemma \ref{Lemma about uniformaization of quotient stacks}.
\end{proof}

\subsection{Degree $1$ locus: Map to $\mb P^1$}
In this section we construct the pendant of the map $d^0 : \FL_3^\internull \to \mb P^1$ for the degree $1$ locus. 

Define
\begin{align*}
    d^1 : \Upsilon^1 &\to \mb P^1 \\
    (u,v,[w_1:w_2]) &\mapsto u.
\end{align*}
$d^1$ is invariant for the $B$-action on $\Upsilon^1$ and hence induces a morphism
    \begin{align*}
        \mathfrak d^1 : \Bun_D^\onegentc \underset{\ref{Lemma theta one}}{\simeq} [\Upsilon^1/B] \to \mb P^1.
    \end{align*}

\begin{Prop} \label{Long Proposition about intereseting locus in degree one compoenent}
     The preimage of $\vert \Bun_D^1 \vert^\inter \cap \vert \Bun_D^\onegentc \vert$ under $\Upsilon^1 \to  \lb \Upsilon^1/B \rb \simeq \Bun_D^\onegentc$ is the constructible, but not locally closed subset
    \begin{align*}
        \Upsilon^\interone := V^\inter \sqcup Z^\inter,
    \end{align*}
    of $\Upsilon^1$, where
    \begin{align*}
         V^\inter := \left\{u-\frac{t-1}{t} \neq 0, \ w_2 \neq 0, \ \l u - \frac{t-1}{t} \r w_1 - v w_2 \neq 0\right\}
    \end{align*}
    is open in $\Upsilon^1$ and
    \begin{align*}
        Z^\inter := \left\{u-\frac{t-1}{t} = 0, \ v \neq 0\right\}
    \end{align*}
    is open inside the fiber $\l d^1 \r^{-1}\l \frac{t-1}{t} \r$, hence locally closed inside $\Upsilon^1$. \\
    We have 
    \begin{align*}
        V^\inter &= \Upsilon^\interone \cap \l d^1 \r^{-1}\l\mb P^1 \setminus \{0,1,\infty,\frac{t-1}{t}\} \r, \\
        Z^\inter &= \Upsilon^\interone \cap  \l d^1 \r^{-1}\l \frac{t-1}{t} \r,
    \end{align*}
    in particular the sets $V^\inter$ and $Z^\inter$ are both $H$-invariant. The morphism $d^1$ induces an isomorphism
    \begin{align*}
        [V^\inter/H] \simeq \mb P^1 \setminus \{0,1,\infty,\frac{t-1}{t}\}.
    \end{align*}
    Any point $u\p \in V^\inter(k)$ over $u \in \mb P^1(k) \setminus \{0,1,\infty,\frac{t-1}{t}\}$ maps under $V^\inter(k) \to \Bun_D^1(k)$ to $P_u^1$. \\    
    The action of $H$ on $Z^\inter$ has two orbits, namely
    \begin{align*}
        Z^{\inter,\gen} &:= \{w_2 \neq 0\} \cap Z^\inter, \\
        Z^{\inter,\spc} &:= \{w_2 = 0\} \cap Z^\inter.
    \end{align*}
    Any point in $Z^{\inter,\gen}(k)$ maps under $Z^{\inter,\gen}(k) \to \Bun_D^1(k)$ to $P_{\frac{t-1}{t}}^1$ and any point in $Z^{\inter,\spc}(k)$ maps under $Z^{\inter,\spc}(k) \to \Bun_D^1(k)$ to $P_{\frac{t-1}{t}}^{1,\spc}$. The set $Z^{\inter,\spc}$ is closed inside $\Upsilon^\interone$.
\end{Prop}
\begin{proof}
    All assertions follow from the classification theorem \ref{Classification Theorem for degree one}.
\end{proof}

\begin{Rem}
One may think of $\Upsilon^\interone$ as constructed as follows: Start with the three-dimensional scheme
\begin{align*}
    \Upsilon^1 = \l(\mb G_m \setminus \{1\}\r) \times \mb A^1 \times \mb P^1
\end{align*}
with coordinates $(u,v,[w_1:w_2])$ and think of it as a scheme over $\mb G_m \setminus \{1\}$ via $d^1$ (i.e. $(u,v,[w_1:w_2]) \mapsto u$). Then consider the two hypersurfaces $H_1 = \{w_2 = 0\}$, $H_2 = \left\{\l u - \frac{t-1}{t} \r w_1 - v w_2 = 0\right\}$. They intersect only over $u = \frac{t-1}{t}$. For each $u \in \mb P^1 \setminus \{0,1,\infty, \frac{t-1}{t}\}$ and $i \in \{1,2\}$, $H_i \cap \l d^1 \r^{-1}(u) \simeq \mb A^1$ is a single $H$-orbit and each point in this orbit has stabilizer $\mb G_m$. \\
On the other hand, $H_2 \cap \l d^1 \r^{-1}\l \frac{t-1}{t} \r$ is the locus $\{vw_2 = 0\}$ inside $\mb A^1 \times \mb P^1$. It has two irreducible components: $\{v = 0\} \simeq \mb P^1$ and $\{w_2 = 0\} \simeq \mb A^1$. They intersect nontrivially in the point $\{v = w_2 = 0\}$ and the component $\{w_2 = 0\}$ agrees with $H_1 \cap \l d^1 \r^{-1}\l \frac{t-1}{t} \r$. The locus $\{vw_2 = 0\}$ inside $\l d^1 \r^{-1}\l \frac{t-1}{t} \r$ decomposes into three $H$-orbits: $\{v = 0, w_2 \neq 0\}$, $\{v \neq 0, w_2 = 0\}$ and $\{v=w_2 = 0\}$.\\
Now $\Upsilon^\interone$ arises by cutting out $H_1 \cup H_2$ from $\Upsilon^1$ and then adding the orbit $\{v \neq 0, w_2 = 0\}$ inside the fiber $\l d^1 \r^{-1}\l \frac{t-1}{t} \r$.
\end{Rem}

\begin{Prop} \label{Proposition agreement of morphisms}
    The morphisms $\mathfrak d^0 \circ \mb T^{(l)}_4$ and $\mathfrak d^1$ agree on $\Bun_D^\onegentc \cap \l \mb T^{(l)}_4\r^{-1}(\Bun_D^\internull)$ (intersection of open substacks in $\Bun_D^1$). 
\end{Prop}
\begin{proof}
    It follows from Proposition \ref{Long Proposition about intereseting locus in degree one compoenent} that the preimage of the open substack $\Bun_D^\onegentc \cap \l \mb T^{(l)}_4\r^{-1}(\Bun_D^\internull)$ under the atlas $\Upsilon^1 \to \Bun_D^\onegentc$ is given by 
    \begin{align*}
        \Upsilon^{1,\min} := \Upsilon^1 \setminus \l \left\{w_2 = 0\right\} \cup \left\{\l u - \frac{t-1}{t} \r w_1 - v w_2 = 0\right\} \r.
    \end{align*}
    We compute the composition 
    \begin{align*}
        \Upsilon^{1,\min} \to \Bun_D^\onegentc \xrightarrow{\mb T_4^{(l)}} \Bun_d^0
    \end{align*}
    using for $(u,v,[w_1:w_2]) \in  \Upsilon^{1,\min}(k)$ the diagram 
     \[
\begin{tikzcd}[column sep=huge, ampersand replacement=\&]
\mc O \oplus \mc O(-1) \oplus \mc O(-1)
\arrow[r, "{\left(\begin{smallmatrix}
1 & -wX & -(v-uw)X \\
 & 1 & -u \\
 &  & Y
\end{smallmatrix}\right)}"]
\&
\mc O \oplus \mc O(-1) \oplus \mc O
\arrow[r, "{\left(\begin{smallmatrix}
1 & 0 & 0 \\
0 & Y & 0 \\
0 & 0 & 1
\end{smallmatrix}\right)}"]
\&
\mc O \oplus \mc O \oplus \mc O
\arrow[r, "{\left(\begin{smallmatrix}
Y & wX & vX  \\
0 & 1 & u \\
0 & 0 & 1
\end{smallmatrix}\right)}"]
\&
\mc O(1) \oplus \mc O \oplus \mc O,
\end{tikzcd}
\] \\
where we set $w := w_1/w_2$. We obtain that $(u,v,[w_1:w_2])$ is mapped to
\begin{align*}
    \l\mc O(\ul 0),  \begin{bmatrix}
        0 \\
        -u \\
        1
    \end{bmatrix}, \begin{bmatrix}
        -w \\
        1 \\
        0
    \end{bmatrix},  \begin{bmatrix}
        -t(w(1-u)+v) \\
        1-u \\
        1
    \end{bmatrix},  \lb \begin{pmatrix}
        1 \\
        0 \\
        0
    \end{pmatrix}, \begin{pmatrix}
        0 \\
        0 \\
        1
    \end{pmatrix} \rb \r.
\end{align*}
Applying the automorphism
\begin{align*}
    \begin{pmatrix}
        1 & w & uw \\
        & 1 & u \\
        & & 1
    \end{pmatrix} \in \Aut(\mc O(\ul 0)),
\end{align*}
this point is isomorphic to
\begin{align*}
     \l\mc O(\ul 0),  \begin{bmatrix}
        0 \\
        0 \\
        1
    \end{bmatrix}, \begin{bmatrix}
        0 \\
        1 \\
        0
    \end{bmatrix},  \begin{bmatrix}
        \beta \\
        1 \\
        1
    \end{bmatrix},  \lb \begin{pmatrix}
        1 \\
        0 \\
        0
    \end{pmatrix}, \begin{pmatrix}
        uw \\
        u \\
        1
    \end{pmatrix} \rb \r,
\end{align*}
where we put
\begin{align*}
    \beta := -t(w(1-u)+v) + w = t\l \l u - \frac{t-1}{t}\r w - v \r
\end{align*}
Note that $\beta$ does not vanish on $\Upsilon^{1,\min}$. Hence, we may apply the automorphism
\begin{align*}
    \begin{pmatrix}
        1 & & \\
        1 & -\beta & \\
        1 & & -\beta
    \end{pmatrix} \in \Aut(\mc O(\ul 0))
\end{align*}
to see that the above is isomorphic to
\begin{align*}
     \l\mc O(\ul 0),  \begin{bmatrix}
        0 \\
        0 \\
        1
    \end{bmatrix}, \begin{bmatrix}
        0 \\
        1 \\
        0
    \end{bmatrix},  \begin{bmatrix}
        1 \\
        0 \\
        0
    \end{bmatrix},  \lb \begin{pmatrix}
        1 \\
        1 \\
        1
    \end{pmatrix}, \begin{pmatrix}
        uw \\
        uw - u\beta \\
        uw - \beta
    \end{pmatrix} \rb \r &= \l\mc O(\ul 0),  \begin{bmatrix}
        0 \\
        0 \\
        1
    \end{bmatrix}, \begin{bmatrix}
        0 \\
        1 \\
        0
    \end{bmatrix},  \begin{bmatrix}
        1 \\
        0 \\
        0
    \end{bmatrix},  \lb \begin{pmatrix}
        1 \\
        1 \\
        1
    \end{pmatrix}, \begin{pmatrix}
        \beta \\
        \beta - u\beta \\
        0
    \end{pmatrix} \rb \r \\
        &= \l\mc O(\ul 0),  \begin{bmatrix}
        0 \\
        0 \\
        1
    \end{bmatrix}, \begin{bmatrix}
        0 \\
        1 \\
        0
    \end{bmatrix},  \begin{bmatrix}
        1 \\
        0 \\
        0
    \end{bmatrix},  \lb \begin{pmatrix}
        1 \\
        1 \\
        1
    \end{pmatrix}, \begin{pmatrix}
        1 \\
        1 - u \\
        0
    \end{pmatrix} \rb \r.
\end{align*}
This maps under $\mf d^0$ to $[(1-u)-1:-1] = [u:1] = d^1(u,v,[w_1:w_2])$. Note that the calculation on $k$-points is enough, since $\Upsilon^{1,\min}$ is reduced and $\mb P^1$ is separated.
\end{proof}

\subsection{Line bundles on the interesting locus} 
We would like to study the Picard group of the interesting locus $\vert \Bun_D^0 \vert^\inter$. Unfortunately, by Proposition \ref{Long Proposition about intereseting locus in degree one compoenent} this locus is not locally closed inside $\vert \Bun_D \vert$, so it does not inherit a natural reduced stack structure from $\Bun_D$. We will "approximate" the degree $1$ component of the interesting locus by an open substack $\Bun_D^\interc$ of $\Bun_D^1$, which contains $\vert \Bun_D^1 \vert^\inter$ as a large subset. We will calculate the Picard group of the rigidification $\Bun_D^\interc \rigid \mb G_m$ for the central embedding of $\mb G_m$ into the automorphism groups of points in $\Bun_D^\interc$.

\begin{Def} \label{Definition Buninterc}
    We define:
    \begin{align*}
        \Bun_D^\interonec &:= \l \mb T_4^{(l)} \r^{-1}\l \Bun_D^\internull \r \cup \Bun_D^\onegentc \\
        \Bun_D^\internullec &:= \mb T_4^{(l)}\l \Bun_D^\interonec \r.
    \end{align*}
\end{Def}

The stack $\Bun_D^{i,\interc}$ admits a natural inclusion into $\Bun_D^i$, $i = 1,2$, along which the morphism $\mb T_4^{(l)}$ restricts to an isomorphism
\begin{align*}
    \Bun_D^\interonec \simeq \Bun_D^\internullec.
\end{align*}

\begin{Lem} \label{Lemma about gluing of morphisms}
    The morphisms $\mathfrak d^0 \circ \mb T^{(l)}_4$ and $\mathfrak d^1$ glue to a morphism
\begin{align*}
    \mf d^\interonec : \Bun_D^\interonec\to \mb P^1.
\end{align*}
\end{Lem}
\begin{proof}
    This follows from Proposition \ref{Proposition agreement of morphisms}.
\end{proof}

We now study the rigidifications of the two stacks in the union in \ref{Definition Buninterc} individually, starting with $\Bun_D^\onegentc$. Recall from Lemma \ref{Lemma theta one} the isomorphism
\begin{align*}
    \Bun_D^\onegentc \rigid \mb G_m \simeq \lb \Upsilon^1 / H \rb,
\end{align*}
where 
\begin{align*}
    \Upsilon^1 &= \l(\mb G_m \setminus \{1\}\r) \times \mb A^1 \times \mb P^1 \\
    H &= \mb G_m \ltimes \mb G_a
\end{align*}
and the action of $H$ on $\Upsilon^1$ is given by
\begin{align*}
    (a,c). (u,v,[w_1:w_2]) = \l u,av+c\l u - \frac{t-1}{t}\r, [aw_1 + cw_2 : w_2]\r.
\end{align*}

\begin{Prop} \label{Prop about Picard group near the special point}
    There is an isomorphism
    \begin{align*}
        \Pic\l\Bun_D^\onegentc \rigid \mb G_m\r \simeq \mb Z^{\oplus 2}.
    \end{align*}
\end{Prop}
\begin{proof}
    Since $\Upsilon^1$ is reduced, we may use the exact sequence 
    \begin{align*}
        0 \to \mc O^{\times}\l\Upsilon^1\r^H \to \mc O^{\times}\l\Upsilon^1\r \to \ul{X^*(H)}\l\Upsilon^1\r \to \Pic\l \lb \Upsilon^1 / H \rb \r \to \Pic\l\Upsilon^1\r \to \Pic\l H \times \Upsilon^1\r/p_2^*\Pic\l  \Upsilon^1 \r
    \end{align*}
    \cite[Proposition~2.10]{BrionLinearization}, where $p_2 : H \times \Upsilon^1 \to \Upsilon^1$ denotes the second projection and $\ul{X^*(H)}$ the constant sheaf on the big étale site for the abelian group $X^*(H) = \Hom_{\text{group}}(H,\mb G_m)$ (here we used the fact that $H$ is split, see \cite[Lemma~2.3]{BrionLinearization}). Noting that $H = \mb G_m \times \mb A^1$ as schemes, we use homotopy invariance of the Picard group \cite[Proposition~27.69]{GW2} to deduce
    \begin{align*}
        \Pic\l H \times \Upsilon^1\r/p_2^*\Pic\l  \Upsilon^1 \r &= 0. \\
        \Pic\l \Upsilon^1 \r &\simeq \Pic(\mb P^1) \simeq \mb Z
    \end{align*}
    Further, since $\Upsilon^1$ is connected, we get $\ul{X^*(H)}\l\Upsilon^1\r = X^*(H) = \mb Z$. Finally, we have
    \begin{align*}
        \mc O^{\times}\l\Upsilon^1\r^H \simeq k^{\times}u^{\mb Z}(u-1)^{\mb Z},
    \end{align*}
    where $u,u-1 \in \mc O^{\times}(\mb G_m \setminus \{1\})$. Since the action of $H$ leaves $u$ and $u-1$ invariant, the first map in the above exact sequence is an isomorphism and we obtain an exact sequence
     \begin{align*}
        0 \to \mb Z \to \Pic\l \lb \Upsilon^1 / H \rb \r \to \Pic\l\Upsilon^1\r \to 0,
    \end{align*}
    which shows the claim, since $\mb Z$ is a projective abelian group.
\end{proof}

\begin{Rem}
    \begin{itemize}
        \item[1)] We could have taken any finite closed subset $T \subset \mb G_m$ and replaced $\mb G_m \setminus \{1\}$ in the definition of $\Upsilon^1$ by $\mb G_m \setminus T$ without changing the calculation in \ref{Prop about Picard group near the special point}.
        \item[2)] We could also replace $\Upsilon^1$ by the $H$-invariant open subscheme
        \begin{align*}
            \Upsilon^{1,\new} := \Upsilon^1 \setminus \left\{u-\frac{t-1}{t} = v = 0 \right\}
        \end{align*}
        without changing the Picard group, since $\Upsilon^1$ is smooth and $\Upsilon^1 \setminus \Upsilon^{1,\new}$ has codimension $2$ in $\Upsilon^1$.
    \end{itemize}
\end{Rem}

Next, we analyze $\Bun_D^\internull$. Consider the three open subschemes $U_i$, $i \in \{0,1,\infty\}$, of $\FL_3$ given in coordinates
\begin{align*}
    \lb \begin{pmatrix}
        a_1 \\ a_2 \\ a_3 
    \end{pmatrix}, \begin{pmatrix}
        b_1 \\ b_2 \\ b_3
    \end{pmatrix} \rb
\end{align*}
by
\begin{align*}
    U_1 &= \{a_2 \neq 0,\ a_3 \neq 0,\ a_2b_1 - a_1b_2 \neq 0, \ a_3b_1 - a_1b_3 \neq 0\} \\
    U_\infty &= \{a_1 \neq 0,\ a_3 \neq 0,\ a_1b_2 - a_2b_1 \neq 0, \ a_3b_2 - a_2b_3 \neq 0\} \\
    U_0 &= \{a_1 \neq 0,\ a_2 \neq 0,\ a_1b_3 - a_3b_1 \neq 0, \ a_2b_3 - a_3b_2 \neq 0\}. 
\end{align*}

\begin{Lem} \label{Lemma about structure of interesting degree zero locus}
    For each $i \in \{0,1,\infty\}$ the following hold:
    \begin{itemize}
        \item[i)] $U_i$ is contained in $\FL_3^\internull$
        \item[ii)] We have $U_i = (d^0)^{-1}\l\mb P^1 \setminus \l\{0,1,\infty\} \setminus \{i\}\r \r$. In particular, $\{U_i\}_{i \in \{0,1,\infty\}}$ is an open cover of $\FL_3^\internull$.
    \item[iii)] Write $\mb A^2 = \Spec(k[a,b])$. There is an isomorphism
    \begin{align*}
        \rho_i : \lb U_i/\mb G_m^{\times 3} \rb \simeq [\l\mb A^2 \setminus \{ab = 1\} \r / \mb G_m^{\times 2}],
    \end{align*}
    where $\mb G_m^{\times 2}$ acts on $\mb A^2 \setminus \{ab = 1\}$ via
    \begin{align*}
        (\lambda, \mu).(a,b) \mapsto \l\frac{\mu}{\lambda} a, \frac{\lambda}{\mu} b\r,
    \end{align*}
    which fits in a commutative diagram \\
    \begin{tikzcd}
 \lb U_i/\mb G_m^{\times 3} \rb  \arrow[rr, "\rho_i"] \arrow[rd, "\mf d^0"'] &         & \lb\l\mb A^2 \setminus \{ab = 1\} \r / \mb G_m^{\times 2}\rb \arrow[ld, "\psi_i"] \\
                                                                       & \mb P^1 &                                               
\end{tikzcd} \\
where $\psi_i$ is the composition of the morphism
\begin{align*}
    \lb \l\mb A^2 \setminus \{ab = 1\} \r / \mb G_m^{\times 2} \rb \to \mb A^1
\end{align*}
induced by
\begin{align*}
    (a,b) \mapsto ab,
\end{align*}
the natural inclusion $\mb A^1 \subset \mb P^1$ and the Möbius transformation $M_i : \mb P^1 \to \mb P^1$ given by
\begin{align*}
    M_i = \begin{cases}
        \begin{pmatrix}
            -1 & 1 \\
            0 & 1
        \end{pmatrix}, & i = 1, \\
         \begin{pmatrix}
            0 & 1 \\
            1 & 0
        \end{pmatrix}, & i = \infty, \\
         \begin{pmatrix}
            1 & 0 \\
            0 & 1
        \end{pmatrix}, & i = 0.
    \end{cases}
\end{align*}
\item[iv)] For distinct $i,j \in \{0,1,\infty\}$ one has $U_i \cap U_j = U_0 \cap U_1 \cap U_\infty = \l d^0 \r^{-1}\l \mb P^1 \setminus \{0,1,\infty\} \r$. The restriction of $\mathfrak d^0$ to $\lb U_i/\mb G_m^{\times 3} \rb \cap \lb U_j/\mb G_m^{\times 3} \rb$ is a neutral $\mb G_m$-gerbe over its image.
\end{itemize}
\end{Lem}
\begin{proof}
    The proof is a direct calculation, which we omit.
\end{proof}

\begin{Prop}
    There is an isomorphism
    \begin{align*}
        \Pic(\Bun_D^\internull \rigid \mb G_m) \simeq \mb Z^{\oplus 4}.
    \end{align*}
\end{Prop}
\begin{proof}
   Set $\mf X := \Bun_D^\internull \rigid \mb G_m$, $\mf U_i := \lb U_i/\mb G_m^{\times 3} \rb \rigid \mb G_m$, $\mf U_{ij} := \mf U_i \cap \mf U_j$, $\mf U_{ijl} := \mf U_i \cap \mf U_j \cap \mf U_l$ (intersection inside $\Bun_D^\internull \rigid \mb G_m$). Denote by $\mc U$ the covering $\{\mf U_i\}_{i \in \{0,1,\infty\}}$ of $\mf X$. Let $\mc F$ be the presheaf $U \mapsto \Pic(U) = H^1_\liset(U,\mb G_m)$ on the lisse-étale site of $\Bun_D^\internull \rigid \mb G_m$. If not otherwise specified, cohomology and \v{C}ech-cohomology in the following is taken on the lisse-étale site. The low-term exact sequence obtained from the \v{C}ech to derived functor spectral sequence for the sheaf $\mb G_m$ on the lisse-étale site of $\mf X$ reads: \begin{align*}
       0 \to \check{H}^1(\mc U,\mb G_m) \to \Pic(\mf X) \to \check{H}^0(\mc U,\mc F) \to \check{H}^2(\mc U,\mb G_m)
   \end{align*}
   Using the isomorphism 
    \begin{align*}
    \mf U_{01\infty} = \mf U_{ij} \simeq \mb P^1 \setminus \{0,1,\infty\}
    \end{align*}
    implied by Lemma \ref{Lemma about structure of interesting degree zero locus}
    one sees that for distinct $i,j,l \in \{0,1,\infty\}$ one has
    \begin{align} \label{Pic of higher intersections}
        \Pic\l \mf U_{ij}  \r \simeq \Pic\l \mf U_{ijl}\r \simeq 0.
    \end{align}
    Further, one sees that the \v{C}ech-complex for $\mc U$ and $\mb G_m$ is identified with the \v{C}ech-complex for the sheaf $\mb G_m$ on the lisse-étale site of $\mb P^1$ and the triple cover of $\mb P^1$ consisting of the three opens $\mb P^1 \setminus \{0,1\}$, $\mb P^1 \setminus \{0,\infty\}$, $\mb P^1 \setminus \{1,\infty\}$. Hence,
     \begin{align*}
        \check{H}^1(\mc U,\mb G_m) &\simeq H^1(\mb P^1,\mb G_m) \simeq \mb Z, \\
        \check{H}^2(\mc U,\mb G_m) &\simeq 0.
    \end{align*}   
    
   Consider the action of $\mb G_m$ on $\mb A^2$ given by $\lambda.(a,b) \mapsto (\lambda^{-1}a,\lambda b)$. Then
   \begin{align*}
        \mf U_i \simeq [\l\mb A^2 \setminus \{ab = 1\} \r / \mb G_m],
   \end{align*}
   hence
    \begin{align*}
        \Pic\l \mf U_i \r \simeq X^*(\mb G_m) \simeq \mb Z.
    \end{align*}
    Together with (\ref{Pic of higher intersections}) we obtain
    \begin{align*}
        \check{H}^0(\mc U,\mc F) \simeq \mb Z^{\oplus 3}.
    \end{align*}
    Thus, the above sequence reads
    \begin{align*}
        0 \to \mb Z \to \Pic\l \Bun_D^\internull \rigid \mb G_m \r \to \mb Z^{\oplus 3} \to 0,
    \end{align*}
    which finishes the proof.
\end{proof}

Finally, we use the Mayer-Vietoris sequence to calculate the Picard group of $\Bun_D^\interonec \rigid \mb G_m$.

\begin{Theo} \label{Theorem Picard group of interc}
    There is an isomorphism
    \begin{align*}
        \Pic\l \Bun_D^\interonec \rigid \mb G_m \r \simeq \mb Z^{\oplus 6}.
    \end{align*}
\end{Theo}
\begin{proof}
    We use the cover of $\Bun_D^\interonec \rigid \mb G_m$ by $\l \mb T_4^{(l)} \r^{-1}\l \Bun_D^\internull \r \simeq  \Bun_D^\internull$ and $\Bun_D^\onegentc \rigid \mb G_m$. Their intersection is $\mb G_m \setminus \{1\}$, which has trivial Picard group. Further, we note that the pullback of the elements $u,u-1 \in \Gamma(\mb G_m \setminus \{1\}, \mc O^{\times})$ to $\Gamma(\Upsilon^1, \mc O^{\times})$ are $H$-invariant, hence lie in $\Gamma\l \Bun_D^\onegentc \rigid \mb G_m, \mc O^{\times} \r$. This implies that the restriction map
    \begin{align*}
        \Gamma\l \Bun_D^\onegentc \rigid \mb G_m, \mc O^{\times} \r \to \Gamma\l \mb G_m \setminus \{1\}, \mc O^{\times} \r
    \end{align*}
    and hence also
    \begin{align*}
        \Gamma\l \Bun_D^\onegentc \rigid \mb G_m, \mc O^{\times} \r \times \Gamma\l \Bun_D^\internull \rigid \mb G_m, \mc O^{\times} \r \to \Gamma\l \mb G_m \setminus \{1\}, \mc O^{\times} \r
    \end{align*}
    is surjective. Hence, the relevant part of the Mayer-Vietoris sequence reads
    \begin{align*}
        0 \to  \Pic\l \Bun_D^\interonec \rigid \mb G_m \r \to \mb Z^{\oplus 2} \oplus \mb Z^{\oplus 4} \to 0,
    \end{align*}
    showing the claim.
\end{proof}

It would be interesting to investigate if one can read off the rank of this Picard group only from the numerical datum $D \in \ParDat$.